\input amstex.tex
\input amsppt.sty
\input epsf.sty

\NoBlackBoxes

\define\w{\operatorname{W}}
\define\const{\operatorname{const}}
\define\Y{\Bbb Y}
\define\m{\,\operatorname{mod}\,}
\define\wt{\widetilde}
\define\wh{\widehat}
\define\id{\text{\bf 1}}
\define\Z{\Bbb Z}
\redefine\L{\Cal L}
\define\R{\Bbb R}
\define\C{\Bbb C}
\define\p{\operatorname{Prob}}
\define\X{\frak X}
\define\K{\Cal K}

\redefine\b{\Cal B}
\redefine\P{\Cal P}
\define\di{\operatorname{dilog}}
\define\shift{\operatorname{shift}}

\topmatter
\title Periodic Schur process and cylindric partitions
\endtitle
\author Alexei Borodin\endauthor

%\date December 27, 2005
%\enddate

\abstract Periodic Schur process is a generalization of the Schur
process introduced in \cite{OR1} (math.CO/0107056). We compute its
correlation functions and their bulk scaling limits, and discuss
several applications including asymptotic analysis of uniform
measures on cylindric partitions, time-dependent extensions of the
discrete sine kernel, and bulk limit behavior of certain measures
on partitions introduced in \cite{NO} (hep-th/0306238) in
connection with supersymmetric gauge theories.

\endabstract

\toc \widestnumber\head{8} \head {} Introduction\endhead

\head 1. Periodic Schur process\endhead

\head 2. Correlation functions \endhead

\head 3. Bulk scaling limit\endhead

\head 4. Extensions of the discrete sine kernel\endhead

\head 5. Cylindric partitions\endhead

\head 6. Bulk of large cylindric partitions with finite or slowly growing
period
\endhead

\head 7.  Bulk of large cylindric partitions with period of intermediate
growth\endhead

\head 8. On a measure of Nekrasov and Okounkov
\endhead

\head{} References \endhead

\endtoc

\endtopmatter

\document

\head Introduction
\endhead

One way to see how the content of this paper is different from
previous works on the subject is to examine the following three
pictures.

$$
\epsffile{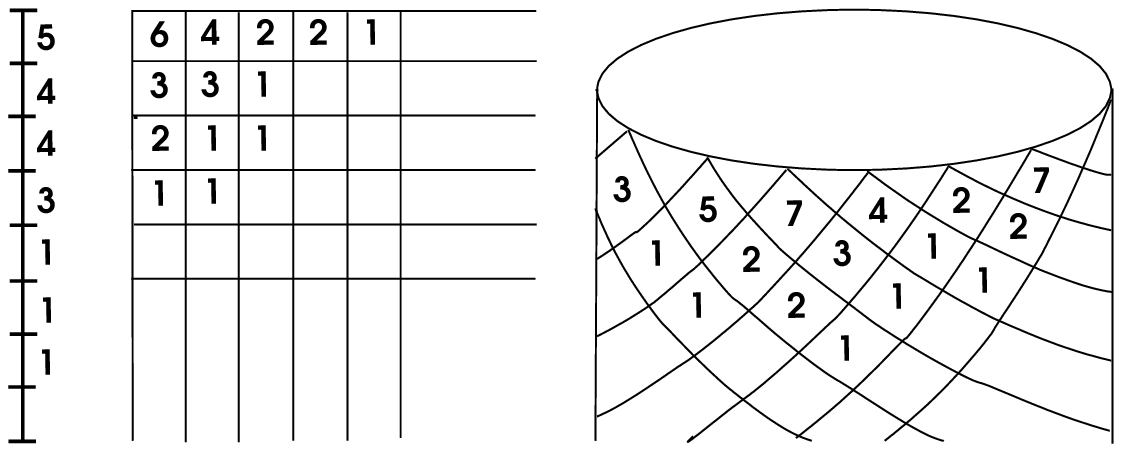}
$$

The leftmost one is a schematic image of an ordinary {\it
partition\/}
--- a way of representing a natural number as an unordered sum of
natural summands.\footnote{In this example the number 19 is
represented as 19=5+4+4+3+1+1+1.} Partitions can also be viewed as
ways of inscribing weakly decreasing nonnegative integers into
unit intervals filling the half-line so that the total number of
nonzero entries is finite.

The middle picture represents a {\it plane partition\/} --- a way
of filling the boxes of the square grid in the quarter plane with
nonnegative integers so that the numbers do not increase as we
move to infinity in the directions of the axes, and the total
number of nonzero entries is also finite.

Finally, the rightmost picture represents a {\it cylindric
partition\/}
--- a way of filling the boxes of the square grid wrapped around a
half-cylinder with nonnegative integers so that the numbers do not
increase as we move away from the border of the cylinder in either
of the two perpendicular directions of the grid lines. The total
number of nonzero entries is also required to be finite.

In this paper we initiate the study of {\it random\/} cylindric
partitions.

Random (ordinary) partitions or, in other words, various
probability measures on partitions, have been extensively studied
since 1940's. The number of references is so large that we will
not even attempt to list them. An excellent survey of uses of
random partitions is available in \cite{O2}.

Random plane partitions are less common, partly because they are
much harder to study. Substantial progress in understanding the
uniform measure on plane partitions with given norm (=sum of
filling numbers) was achieved only recently, see \cite{CK},
\cite{OR1}, \cite{OR2}. In particular, the authors of \cite{OR1}
introduced new techniques which allowed them to derive
determinantal formulas for the correlation functions of random
plane partitions with weights proportional to $q^{\text{norm}}$,
$0<q<1$.\footnote{For such measures the conditional distribution
of the plane partitions with fixed norm is always uniform (and
independent of $q$). For this reason these measures are often also
called uniform.}

The main object of \cite{OR1} called the {\it Schur process\/} is
a generalization of an earlier concept called the {\it Schur
measure} introduced in \cite{O1} to handle certain measures on
(ordinary) partitions. The range of applications of Schur measures
and Schur processes, apart from uniform measures on plane
partitions, is remarkably broad; examples include harmonic
analysis on the infinite symmetric group \cite{O1, \S2.1.4},
\cite{BO1}, Szeg\"o-type formulas for Toeplitz determinants
\cite{BOk}, relative Gromov-Witten theory of $\C^*$ \cite{OP},
random domino tilings of the Aztec diamond \cite{J2}, discrete and
continuous polynuclear growth processes in one space and one time
dimensions \cite{PS}, \cite{J1}, etc.

In this paper we introduce and study a generalization of the Schur
process which we call the {\it periodic\/} Schur process.  We
derive explicit formulas for the correlation functions of this new
process and use them to compute the correlation functions of the
random cylindric partitions with weights proportional to
$q^{\text{norm}}$. This result allows us to obtain various (bulk)
limits of these correlation functions as $q\to 1$. The limiting
cases differ by how fast the radius of the cylinder grows
comparing to $|\ln q|^{-1}$, and the results depend in a
nontrivial way on the angle between the grid lines and the axis of
the cylinder.

We also present two other applications of the periodic Schur
process.

First, we use it to construct an infinite-dimensional family of
determinantal point processes on $\Z^2$ which extend the
well-known one-dimensional discrete sine process. This family
includes two previously obtained in \cite{OR1} and \cite{BO2}
extensions as special cases. Such abundance of two-dimensional
extensions is rather unexpected: In all previously known examples,
probabilistic models yielded only one extension per model, and for
most one-dimensional determinantal point processes no more than
one extension is known.

The second application is a computation of the correlation
functions and their bulk scaling limits for a measure on
(ordinary) partitions introduced in \cite{NO} in connection with
certain supersymmetric gauge theories.

Let us describe our results in more detail.

The periodic Schur process depends on a natural number $N$ (the
period), a parameter $t$, $|t|<1$, and $2N$ specializations
$a[1],b[1],\dots,a[N],b[N]$ of the algebra $\Lambda$ of symmetric
functions.\footnote{A specialization of $\Lambda$ is an algebra
homomorphism of $\Lambda$ to $\C$.}  The process lives on periodic
sequences of $2N$ partitions
$$
\lambda^{(N)}=\lambda^{(0)}\supset \mu^{(1)}\subset
\lambda^{(1)}\supset\dots \subset\lambda^{(N-1)}\supset
\mu^{(N)}\subset\lambda^{(N)}=\lambda^{(0)},
$$
and it assigns to such a sequence the weight proportional to
$$
t^{|\lambda^{(0)}|}s_{\lambda^{(0)}/\mu^{(1)}}(a[1])\,
s_{\lambda^{(1)}/\mu^{(1)}}(b[1])\cdots
s_{\lambda^{(N-1)}/\mu^{(N)}}(a[N])\,
s_{\lambda^{(N)}/\mu^{(N)}}(b[N]).
$$
Here $s_{\lambda/\mu}$'s are the skew Schur functions. The
proportionality coefficient (which we explicitly compute) is
chosen so that the sum of all weights is equal to 1. The weights
can be viewed either as complex numbers or as formal series in
$\Lambda^{\otimes 2N}[t]$.

When $t=0$ the partition $\lambda^{(0)}=\lambda^{(N)}$ must be
empty in order for the sequence to have a nonzero weight, and the
periodic Schur process turns into the conventional Schur process
of \cite{OR1}. On the other hand, if all the specializations
$a[k], b[k]$ are trivial, the periodic Schur process turns into
the uniform measure on partitions, which assigns to a partition
$\lambda$ the weight proportional to $t^{|\lambda|}$.

Denote by $\Z'$ the set $\Z+\frac 12$. To any sequence of
partitions as above it is convenient to associate a point
configuration (subset) in
$\{1,\dots,N\}\times\Z'={\Z'\sqcup\ldots\sqcup\Z'}$ given by
$$
\bigl\{\lambda_i^{(1)}-i+\tfrac 12\bigr\}_{i\ge
1}\sqcup\ldots\sqcup \bigl\{\lambda_i^{(N)}-i+\tfrac
12\bigr\}_{i\ge 1}\,.
$$
This set determines the sequence
$(\lambda^{(1)},\dots,\lambda^{(N)})$ uniquely.

Correlation functions are defined as probabilities, with respect
to the periodic Schur process,\footnote{One should remember that
these ``probabilities'' do not have to be positive and, moreover,
may be just formal series of symmetric functions.} that this
random point configuration contains a fixed finite set of points:
$$
\rho_n(\tau_1,x_1,\dots,\tau_n,x_n)=
\p\left\{x_j\in\bigl\{\lambda_i^{(\tau_j)}-i+\tfrac
12\bigr\}_{i\ge 1}\mid j=1,\dots,n\right\}.
$$

It turns out that the algebraic structure of the correlation
functions substantially simplifies if one considers a modification
of the periodic Schur process which we call the {\it
shift-mixed\/} periodic Schur process. It can be viewed as the
product measure of the periodic Schur process and a measure on
$\Z$ given by
$$
\p\{S\}=\const\cdot z^St^{\frac {S^2}2},\qquad S\in\Z,\quad
z\in\C,
$$
mapped to the space of point configurations in
$\{1,\dots,N\}\times\Z'$ via
$$
(\lambda,S)\mapsto \bigl\{S+\lambda_i^{(1)}-i+\tfrac
12\bigr\}_{i\ge 1}\sqcup\ldots\sqcup
\bigl\{S+\lambda_i^{(N)}-i+\tfrac 12\bigr\}_{i\ge 1}\,.
$$
In other words, all points of the random point configuration of
the periodic Schur process are shifted by an independent integral
valued random variable $S$ distributed as above. Here $z$ is a new
complex parameter.

The normalization constant of the $S$-distribution is the inverse
of one of the Jacobi theta-functions
$$
\theta_3(z;t)=\sum_{S\in\Z}z^S t^{\frac {S^2}2}=\prod_{n\ge
1}(1-t^n)\prod_{n=\frac 12,\frac 32,\frac 52,\dots}
(1+t^{n}z)(1+t^{n}/z),
$$
see e.g. \cite{Er, 13.19(16)}. We assume that $z\ne -t^{\pm \frac
12}, -t^{\pm \frac 32},\dots$, so that $\theta_3(z;t)\ne 0$ .

The correlation functions of the shift-mixed process are defined
in the same way as those of the initial process, and we denote
these new functions as $\rho_n^{\shift}$. It is not hard to see
that $\rho_n$ is equal to the constant term in $z$ of
$\theta_3(z;t)\,\rho_n^{\shift}$.

In order to state our first result, we need to introduce more
notation. Set $a_m[k]:=\frac 1m\,p_m(a[k])$, $b_m[k]:=\frac
1m\,p_m(b[k])$, where $p_k$'s are the Newton power sums. For any
$\tau=1,\dots,N$ denote
$$
\multline F(\tau,\zeta)=\exp\sum_{n\ge
1}\Biggl(\frac{\zeta^{n}}{1-t^n}\sum_{k=1}^{\tau}b_n[k]+
\frac{(t\zeta )^n}{1-t^n}\sum_{k=\tau+1}^N b_n[k]\\
-\frac{(t/\zeta)^{n}}{1-t^n}\sum_{k=1}^{\tau}a_n[k]
-\frac{(1/\zeta)^{n}}{1-t^n}\sum_{k=\tau+1}^N a_n[k] \Biggr).
\endmultline
$$
\proclaim{Theorem A} The correlation functions of the shift-mixed
periodic Schur process have determinantal form: For any $n\ge 1$
and $(\tau_1,x_1),\dots,(\tau_n,x_n)\in\{1,\dots,N\}\times\Z'$ we
have
$$
\rho_n^{\shift}(\tau_1,x_1;\dots;\tau_n,x_n)=\det[K(\tau_i,x_i;\tau_j,x_j)]_{i,j=1}^n
$$
where the generating series of the correlation kernel
$K(\sigma,x;\tau,y)$ has the form
$$
\sum_{x,y\in\Z'}K(\sigma,x;\tau,y)\,\zeta^x\eta^y=\cases
\dfrac{F(\sigma,\zeta)}{F(\tau,\eta^{-1})}\, \sum\limits_{m\in\Z'}
\dfrac{(\zeta\eta)^m}{1+(zt^{m})^{-1}}\,,& \sigma\le \tau,\\ &\\
-\dfrac{F(\sigma,\zeta)}{F(\tau,\eta^{-1})}\,\sum\limits_{m\in\Z'}
\dfrac{(\zeta\eta)^m}{1+zt^m}\,,&\sigma>\tau.
\endcases
$$
\endproclaim

This statement can be understood in two different ways --- as a
formal identity of series in $\Lambda^{\otimes 2N}[t]$ or as a
numeric equality under suitable convergence conditions.

Ramanujan's summation formula for ${}_1\psi_1$-series shows that
the two series in the formula above are expansions of one and the
same holomorphic function in two disjoint annuli. This function
can be expressed in terms of Jacobi theta-functions, see Remark
2.4 below for details.

For $t=0$ Theorem A yields the determinantal formula for the
correlation functions of the conventional Schur process initially
proved in \cite{OR1}, see also \cite{J1} and \cite{BR} for other
proofs. Note that the argument presented in this paper provides an
independent proof of this result.

Let us point out that Theorem A is not particularly trivial even
in the simplest case of the uniform measure on partitions, which
arises when all specializations $a[k]$ and $b[k]$ are trivial.
Then one readily sees that off-diagonal values of the correlation
kernel vanish, and the statement reduces to the fact that the
shift-mixed version of the uniform measure on partitions is
equivalent to the product of countably many independent Bernoulli
measures, see Corollary 2.6 below. The author was not able to
locate this fact in the literature, although certain formulas
equivalent to it can be found in \cite{O1}.

The correlation functions of the initial periodic Schur process
are {\it not\/} determinantal. Nevertheless, they possess a nice
multivariate integral representation given in Corollary 2.8 below,
apart from the fact that they can be extracted from the formula of
Theorem A by taking the constant term in $z$ as mentioned above.

The proof of Theorem A that we present in this paper is a
verification rather than a derivation of the formula for the
correlation functions. The initial proof involved the formalism of
the Fock space and was similar in spirit to the derivations given
in \cite{O1} and \cite{OR1}. However, we decided to leave it out
of this paper because of its length and certain technical
difficulties in justification of formal manipulations with
operators in Fock spaces. As a matter of fact, our initial
inspiration came from the work \cite{Ts}, where the {\it universal
characters\/} -- analogs of the Schur symmetric functions for
nonpolynomial representations of the unitary groups -- were
represented as matrix elements of certain operators in Fock
spaces. We hope to return to the Fock space formalism in a
subsequent publication.

Our second result describes the ``bulk limit'' of the correlation
functions of the periodic Schur process and its shift-mixed
version as $t\to 1$ and the period $N$, as well as the
specializations $a[k]$ and $b[k]$, remain fixed.

\proclaim{Theorem B} Assume that $z\notin\R_{\le 0}$;
$a_m[k],b_m[k]=O(R^m)$ as $m\to\infty$ for some $0<R<1$ and all
$k=1,\dots,N$; and
$$
A_m:=\sum_{k=1}^N
a_m[k]=\sum_{k=1}^N\overline{b_m[k]}=:\overline{B}_m,\qquad
m=1,2,\dots\,.
$$
Then as $t\to 1$, the  correlation functions of the periodic Schur
process and its shift-mixed version have a limit in the following
sense: Choose $x_1(t),\dots,x_n(t)\in\Z'$ such that $|\ln t|\cdot
x_k(t)\to \gamma$ for all $k=1,\dots,n$ and some $\gamma\in\R$,
and all pairwise distances $x_i-x_j$ are independent of $t$. Then
for any $1\le \tau_1,\dots,\tau_n\le N$
$$
\gathered \lim_{t\to
1}\rho_n^{\shift}(\tau_1,x_1(t);\dots,\tau_n,x_n(t))=\det\bigl[
 \K_{\tau_i,\tau_j}^{(z,\gamma)}(x_i-x_j)\bigr]_{i,j=1}^n,\\
\lim_{t\to 1}\rho_n(\tau_1,x_1(t);\dots,\tau_n,x_n(t))=\det\bigl[
 \K_{\tau_i,\tau_j}^{(1,\gamma)}(x_i-x_j)\bigr]_{i,j=1}^n,
\endgathered
$$
where the limit correlation kernel has the form
$$
\K_{\sigma,\tau}^{(z,\gamma)}(d)=\cases\frac 1{2\pi
i}\oint\limits_{|\zeta=1|} \frac{\exp\bigl({-\sum_{m\ge
1}\sum_{k=\sigma+1}^{\tau}(a_m[k]\zeta^{-m}+b_m[k]\zeta^m)}\bigr)}
{1+z^{-1}\exp\bigl({\gamma-\sum_{m\ge
1}(A_m\zeta^{-m}+B_m\zeta^m)}\bigr)}\,\frac{d\zeta}{\zeta^{d+1}},&\sigma\le
\tau,\\
-\frac 1{2\pi i}\oint\limits_{|\zeta|=1}
\frac{\exp\bigl({\sum_{m\ge
1}\sum_{k=\tau+1}^{\sigma}(a_m[k]\zeta^{-m}+b_m[k]\zeta^m)}\bigr)}
{1+z\exp\bigl({-\gamma+\sum_{m\ge
1}(A_m\zeta^{-m}+B_m\zeta^m)}\bigr)}\,\frac{d\zeta}{\zeta^{d+1}},&\sigma>\tau.
\endcases
$$
\endproclaim

Note that the limit correlation functions are invariant with
respect to simultaneous shifts of the ``space variables'' $x_i$.

The assumption $A_m=\overline{B}_m$ is crucial here; without this
condition being satisfied even the integrals above may become
meaningless because the denominators would be allowed to vanish on
the integration contours.

It is rather unusual that the limit can be computed in such
generality.

Before proceeding to cylindric partitions, let us describe in more
detail the two other applications of Theorems A and B mentioned
earlier.

Recall that a classical theorem proved independently by
Aissen--Edrei--Schoen\-berg--Whitney in 1951 \cite{AESW},
\cite{Ed}, and by Thoma in 1964 \cite{Th}, states that a sequence
$\{c_n\}_{n=0}^\infty$, $c_0=1$, is totally positive\footnote{By
definition, this sequence is totally positive if all minors of the
matrix $[c_{i-j}]_{i,j\ge 0}$ are nonnegative.} if and only if its
generating series has the form
$$
\sum_{n=0}^\infty c_nu^n=e^{\gamma u}\frac{\prod_{i\ge 1}
(1+\beta_iu)}{\prod_{i\ge
1}(1-\alpha_iu)}=:TP_{\alpha,\beta,\gamma}(u)
$$
for certain nonnegative parameters $\{\alpha_i\}$, $\{\beta_i\}$
and $\gamma$ such that $\sum_i(\alpha_i+\beta_i)<\infty$.

\proclaim{Corollary 1} For any doubly infinite sequences of
totally positive parameter sets
$\{(\alpha^{(k)},\beta^{(k)},\gamma^{(k)})\}_{k\in\Z}$ and
$\{(\wt\alpha^{(k)},\wt\beta^{(k)},\wt\gamma^{(k)})\}_{k\in\Z}$
satisfying the additional conditions
$$\alpha^{(k)}_i,\ \wt\alpha^{(k)}_i\le 1,\qquad
\beta^{(k)}_i,\ \wt\beta^{(k)}_i<1,\qquad\  i\ge 1,\quad k\in\Z,
$$
and any $c\in(0,\pi)$, there exists a determinantal point
process\footnote{This time all the probabilities are nonnegative!}
on $\Z\times\Z$ with the correlation kernel
$$
\Cal K(\sigma,x;\tau,y)=\cases\dfrac 1{2\pi
i}\int\limits_{e^{-ic}}^{e^{ic}}
\prod\limits_{k=\sigma+1}^\tau{\left(TP_{\alpha^{(k)},\beta^{(k)},\gamma^{(k)}}(
\zeta)\,TP_{\wt\alpha^{(k)},\wt\beta^{(k)},\wt\gamma^{(k)}}(
\zeta^{-1})\right)^{-1}}\,
\dfrac{d\zeta}{\zeta^{x-y+1}}\\
-\dfrac 1{2\pi i}\int\limits_{e^{ic}}^{e^{-ic}}
\prod\limits_{k=\tau+1}^\sigma
{\left(TP_{\alpha^{(k)},\beta^{(k)},\gamma^{(k)}}(
\zeta)\,TP_{\wt\alpha^{(k)},\wt\beta^{(k)},\wt\gamma^{(k)}}(
\zeta^{-1})\right)}\, \dfrac{d\zeta}{\zeta^{x-y+1}}
\endcases
$$
where the first formula is used for $\sigma\le \tau$, the second
formula is used for $\sigma>\tau$, and both integrals are taken
over positively oriented arches of the unit circle.
\endproclaim

The equal time values of the kernel above are exactly those of the
discrete sine kernel on $\Z$: For any $\tau\in\Z$
$$
\Cal K(\tau,x;\tau,y)=\dfrac 1{2\pi i}\int_{e^{-ic}}^{e^{ic}}
\dfrac{d\zeta}{\zeta^{x-y+1}}=\frac{\sin(c(x-y))}{\pi
(x-y)}\,,\qquad x,y\in\Z.
$$
Thus, the kernels $\Cal K(\sigma,x;\tau,y)$ are extensions of the
discrete sine kernel.

The choice of $\alpha_1^{(k)}\equiv 1$ and all other parameters
being zero brings us to the incomplete beta kernel of \cite{OR1}.
On the other hand, taking $\gamma^{(k)}=\wt\gamma^{(k)}$ with all
other parameters being zero yields the extension of the discrete
sine kernel obtained in \cite{BO2, Theorem 4.2}.

As for the second application, we consider a probability measure
on the set of all partitions given by the formula, see \cite{NO,
\S6.2},
$$
M_{\mu,t}(\lambda)=\prod_{n\ge 1}(1-t^n)^{1-\mu^2}\cdot
t^{|\lambda|}\prod_{\square\in\lambda} \frac
{h(\square)^2-\mu^2}{h(\square)^2}\,,\qquad \lambda\in\Y.
$$

Here $\mu\in i\R$ and $t\in (0,1)$ are the parameters, the product
is taken over all boxes of the Young diagram $\lambda$, and
$h(\square)$ denotes the length of the hook rooted at the box
$\square$.

One remarkable feature of this measure is that it interpolates
between the uniform measure on partitions, which appears at
$\mu=0$, and the (poissonized) Plancherel measure on partitions
(see, e.g., \cite{BOO}), which is obtained from $M_{\mu,t}$ by the
limit transition $t\to 0$, $\mu\to\pm i\infty$, $t|\mu|^2\to
\theta>0$.

The measures $M_{\mu,t}$ may be viewed as special cases of the
periodic Schur process with period $N=1$, and Theorem B leads to
the following statement.

\proclaim{Corollary 2} As $t\to 1$, the correlation functions of
$M_{\mu,t}$ have the following limit: Choose
$x_1(t),\dots,x_n(t)\in\Z'$ as in Theorem B. Then the correlation
functions converge to determinants of the limit correlation kernel
$$
\K^{(\gamma,\mu)}(x,y)=\frac 1{2\pi i}\oint_{|\zeta=1|} \frac{1}
{1+e^{\gamma}(1-\zeta)^{-\mu}(1-\zeta^{-1})^{\mu}}\,
\frac{d\zeta}{\zeta^{x-y+1}}\,,\qquad x,y\in\Z.
$$
\endproclaim

Now let us return to cylindric partitions.

It is convenient for us to represent cylindric partitions as
periodic sequences of ordinary partitions by reading the filling
numbers along the diagonal rays which form $45^\circ$ angle with
 grid lines. For example, the visible part of the cylindric
partition represented by the picture in the beginning of this
introduction gives the sequence of partitions
$$
\dots(3)\supset(1)\subset(5,1)\supset (2)\subset(7,2)\supset
(3,1)\subset(4,1)\supset(1)\subset(2,1)\supset(2)\subset (7)\dots
$$

The condition of filling numbers not increasing along the grid
lines is equivalent to neighboring partitions in such a sequence
having {\it interlacing parts\/}. It is also equivalent to saying
that the Young diagrams of any pair of neighboring partitions are
different by either adding or removing a horizontal strip. Such a
relation between two partitions $(\kappa,\nu)$ is denoted as
$\kappa\succ\nu$ or $\kappa\prec\nu$, depending on which of these
two partitions is larger.

The choice of $\succ$'s and $\prec$'s between neighboring
partitions is exactly the choice of the boundary profile of our
cylindric partition near the cut of the cylinder. We will fix such
a profile by providing two periodic sequences $\{A[k]\}$ and
$\{B[k]\}$ of 0's and 1's such that $A[k]+B[k]\equiv 1$; the
$\succ$'s correspond to $A[k]=1$ and $B[k]=0$, and $\prec$'s
correspond to $A[k]=0$, $B[k]=1$. The ratio $\varkappa$ of the
total number of $\prec$'s in a period over the total number of
$\succ$'s in a period will be called the {\it slope\/} of the
profile. The slope depends only on the angle between the grid
lines and the axis of the cylinder. The case $\varkappa=1$
corresponds to the diagonal rays being parallel to the cylinder
axis.

The following picture represents a cylindric partition with slope
$\varkappa<1$, and the visible part of the boundary profile
corresponds to the sequences
$$
\gathered \{A[k]\}=(\dots,1,0,1,1,1,1,1,0,1,\dots), \\
\{B[k]\}=(\dots,0,1,0,0,0,0,0,1,0,\dots).
\endgathered
$$

$$
\epsffile{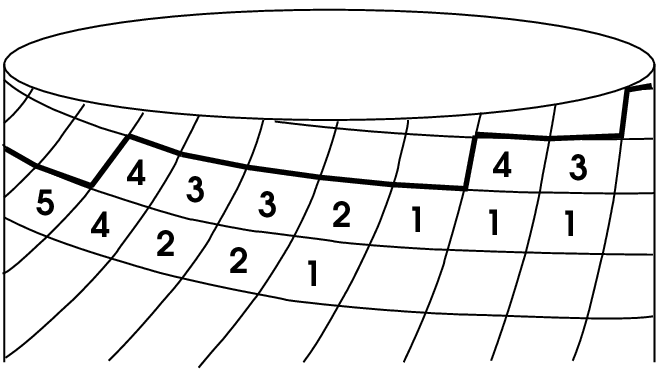}
$$

We will denote by $N$ the period of the sequences above. For any
integer $m$, let $m(N)$ be the smallest positive integer such that
$m\equiv m(N)\mod N$.

\proclaim{Proposition} For any profile $\{A[k]\}_{k=1}^N$,
$\{B[k]\}_{k=1}^N$ as above we have
$$
\sum_{\Sb\text{ cylindric partitions $\pi$} \\ \text{of the given
profile}\endSb }s^{|\pi|}=\prod_{n\ge 1}\frac
1{1-s^{nN}}\prod_{\Sb
p\in \overline{1,N}:\,A[p\,]=1\\
q\in \overline{1,N}:\,B[q\,]=1\endSb}
\frac1{1-s^{(p-q)(N)+(n-1)N}}\,.
$$
Here $|\,\cdot\,|$ denotes the norm {\rm (}=the sum of filling
numbers{\rm )} of cylindric partitions.
\endproclaim

We were unable to find this formula in the literature, although it
may well follow from more refined results of \cite{GK}. In the
limit when the radius of the cylinder becomes large, if
$A=(1,\dots,1,0,\dots,0)$, $B=(0,\dots,0,1,\dots,1)$, and
$\varkappa$ remains bounded away from 0 and $\infty$ (this means
that the boundary of the cylindric partitions locally looks like
the boundary of the quarter plane), the formula of the proposition
reproduces the celebrated formula of MacMahon for the sum of the
weights $s^{\text{norm}}$ over all plane partitions.

The applicability of the periodic Schur process to cylindric
partitions follows from the following basic property of the skew
Schur functions evaluated at a single indeterminate $x$:
$s_{\kappa/\nu}(x)=x^{|\kappa|-|\nu|}$ if $\kappa\succ\nu$ and 0
otherwise. One readily checks that the periodic Schur process with
$t=s^N$ and specializations $a[k]$ and $b[k]$ being the
evaluations at $s^kA[k]$ and $s^{-k}B[k]$, respectively, for all
$k=1,\dots,N$, is exactly the measure $\const\cdot
s^{\text{norm}}$ on cylindric partitions with a fixed profile
described by $\{A[k]\}$ and $\{B[k]\}$. From now on we will use
the term ``uniform measure'' for this distribution.

The above observation implies that Theorem A provides a
determinantal formula for correlation functions of the shift-mixed
modification of the uniform measure on cylindric
partitions.\footnote{By correlation functions of the uniform
measure on cylindric partitions we mean the correlation functions
of the corresponding periodic Schur process, and similarly for the
shift-mixed versions.} Our next goal is to explain what happens to
Theorem B.

In order to state our next result, we need to introduce a curve
$\Gamma_{\varkappa}$ in the complex plane via (here $\varkappa$ is
a positive parameter)
$$
\Gamma_\varkappa=\left\{-\frac{\sin\frac\varphi{1+\varkappa}}{\sin\frac\varphi{1+\varkappa^{-1}}}
\,e^{i\varphi}\mid \varphi\in[-\pi,\pi]\right\}.
$$
This is a piecewise smooth closed curve which has a corner-like
singularity at the point $1$. We orient $\Gamma_{\varkappa}$
counterclockwise.

For $\sigma\le \tau$ set $A(\sigma,\tau]=\sum_{k=\sigma+1}^\tau
A[k]$, and similarly for $B(\sigma,\tau]$. The the slope
$\varkappa$ is equal to $B(0,N]/A(0,N]$.

\proclaim{Theorem C} In the limit $s\to 1$, the correlation
functions of the uniform measure on cylindric partitions with a
given profile $\{A[k],B[k]\}_{k=1}^N$ have a limit in the
following sense: Choose $x_1(s),\dots,x_n(s)\in\Z+\frac 12$ such
that as $s\to 1$, $N|\ln s|\cdot x_k(s)\to \gamma$ for all
$k=1,\dots,n$ and some $\gamma\in\R$, and all pairwise distances
$x_i-x_j$ are independent of $s$. Then for any
$\tau_1,\dots,\tau_n\in\{1,\dots,N\}$
$$
\lim_{s\to 1-}\rho_n(\tau_1,x_1(s);\dots,\tau_n,x_n(s))=\det\bigl[
\K^{(\gamma)}_{\tau_i,\tau_j}(x_i-x_j)\bigr]_{i,j=1}^n
$$
where the correlation kernel has the form
$$
\K^{(\gamma)}_{\sigma,\tau}(d)=\cases\frac 1{2\pi
i}\int\limits_{\Gamma_{\varkappa}}
\dfrac{(1-\zeta)^{B(\sigma,\tau]}(1-\zeta^{-1})^{A(\sigma,\tau]}}
{1+e^\gamma(1-\zeta)^{B(0,N]}(1-\zeta^{-1})^{A(0,N]}}\,\dfrac{d\zeta}{\zeta^{d+1}}\,,&\sigma\le
\tau,\\
-\frac 1{2\pi i}\int\limits_{\Gamma_\varkappa}
\dfrac{(1-\zeta)^{-B(\tau,\sigma]}(1-\zeta^{-1})^{-A(\tau,\sigma]}}
{1+e^{-\gamma}(1-\zeta)^{-B(0,N]}(1-\zeta^{-1})^{-A(0,N]}}\,
\dfrac{d\zeta}{\zeta^{d+1}}\,, &\sigma>\tau.
\endcases
$$
\endproclaim
The function $(1-\zeta)^{B(0,N]}(1-\zeta^{-1})^{A(0,N]}$ takes
nonnegative values on $\Gamma_\varkappa$, and thus the integrals
are correctly defined.

The limit density function $\K^{(\gamma)}_{\tau,\tau}(0)$ does not
depend on $\tau$, which means that it is invariant with respect to
rotations of the cylindric partitions. Note also that it depends
on the profile only through $A(0,N]$ and $B(0,N]$, or, in other
words, through the period $N$ and the slope $\varkappa$.

Interestingly enough, a formal application of Theorem B to the
case of cylindric partitions may produce an incorrect answer if
$\varkappa\ne 1$. This happens because for $\varkappa\ne 1$ the
condition $A_m=\overline{B}_m$ of Theorem B is violated, and we
need to deform the integration contours in order to perform the
asymptotic analysis.

In Theorem C we kept the period $N$ finite while sending $s$ to 1.
The next level of difficulty is to consider the case of periods
growing together with $|\ln s|^{-1}$.

If the growth of the period is slow in the sense that $N|\ln s|$
still tends to 0, then we prove that the limit behavior of the
correlation functions can be read off Theorem C above by taking
the limit $N\to\infty$ and keeping the slope $\varkappa$ fixed. In
the limit one sees extensions of the discrete sine kernel as in
Corollary 1 above with the parameters $\alpha_1^{(k)}$ or
$\wt\alpha_1^{(k)}$ taking values between 0 and 1 and all other
parameters being zero. Details can be found in \S 6.

The case of the period $N$ growing in such a way that the product
$N|\ln s|$ has a finite limit, is substantially more complicated.
The reason is simple -- the limiting behavior depends on the
details of the profile rather than just on its slope. In this
paper we only consider the case of the slope being equal to 1, and
we prove two results.

First, we show that if the profile sequences $\{A[k]\}$ and
$\{B[k]\}$ are periodic with a finite period (in addition to being
periodic with the growing period $N$) then the limit behavior is
just the same as in the case of the slowly growing periods.

Second, we consider the corner-like profiles with $\{A[k]\}$ and
$\{B[k]\}$ consisting of one block of 0's and one block of 1's. We
compute the limit of the correlation functions near two
``corners'', where the partitions are the largest and the
smallest. The results in both cases are governed by the incomplete
beta kernel with the density functions given by certain analytic
expressions involving elliptic functions. Details can be found in
\S7.

To conclude the introduction, let us mention two circles of
questions which we do not discuss in this paper, but which are
certainly very interesting.

The existence of limit correlation functions which decay fast
enough when the distance between the arguments grows, is a strong
indication for the existence of a {\it limit shape\/} of the
corresponding random (ordinary, plane or cylindric) partitions.
Such a decay is present in all the cases we considered. The limit
density function allows one to predict what the limit shape would
look like (see, for example, Comment 3 after Theorem 3.1 below),
but one needs additional arguments to actually prove the
concentration phenomenon. For example, for the measures
$M_{\mu,t}$ considered in Corollary 2 above (and, in fact, for a
substantially larger class of measures on partitions), the
existence of the limit shape as $\mu\to \pm i\infty$ was proved in
\cite{NO} by variational techniques.

Also observe that in this paper we consider only what is usually
called ``the bulk scaling limit'' of the correlation functions. It
would be very interesting to study the ``edge scaling limit'' as
well. In the case of the uniform measure on (skew) plane
partitions different edge scaling limits were computed in
\cite{OR2}.

\subhead Acknowledgements\endsubhead The author is very grateful
to Grigori Olshanski for numerous discussions which were extremely
useful and stimulating. In particular, a large part of \S8 below
consists of results of such discussions. The author would also
like to thank Eric Rains for pointing out a Cauchy type
determinant formula for theta functions due to Frobenius; it is
used in the proof of Corollary 2.8 below.

This research was partially supported by the NSF grant DMS-0402047
and the CRDF grant RIM1-2622-ST-04.

\head 1. Periodic Schur process
\endhead

Fix a natural number $N$ and consider a measure on periodic sequences of $2N$
partitions
$$
\lambda^{(N)}=\lambda^{(0)}\supset \mu^{(1)}\subset \lambda^{(1)}\supset\dots
\subset\lambda^{(N-1)}\supset \mu^{(N)}\subset\lambda^{(N)}=\lambda^{(0)}
$$
by specifying the weight of such a sequence to be equal to
$$
\w(\lambda,\mu)=t^{|\lambda^{(0)}|}s_{\lambda^{(0)}/\mu^{(1)}}(a[1])\,
s_{\lambda^{(1)}/\mu^{(1)}}(b[1])\cdots s_{\lambda^{(N-1)}/\mu^{(N)}}(a[N])\,
s_{\lambda^{(N)}/\mu^{(N)}}(b[N]).
$$
Here $t$ is a parameter and $a[m],b[m]$, $m=1,\dots,N$, are arbitrary
specializations of the algebra $\Lambda$ of symmetric functions.\footnote{A
specialization of $\Lambda$ is an algebra homomorphism of $\Lambda$ to $\C$.}
We will use the notation
$$
a_k[m]:=\tfrac 1k\,p_k(a[m]), \quad b_k[m]:=\tfrac 1k\,p_k(b[m]),
$$
where $p_k$'s are the Newton power sums.

The weights $\w(\lambda,\mu)$ may also be viewed as elements of
$\Lambda^{\otimes 2N}[t]$; in that case the notation $f(a[k])$ or $f(b[k])$ for
$f\in\Lambda$ just indicates to which of the copies of $\Lambda$ in
$\Lambda^{\otimes 2N}$ this symmetric function belongs.

The partition function of such a measure will be denoted as (the notation $\Y$
below stands for the set of all partitions including the empty one)
$$
Z(N,t,a,b):=\sum_{\lambda^{(1)},\mu^{(1)},\dots,\lambda^{(N)},\mu^{(N)} \in\Y}
\w(\lambda,\mu).
$$
We will use the term {\it periodic Schur process\/} for this measure.

It is clear that without loss of generality we may consider only the
distribution of $\lambda$'s. Indeed, by making some of the specializations
$a[k],b[k]$ trivial\footnote{The trivial specialization of $\Lambda$ is
characterized by the fact that all power sums $p_k$, $k\ge 1$, specialize to
zero.}  we can force any given $\mu^{(\,\cdot\,)}$ to coincide with a
neighboring $\lambda^{(\,\cdot\,)}$.

The periodic Schur process as defined above is not symmetric with respect to
the circular shifts $\lambda^{(k)}\mapsto \lambda^{(k+l \m N)}$. This lack of
symmetry can be easily repaired. Introduce new specializations $\wt a[m], \wt
b[m]$ of $\Lambda$ by assigning the following values to the power sums:
$$
\wt a_k[m]:=s^{-km} a_k[m],\quad \wt b_k[m]:=s^{km} b_k[m],
$$
where $s^N=t$. Then we obtain
$$
\multline \w(\lambda,\mu)
=s^{|\lambda^{(1)}|+|\lambda^{(2)}|+\dots+|\lambda^{(N)}|}\\
\times s_{\lambda^{(0)}/\mu^{(1)}}(\wt a[1])\, s_{\lambda^{(1)}/\mu^{(1)}}(\wt
b[1])\cdots s_{\lambda^{(N-1)}/\mu^{(N)}}(\wt a[N])\,
s_{\lambda^{(N)}/\mu^{(N)}}(\wt b[N]),
\endmultline
$$
and this expression already possesses the rotational symmetry. The index $m$ of
$\wt a[m]$ and $\wt b[m]$ can now be viewed as an element of $\Z/N\Z$. We will
use both this form of the measure and the initial non-symmetric one.

When all the specializations are trivial the measure concentrates on the
sequences with coinciding terms: all $\lambda$'s and $\mu$'s become equal. The
resulting distribution on one copy of $\Y$ is the so--called uniform measure:
$\w(\lambda)=t^{|\lambda|}$. This observation shows, in particular, that if one
sets $t=1$ then the partition function may become infinite.

If $t=0$ then $\lambda^{(0)}=\lambda^{(N)}$ must be empty, and the
periodic Schur process coincides with the conventional Schur
process introduced in \cite{OR1}, see also \cite{BR}.

In general, the periodic Schur process may be viewed as the
conventional Schur process on $\Y$ started at the uniform
distribution instead of the empty partition, and conditioned to
yield periodic trajectories.

When the period $N$ is equal to 2, the periodic Schur process can
be viewed as an analog of the Schur measure of \cite{O1} for the
universal characters of the unitary groups (see \cite{Ko} for
details on universal characters):
$$
\sum_{\mu^{(1)},\mu^{(2)}\in\Y}
\w(\lambda^{(1)},\lambda^{(2)};\mu^{(1)},\mu^{(2)})=
S_{[\lambda^{(1)},{\lambda^{(2)}}']}(a[1],b[1]^{-})
\,S_{[\lambda^{(1)},{\lambda^{(2)}}']}(a[2]^{-},b[2]).
$$
Here the superscript ``$-$'' in a specialization $c$ stands for the change of
signs of all the power sums: $p_k(c^-):=-p_k(c)$.

Introduce the notation
$$
A_k:=\sum_{m=1}^N a_k[m],\quad B_k:=\sum_{m=1}^N b_k[m].
$$

\proclaim{Proposition 1.1} The partition function of the periodic Schur process
has the form
$$
\multline Z(N,t,a,b)=\prod_{n\ge 1}\frac 1{1-t^n}\,
\exp\left(\sum\limits_{n=1}^\infty n\Biggl(\sum\limits_{1\le l<k\le
N}a_n[k]b_n[l]+\frac{t^nA_nB_n}{1-t^n}\Biggr)\right)\\
=\prod_{n\ge 1}\frac 1{1-s^{nN}}\, \exp\left(\sum\limits_{n=1}^\infty
\frac{n}{1-s^{nN}}\sum\limits_{k,l=1}^N \Bigl(\id_{k>l}\, s^{(k-l)n}+\id_{k\le
l}\, s^{(N+k-l)n}\Bigr)\wt a_n[k]\wt b_n[l]\right).
\endmultline
$$
\endproclaim

\example{Remark 1.2} The equality above can be viewed either as an identity of
formal series in $\Lambda^{\otimes 2N}[t]$ or as a numeric equality under the
assumption that $|t|<1$ and the series $\sum_{n\ge 1} n a_n[k]b_n[l]$ are
absolutely convergent for any $k,l=1,\dots,N$.
\endexample

\example{Remark 1.3} Here is a different way to write the formula for the
partition function. For two specializations $\alpha$ and $\beta$ of $\Lambda$
set
$$
H(\alpha;\beta)=\sum_{\lambda\in\Y}s_\lambda(\alpha)s_\lambda(\beta).
$$
If $\alpha$ and $\beta$ are evaluations of symmetric functions at
variables $(\alpha_n)$ and $(\beta_n)$ then the Cauchy identity,
see e.g. \cite{Macd, \S{I} (4.3)}, reads $
H(\alpha;\beta)=\prod_{i,j}(1-\alpha_i\beta_j)^{-1}. $ It is also
not hard to verify that $H(\alpha;\beta)=\exp\sum_{n\ge 1}
p_n(\alpha)p_n(\beta)/n$. Hence,
$$
\aligned Z(N,t,a,b)&=\prod_{n\ge 1}\frac 1{1-t^n}\prod\limits_{1\le l<k\le
N}H(a[k];b[l])\,\prod_{n\ge 1}\prod\limits_{k,l=1}^N H(t^n a[k]; b[l])\\
&=\prod_{n\ge 1}\left(\frac 1{1-s^{nN}}\prod_{ k>l} H\bigl(s^{(k-l)n}\wt a[k];
\wt b[l]\bigr)\prod_{ k\le l} H\bigl(s^{(N+k-l)n}\wt a[k]; \wt
b[l]\bigr)\right),
\endaligned
$$
where the multiplication of a specialization $\alpha$ by a scalar $q$ is
defined via $$p_n(q\cdot\alpha):=q^np_n(\alpha), \qquad n\ge 1.$$
\endexample

\demo{Proof of Proposition 1.1} The argument uses the same idea as
\cite{Macd, Ex.~I.5.28}. It is more convenient to work in terms of
$s,\wt a,\wt b$. We have
$$
\multline
Z(N,t,a,b)=\sum_{\lambda,\mu} s^{|\lambda^{(1)}|+|\lambda^{(2)}|+\dots+|\lambda^{(N)}|}\\
\times s_{\lambda^{(N)}/\mu^{(1)}}(\wt a[1])\, s_{\lambda^{(1)}/\mu^{(1)}}(\wt
b[1])\cdots s_{\lambda^{(N-1)}/\mu^{(N)}}(\wt a[N])\,
s_{\lambda^{(N)}/\mu^{(N)}}(\wt b[N])\\
=\sum_{\lambda,\mu} s^{|\mu^{(1)}|+|\mu^{(2)}|+\dots+|\mu^{(N)}|}\qquad\qquad
\qquad\qquad\qquad\qquad
\qquad\qquad \qquad\qquad\qquad\qquad\\
\times s_{\lambda^{(N)}/\mu^{(1)}}(s\wt a[1])\, s_{\lambda^{(1)}/\mu^{(1)}}(\wt
b[1])\cdots s_{\lambda^{(N-1)}/\mu^{(N)}}(s\wt
a[N])\, s_{\lambda^{(N)}/\mu^{(N)}}(\wt b[N])\\
=H(s\wt a[1];\wt b[N])H(s\wt a[2];\wt b[1])\dots H(s\wt a[N];\wt b[N-1])
 \sum_{\kappa,\mu} s^{|\mu^{(1)}|+|\mu^{(2)}|+\dots+|\mu^{(N)}|}\qquad\\
\times s_{\mu^{(1)}/\kappa^{(N)}}(\wt b[N])\, s_{\mu^{(1)}/\kappa^{(1)}}(s \wt
a[2])\cdots s_{\mu^{(N)}/\kappa^{(N-1)}}(\wt b[N-1])\,
s_{\mu^{(N)}/\kappa^{(N)}}( s\wt a[1])
\endmultline
$$
where to sum over $\lambda$'s we used the well-known formula
$$
\sum_{\lambda} s_{\lambda/\mu}(\alpha)s_{\lambda/\widehat\mu}(\beta)=
H(\alpha;\beta) \sum_{\kappa}s_{\mu/\kappa}(\beta)s_{\widehat
\mu/\kappa}(\alpha),
$$
see \cite{Macd, Ex.~I.5.26}. Applying the same trick to sum over
$\mu$'s we get
$$
\multline Z(N,t,a,b)=H(s\wt a[1];\wt b[N])H(s\wt a[2];\wt b[1])\dots
H(s\wt a[N];\wt b[N-1])\\
\qquad\  \times H(s^2\wt a[1];\wt b[N-1])H(s^2\wt a[2]; \wt b[N])\dots
H(s^2\wt a[N];\wt b[N-2])\\
\times
 \sum_{\kappa,\rho}s^{|\kappa^{(1)}|+|\kappa^{(2)}|+\dots+
|\kappa^{(N)}|}\qquad\qquad \qquad\qquad\qquad\qquad
\qquad\qquad \qquad\qquad\qquad\qquad\\
\times s_{\kappa^{(N)}/\rho^{(1)}}(s^2\wt a[2])\, s_{\kappa^{(1)}/\rho^{(1)}}(
\wt b[N])\cdots s_{\kappa^{(N-1)}/\rho^{(N)}}(s^2\wt a[1])\,
s_{\kappa^{(N)}/\rho^{(N)}}(\wt b[N-1]).
\endmultline
$$
Continuing in the same fashion after $N$ summations we obtain
$$
Z(N,t,a,b)=\prod_{ k>l} H\bigl(s^{k-l}\wt a[k]; \wt b[l]\bigr)\prod_{ k\le l}
H\bigl(s^{N+k-l}\wt a[k]; \wt b[l]\bigr)\cdot Z(N,t,s^Na,b).
$$
Iterating the above procedure we arrive at the final formula noting that both
in the formal and analytic settings
$$
\lim_{n\to\infty}Z(N,t,t^na,b)=Z(N,t,\operatorname{trivial},b)=\prod_{n\ge
1}\frac 1{1-t^n}\,.\qed
$$
\enddemo

\example{Remark 1.4} Using the formula
$$
\multline \sum_{\lambda,\mu,\nu} s_{\kappa/\mu}(\alpha) s_{\lambda/\mu}(\beta)
s_{\lambda/\nu}(\gamma)s_{\sigma/\nu}(\delta)=H(\beta;\gamma)\sum_{\mu,\nu,\rho}s_{\kappa/\mu}(\alpha)
s_{\mu/\rho}(\gamma)s_{\nu/\rho}(\beta)s_{\sigma/\nu}(\delta)\\=H(\beta;\gamma)\sum_{\rho}s_{\kappa/\rho}
(\alpha,\gamma) s_{\sigma/\rho}(\beta,\delta)
\endmultline
$$
one readily sees that the projection of the periodic Schur process to
$$
\lambda^{(0)},\lambda^{(k_1)},\dots,\lambda^{(k_{M-1})},
\lambda^{(k_M)}=\lambda^{(N)}
$$
is again a periodic Schur process with a shorter period $M$ and modified
specializations
$$
\gather \wh a_n[1]=\sum_{i=1}^{k_1} a_n[i],\quad \wh
a_n[2]=\sum_{i=k_1+1}^{k_2} a_n[i],\quad \dots,\quad \wh
a_n[M]=\sum_{i=k_{M-1}+1}^{N} a_n[i],\\ \wh b_n[1]=\sum_{i=1}^{k_1}
b_n[i],\quad \wh b_n[2]=\sum_{i=k_1+1}^{k_2} b_n[i],\quad \dots,\quad \wh
b_n[M]=\sum_{i=k_{M-1}+1}^{N} b_n[i].
\endgather
$$
This fact can be used to define periodic Schur processes with {\it continuous
time}.

\endexample

\head 2. Correlation functions
\endhead

For any $\lambda\in\Y$ set
$$
\L(\lambda)=\bigl\{\lambda_i-i+\tfrac 12\bigr\}_{i\ge 1}\subset\Z'=\Z+\tfrac
12\,.
$$
By definition, the $n$th dynamical correlation function of the
periodic Schur process is the probability that the random point
configurations $\L(\lambda^{(\tau_k)})$ contain some fixed points
$x_k\in\Z'$ for all $k=1,\dots,n$:
$$
\rho_n(\tau_1,x_1;\dots;\tau_n,x_n)=\frac 1{Z(N,t,a,b)}
\sum_{(\lambda,\mu)\,:\,x_k\in\L(\lambda^{(\tau_k)}),\,k=\overline{1,n}}
\w(\lambda,\mu).
$$
Some of the time moments $\tau_k$ may coincide, but if
$(\tau_i,x_i)=(\tau_j,x_j)$ for $i\ne j$ then the correlation function
vanishes.

In order to compute the correlation functions it is convenient to introduce a
modification of the periodic Schur process which we call the {\it shift-mixed}
periodic Schur process.

The shift-mixed process is a measure on point configurations in
the disjoint union of $N$ copies of $\Z'$ obtained as follows. Let
us take in the $i$th copy of $\Z'$ the point configuration
$\L(\lambda^{(i)})$, where $(\lambda^{(1)},\dots,\lambda^{(N)})$
form the periodic Schur process, and let us shift all these $N$
point configurations simultaneously by a random integer $S$
distributed according to
$$
\p\{S\}={z^S t^{\frac {S^2}2}}/{\theta_3(z;t)}\,.
$$

Here $z\ne 0$ is a new parameter, and
$\theta_3(z;t)=\sum_{S\in\Z}z^S t^{\frac {S^2}2}$ is the partition
function of the weights $z^S t^{\frac {S^2}2}$ which also happens
to be one of the Jacobi theta-functions, see e.g. \cite{Er,13.19}.

To summarize, in the shift-mixed periodic Schur process the weight of the point
configuration of the form
$$
\bigl\{S+\lambda_i^{(1)}-i+\tfrac 12\bigr\}_{i\ge
1}\sqcup\dots\sqcup\bigl\{S+\lambda_i^{(N)}-i+\tfrac 12\bigr\}_{i\ge 1}
$$
is equal to $z^S t^{\frac {S^2}2}\sum_{\mu}\w(\lambda,\mu)$, all point
configurations not of this form have weight zero, and the partition function is
equal $\theta_3(z;t)Z(N,t,a,b)$.

Since
$$
\theta_3(z;t)=\prod_{n\ge 1}(1-t^n)\cdot\prod_{n=\frac 12,\frac 32,\frac
52,\dots} (1+t^{n}z)(1+t^{n}/z),
$$
see e.g. \cite{Er, 13.19(16)}, we will always assume that $z\ne
-t^{\pm \frac 12}, -t^{\pm \frac 32},\dots$, so that the partition
function is never zero.

The dynamical correlation functions of the shift-mixed process are
defined in the same way as those of the initial process; we will
denote them by $\rho_n^{\shift}$.

\proclaim{Proposition 2.1} The dynamical correlation functions of the periodic
Schur process and its shift-mixed modification are related as follows:
$$
\gathered \rho_n(\tau_1,x_1;\dots;\tau_n,x_n)=\text{\rm constant term in $z$ of
}
\left\{\,\theta_3(z;t)\,\rho_n^{\shift}(\tau_1,x_1;\dots;\tau_n,x_n)\,\right\},\\
\rho_n^{\shift}(\tau_1,x_1;\dots;\tau_n,x_n)=\frac1{\theta_3(z;t)}
\sum_{S\in\Z}z^St^{\frac{S^2}2}\rho_n(\tau_1,x_1-S;\dots;\tau_n,x_n-S).
\endgathered
$$
\endproclaim
\demo{Proof} Follows directly from the definition of the shift-mixed process.
\qed
\enddemo

The reason for introducing the shift-mixed process is the fact that this
process is {\it determinantal\/}, i.e., its correlation functions can be
written as certain determinants. To make an exact statement we need additional
notation.

For any $\tau=1,\dots,N$ set
$$
\align F(\tau,\zeta)&=\exp\sum_{n\ge
1}\Biggl(\frac{\zeta^{n}}{1-t^n}\sum_{k=1}^{\tau}b_n[k]+ \frac{(t\zeta
)^n}{1-t^n}\sum_{k=\tau+1}^N b_n[k]\\& \qquad\qquad\qquad\qquad\qquad\qquad
-\frac{(t/\zeta)^{n}}{1-t^n}\sum_{k=1}^{\tau}a_n[k]
-\frac{(1/\zeta)^{n}}{1-t^n}\sum_{k=\tau+1}^N a_n[k] \Biggr)
\\&=\exp\sum_{n\ge 1}\Biggl( \frac{(s^N\zeta)^{n}}{1-s^{nN}} \sum_{k=\tau+1}^{N+\tau}
\wt b_n[k]s^{-kn}-\frac{(1/\zeta)^{n}}{1-s^{nN}} \sum_{k=\tau+1}^{N+\tau} \wt
a_n[k]s^{kn} \Biggr).
\endalign
$$

\proclaim{Theorem 2.2} The dynamical correlation functions of the shift-mixed
periodic Schur process have determinantal form
$$
\rho_n^{\shift}(\tau_1,x_1;\dots;\tau_n,x_n)=\det[K(\tau_i,x_i;\tau_j,x_j)]_{i,j=1}^n
$$
where the generating series of the correlation kernel $K(\sigma,x;\tau,y)$ has
the form
$$
\sum_{x,y\in\Z'}K(\sigma,x;\tau,y)\,\zeta^x\eta^y=\cases
\dfrac{F(\sigma,\zeta)}{F(\tau,\eta^{-1})}\, \sum\limits_{m\in\Z'}
\dfrac{(\zeta\eta)^m}{1+(zt^{m})^{-1}}\,,& \sigma\le \tau,\\ &\\
-\dfrac{F(\sigma,\zeta)}{F(\tau,\eta^{-1})}\,\sum\limits_{m\in\Z'}
\dfrac{(\zeta\eta)^m}{1+zt^m}\,,&\sigma>\tau.
\endcases
$$
\endproclaim

\example{Remark 2.3} Similarly to Proposition 1.1, the formula above carries
two statements: One holds in the algebra of formal series in $\Lambda^{\otimes
2N}[t]$ with $z$ being an arbitrary nonzero complex number. To decompose the
right-hand side into a series in $t$ one uses the expansions
$$
\frac1{1+zt^m}=\cases\sum_{i\ge 0} (-zt^m)^i,& m>0,\\
-\sum_{i\ge 1}(-zt^m)^{-i},&m<0,
\endcases\qquad \frac 1{1+(zt^m)^{-1}}=1-\frac1{1+zt^m}\,.
$$

The other statement is a numeric equality which we prove under the following
convergence conditions:
$$
\gathered a_n[k],b_n[l]=O(R^n)  \, \text{  for some $0<R<1$ and all  } \, 1\le
k,l\le N;\\ \quad |t|<1; \qquad z\ne 0, -t^{\pm \frac 12}, -t^{\pm \frac 32},
\dots\,.
\endgathered
$$
These conditions guarantee that the generating function above is
an analytic function in $\zeta$ and $\eta$ varying in an annulus
either slightly outside or slightly inside the unit circle
(depending on whether $\sigma\le \tau$ or $\sigma>\tau$), and the
values of the kernel are obtained as the Laurent coefficients of
this function; see also the next remark.
\endexample

\example{Remark 2.4} The extraction of coefficients of the generating series
can be performed by computing the corresponding contour integrals. Also, both
series in the formula for the kernel above sum up to one and the same analytic
function in disjoint domains:
$$
\prod_{n\ge 1}(1-t^n)^3\,\frac{-\sqrt{y}\,\theta_3(yz;t)}{\theta_3(-y t^{-\frac
12};t)\,\theta_3(z;t)}=\cases\sum\limits_{m\in\Z'}
\dfrac{y^m}{1+(zt^{m})^{-1}}\,,&1<|y|<|t|^{-1},\\
-\sum\limits_{m\in\Z'} \dfrac{y^m}{1+zt^m}\,,&|t|<|y|<1. \endcases
$$
These are special cases of Ramanujan's summation formula for
${}_1\psi_1$-series, see e.g. \cite{GR, \S5.2}. Indeed, the second
sum equals
$$
\multline \sum\limits_{m\in\Z}\dfrac{y^{m+\frac 12}}{1+zt^{m+\frac
12}}=\frac{\sqrt{y}}{1+zt^\frac12}\sum_{m\in\Z}\frac{(-zt^{\frac
12};t)_m\,y^{m}}{(-zt^{\frac 32};t)_m}
=\frac{\sqrt{y}}{1+zt^\frac12}{}_1\psi_1(-zt^\frac12;-zt^\frac32;t,y)\\
=\frac{\sqrt{y}}{1+zt^\frac12}\,\frac{(t,t,-t^\frac12zy,-t^\frac12/(zy);t)_\infty}
{(-t^\frac32z,-t^\frac12/z,y,t/y;t)_\infty}=\frac{(t;t)_\infty^3\sqrt{y}\,
\theta_3(yz;t)}{\theta_3(-y t^{-\frac 12};t)\,\theta_3(z;t)}
\endmultline
$$
where we use the conventional notation
$$
(a;t)_m=\prod_{n=0}^{m-1}(1-at^n),\qquad
(a;t)_\infty=\prod_{n=0}^{\infty}(1-at^n),
$$
The first sum is obtained by the change $(y,z)\mapsto (y^{-1},z^{-1})$.

Hence, in the analytic setting the formula
for the correlation kernel above can be rewritten as follows:
$$
K(\sigma,x;\tau,y)=-\frac{\prod_{n\ge 1}(1-t^n)^3}{\theta_3(z;t) \,(2\pi
i)^2}\oint_\zeta\oint_\eta \dfrac{F(\sigma,\zeta)}{F(\tau,\eta^{-1})}
\,\frac{\theta_3(z\zeta\eta;t)}{\theta_3(-\zeta\eta\, t^{-\frac 12};t)}\,
\frac{d\zeta d\eta}{\zeta^{x+\frac12}\eta^{y+\frac12}}\,,
$$
where both integration contours are simple loops going around the origin in
positive direction such that $R<|\zeta|,|\eta|<R^{-1}$, and for $\sigma\le
\tau$ we have $1<|\zeta\eta|<\min\{R^{-1},|t|^{-1}\}$, while for $\sigma>\tau$
we have $\max\{R,|t|\}<|\zeta\eta|<1$.
\endexample

\example{Remark 2.5} In the limit $t\to 0$ both the periodic Schur
process and its shift-mixed version turn into the conventional
Schur process of \cite{OR1}. Accordingly, Theorem 2.2 yields a
determinantal formula for the correlation functions of the Schur
process derived in \cite{OR1} (see also \cite{J1} and \cite{BR}
for different proofs). Note that the summation formulas of Remark
2.4 just turn into geometric series
$$
\frac {\sqrt{y}}{y-1}=\cases\sum\limits_{m\in\Z'_-}
y^m,&|y|>1,\\
-\sum\limits_{m\in\Z'_+} {y^m},& |y|<1.
\endcases
$$
\endexample

\proclaim{Corollary 2.6} The shift-mixed uniform measure on partitions is
isomorphic to the product of independent Bernoulli random variables on
$\{0,1\}^{\Z'}=\{(x_m)_{m\in\Z'}\}$ with
$$
\p\{x_m=1\}=\frac{zt^m}{1+zt^m}\,,\qquad m\in\Z',
$$
in the following sense: Pairs $(\lambda,S)\in\Y\times\Z$ are in one-to-one
correspondence with the sequences $(x_m)\in \{0,1\}^{\Z'}$ of nonzero weights
via
$$
x_m=\cases 1,& \text{if  }\  m\in \{S+\lambda_i-i+\frac 12\},\\
0,&\text{otherwise},\endcases
$$
and the weight of such a sequence is equal to
$$
\p\{\lambda,S\}=\frac{t^{|\lambda|+\frac {S^2}2}
z^S}{\prod_{n=\frac12,\frac32,\dots}(1+t^{n}z)(1+t^{n}/z)}\,.
$$
\endproclaim

\demo{Proof} Follows from the fact that for the trivial specializations $a$ and
$b$ the periodic Schur process turns into the uniform measure (see \S1), the
function $F(\,\cdot\,,\,\cdot\,)$ becomes identically equal to 1, and the equal
time values of the correlation kernel are readily seen to be equal to
$$
K(\tau,x;\tau,y)=\frac{zt^x}{1+zt^{x}}\,\delta_{x,y}\,.\qed
$$
\enddemo

\example{Remark 2.7} Corollary 2.6 is equivalent to \cite{O1,
(3.14)}. It is also fairly easy to prove this statement
independently by explicitly computing the weight of a sequence in
$\{0,1\}^{\Z'}$. This will essentially be done in the proof of
Proposition 2.12 below.
\endexample

Before proceeding to the proof of Theorem 2.2 let us draw one more
corollary. Proposition 2.1 explains that one can obtain the
correlation functions of the periodic Schur process by extracting
the constant coefficient in $z$ from the determinantal formula of
Theorem 2.2. In fact, this extraction can be performed explicitly
yielding a multivariate integral representation. To state the
result it is more convenient to work with another Jacobi theta
function $\theta_1(x;t)$ defined as follows, cf. \cite{Er, 13.19}:
$$
\theta_1 (x;t)=\sum_{n=-\infty}^{\infty} (-1)^n\, t^{\frac{n(n+1)}2} x^{n+\frac
12}=(t,t)_\infty (x^{\frac 12}-x^{-\frac 12})\, (tx;t)_\infty (t/x;t)_\infty.
$$
Since, as was mentioned above, $\theta_3(x;t)=(t;t)_\infty
(-\sqrt{t}\,x;t)_\infty (-\sqrt{t}/x;t)_\infty$, we have
$$
\theta_1(x;t)=-\tfrac 1{\sqrt{x}}\,\theta_3\left(-\tfrac x{\sqrt{t}};t\right).
$$

\proclaim{Corollary 2.8} The dynamical correlation functions of the periodic
Schur process in the analytic setting can be written as $$ \multline
\rho_n(\tau_1,x_1;\dots;\tau_n.x_n)\\ =\frac { (t;t)_\infty^{3n}}{(2\pi
i)^{2n}} \oint\cdots\oint\frac{\prod_{1\le i<j\le
n}\theta_1(\zeta_i/\zeta_j;t)\,\theta_1(\eta_i/\eta_j;t)}
{\prod_{i,j=1}^n\theta_1(\zeta_i\eta_j;t)}\,\prod_{i=1}^n
\frac{F(\tau_i,\zeta_i)}{F(\tau_i,\eta_i^{-1})}\,
\frac{d\zeta_id\eta_i}{(\zeta_i\eta_i)^{x_i+1}}\,,
\endmultline
$$
where the integration variables $\zeta_i$ and $\eta_i$ range over circles
$|\zeta_i|=\alpha_i$, $|\eta_i|=\beta_i$ such that
$$
\alpha_1>\frac 1{\beta_1}>\alpha_2>\frac 1{\beta_2}>\dots>\alpha_n>\frac
1{\beta_n}
$$
and all the radii $\alpha_i$ and $\beta_i$ are close enough to 1.
\endproclaim
\demo{Proof} We start with the interpretation of Theorem 2.2 given in Remark
2.4 and note that
$$
-\frac{\theta_3(z\zeta\eta;t)}{\theta_3(z;t)\,\theta_3(-\zeta\eta t^{-\frac
12};t)}=\frac{\theta_1(-z\sqrt{t}\,\zeta\eta;t)}{\theta_1(-z\sqrt{t};t)\,
\theta_1(\zeta\eta;t)}=\frac{\theta_1(\wh z\,\zeta\eta;t)}{\theta_1(\wh z;t)\,
\theta_1(\zeta\eta;t)}\,,
$$
where we used the notation $\wh z=-z\sqrt{t}$. The following
Cauchy-type determinantal formula is due to Frobenius, see
\cite{R, Lemma 4.3}, \cite{F}:
$$
\det\left[\frac{\theta_1(\wh z\,\zeta_i\eta_j;t)}{\theta_1(\wh z;t)\,
\theta_1(\zeta_i\eta_j;t)}\right]_{i,j=1}^n=\frac{\theta_1(\wh z\prod_{i=1}^n
\zeta_i\eta_i;t)}{\theta_1(\wh z;t)}\, \frac{\prod_{1\le i<j\le
n}\theta_1(\zeta_i/\zeta_j;t)\,\theta_1(\eta_i/\eta_j;t)}
{\prod_{i,j=1}^n\theta_1(\zeta_i\eta_j;t)}\,.
$$
Since
$$
\frac{\theta_1(\wh z\prod_{i=1}^n \zeta_i\eta_i;t)}{\theta_1(\wh
z;t)}=\frac{\theta_3(z\prod_{i=1}^n
\zeta_i\eta_i;t)}{\theta_3(z;t)\sqrt{\zeta_1\eta_1\dots\zeta_n\eta_n}\,},
$$
and the constant term in $z$ of $\theta_3(z\prod_{i=1}^n \zeta_i\eta_i;t)$ is
equal to 1, the application of Proposition 2.1 concludes the proof.\qed
\enddemo

Note that Corollary 2.8 implies the generating series of the density function
of the uniform measure on partitions is, up to a constant, the inverse of the
first Jacobi theta-function:
$$
\sum_{x\in\Z'} \rho_1(x)\,\xi^x=(t;t)_\infty/\theta_1(\xi;t).
$$
\demo{Proof of Theorem 2.2} We will provide a proof for the algebraic variant
of the theorem when both sides are considered as formal series, see Remark 2.3.
The numeric equality in the analytic setting is a mere corollary: Under the
convergence conditions stated in Remark 2.3 the formal series in both sides are
absolutely convergent, and hence the fact that they coincide termwise implies
that their sums are equal.

Our proof is based on the following well known statement:

Let $\X$ be a finite set; let $L$ be a $|\X|\times|\X|$ matrix with rows and
columns marked by the points of $\X$ with matrix elements from an algebra $\Cal
A$, and assume that $\det(1+L)$ is an invertible element of $\Cal A$. Consider
an $\Cal A$--valued measure on the set $2^\X$ of all subsets of $\X$ given by
$$
\p(X)=\frac {\det L_X}{\det(1+L)}\,,
$$
where $L_X$ is the symmetric submatrix of $L$ corresponding to $X$:
$$
L_X=\Vert L(x_i,x_j)\Vert_{x_i,x_j\in X}.
$$
Then the correlation functions of this measure are also given by minors of a
matrix: For any $Y\subset \X$
$$
\rho(Y)=\p\{X\in 2^\X\mid Y\subset X\}=\det K_Y,
$$
where $K=L(1+L)^{-1}$.

Proofs of this statement can be found in \cite{Macc, DVJ, BR}.
Measures of the form above are often called {\it $L$-ensembles}.

Interesting examples of $L$-ensembles usually involve an infinite set $\X$, and
thus the above linear algebraic statement needs to be adjusted to the specific
situation at hand. Our case is of the same nature.

Let us take $\X=\Z'\sqcup\dots\sqcup\Z'$ (total of $N+1$ copies) and consider
the matrix $L$ which has the following block form corresponding to this
splitting of $\X$:
$$
L=\bmatrix 0&0&0&\dots&0& Q\\
-L[1]&0&0&\dots&0&0\\
0&-L[2]&0&\dots&0&0\\
\dots&\dots&\dots&\dots&\dots&\dots\\
0&0&0&\dots&-L[N]&0
\endbmatrix.
$$
The $\Z'\times\Z'$ matrices $L_1,\dots,L_N$ are Toeplitz and their matrix
elements are given by
$$
(L[k])_{xy}=\text{  coefficient of $\zeta^{x-y}$ in  } \exp\sum_{n\ge 1}
\Bigl(a_n[k]\zeta^{-n}+b_n[k]\zeta^n\Bigr),
$$
while the matrix $Q$ is diagonal: $Q_{xy}=zt^x\delta_{xy}$.

As the algebra $\Cal A$ we choose the algebra of formal series in
$\Lambda^{\otimes(2N)}[t^{\pm 1}]$ which have at most a finite order pole at
$t=0$ (in other words, the degrees of $t$ entering any element of $\Cal A$ must
be bounded from below). We will also be using the subalgebra $\Cal A_{hol}$ of
$\Cal A$ which consists of series holomorphic at $t=0$ (i.e., the series which
do not contain negative powers of $t$). The algebra $\Cal A_{hol}$ has a
natural $\Z_+$-filtration induced by the degrees of symmetric functions and
polynomials in $t$. We will denote its filtered components by $\Cal
A_{hol}(d)$, $d\ge 0$. That is, $\Cal A_{hol}(d)$ consists of series whose
terms are of degree at least $d$. The algebra $\Cal A$ also has a natural
topology: Two series are close if their difference is in $\Cal A_{hol}(d)$ for
$d$ large enough. With respect to this topology the algebra $\Cal A$ is
complete.

Recall also that the parameter $z$ entering $Q$ is considered numeric, which
means that it does not contribute to the degree count.

The connection between the matrix $L$ above and the periodic Schur process is
explained by the following statement.

\proclaim{Lemma 2.9} For any partitions $\lambda$ and $\nu$, an integer $l\ge
\max\{\ell(\lambda),\ell(\nu)\}$, and any $k\in\{1,\dots,N\}$ we have
$$
\multline \det \bigl[(L[k])_{\nu_i-i,\lambda_j-j}\bigr]_{i,j=1}^l\\
=H(a[k];b[k]) \sum_{\mu\in\Y}s_{\lambda/\mu}(a[k])s_{\nu/\mu}(b[k])+\Cal
O\bigl(2l-\ell(\lambda)-\ell(\nu)+2\bigr),
\endmultline
$$
where we use the notation $\Cal O(d)$ for elements of $\Cal A_{hol}(d)$.
\endproclaim
\demo{Proof} To simplify the notation we will omit the index
``$[k]$'' in the formulas below. By \cite{Macd, Ex.~I.5.26} we
have
$$
H(a;b) \sum_{\mu}s_{\lambda/\mu}(a)s_{\nu/\mu}(b)=
\sum_{\rho}s_{\rho/\lambda}(b)s_{\rho/\nu}(a).
$$
Let us split the sum in the right-hand side into two parts ---
over partitions $\rho$ of length $\le l$ and $>l$. The second part
has terms of degree at least $2(l+1)-\ell(\lambda)-\ell(\nu)$ and
thus can be ignored. The first part can be rewritten using the
Jacobi-Trudi formula (see \cite{Macd, \S{I} (5.4)}) as:
$$
\sum_{\rho_1\ge\rho_2\ge\dots\ge \rho_l\ge
0}\det[h_{\rho_i-i-\lambda_j+j}(b)]_{i,j=1}^l
\det[h_{\rho_i-i-\nu_j+j}(a)]_{i,j=1}^l\,,
$$
where $h_n$'s are the homogeneous symmetric functions; $h_n=0$ for $n<0$. By
the Cauchy-Binet formula the last sum is readily seen to be equal to
$$
\det \Biggl[\sum_{k\ge-l}
h_{k-\lambda_i+i}(b)h_{k-\nu_j+j}(a)\Biggr]_{i,j=1}^l,
$$
and the sum inside the determinant is equal to the needed matrix element of $L$
because for $p,q\ge -l$ we have
$$
\sum_{k\ge-l} h_{k-p}(b)h_{k-q}(a)=\sum_{k\in\Z}
h_{k-p}(b)h_{k-q}(a)=\sum_{k\in\Z} h_{k+q-p}(b)h_{k}(a)
$$
and
$$
\multline \sum_{n\in\Z}\sum_{k\in\Z}
h_{k+n}(b)h_{k}(a)\zeta^{n}=\sum_{k_1\in\Z}h_{k_1}(b)\zeta^{k_1}
\sum_{k_2\in\Z}h_{k_2}(a)\zeta^{-k_2}\\=\exp\sum_{n\ge 1}
\left(p_n(a)\zeta^{-n}+p_n(b)\zeta^{n}\right). \qed
\endmultline
$$
\enddemo

It is convenient to introduce separate notations for the two values of the
kernel entering the statement of Theorem 2.2: Define $K_+(\sigma,x;\tau,y)$ and
$K_-(\sigma,x;\tau,y)$ through the generating functions
$$
\gathered \sum_{x,y\in\Z'}K_+(\sigma,x;\tau,y)\,\zeta^x\eta^y=
\dfrac{F(\sigma,\zeta)}{F(\tau,\eta^{-1})}\, \sum\limits_{m\in\Z'}
\dfrac{(\zeta\eta)^m}{1+(zt^{m})^{-1}}\,,\\
\sum_{x,y\in\Z'}K_-(\sigma,x;\tau,y)\,\zeta^x\eta^y=
-\dfrac{F(\sigma,\zeta)}{F(\tau,\eta^{-1})}\,\sum\limits_{m\in\Z'}
\dfrac{(\zeta\eta)^m}{1+zt^m}\,.
\endgathered
$$
We will denote by $K_\pm[\sigma,\tau]$ the $\Z'\times\Z'$ matrices with matrix
elements $$(K_\pm[\sigma,\tau])_{xy}=K_\pm(\sigma,x;\tau,y).$$ Here $\sigma$
and $\tau$ are allowed to take values between $0$ and $N$.

\proclaim{Lemma 2.10} (i) For any $\sigma\in\{0,\dots,N\}$ we have
$K_+[\sigma,\sigma]=\id+K_-[\sigma,\sigma]$. \smallskip

(ii) For any $\sigma,\tau\in\{0,\dots,N-1\}$ we have
$L[\sigma+1]K_\pm[\sigma,\tau]=K_\pm[\sigma+1,\tau]$. \smallskip

(iii) For any $\tau\in\{0,\dots,N\}$ we have $QK_-[N,\tau]=-K_+[0,\tau]$.
\endproclaim
\demo{Proof} All these statements are proved by simple algebraic manipulations.

For (i) we have
$$
\sum_{x,y\in\Z'}(K_+(\sigma,x;\sigma,y)-K_-(\sigma,x;\sigma,y))\,\zeta^x\eta^y=
\dfrac{F(\sigma,\zeta)}{F(\sigma,\eta^{-1})}\, \sum\limits_{m\in\Z'}
{(\zeta\eta)^m}=\sum\limits_{m\in\Z'} {(\zeta\eta)^m}
$$
whence $K_+[\sigma,\sigma]-K_-[\sigma,\sigma]=\id$.

The formula (ii) follows from the definitions of $L[\sigma]$,
$K_\pm[\sigma,\tau]$, and the fact that
$$
\Bigl(\exp\sum_{n\ge 1}
\bigl(a_n[\sigma+1]\zeta^{-n}+b_n[\sigma+1]\zeta^n\bigr)\Bigr)
F(\sigma,\zeta)=F(\sigma+1,\zeta).
$$

Finally, for (iii) we have
$$
\multline \sum_{x,y\in\Z'}(QK_-[N,\tau])_{xy}\,\zeta^x\eta^y=
-\dfrac{zF(N,t\zeta)}{F(\tau,\eta^{-1})}\,\sum\limits_{m\in\Z'}
\dfrac{(t\zeta\eta)^m}{1+zt^m}\\ =-\dfrac{F(0,\zeta)}
{F(\tau,\eta^{-1})}\,\sum\limits_{m\in\Z'} \dfrac{(\zeta\eta)^m}{1+(zt^m)^{-1}}
=-\sum_{x,y\in\Z'}(K_+[0,\tau])_{xy}\,\zeta^x\eta^y.\qed
\endmultline
$$
\enddemo

The relations of Lemma 2.10 immediately imply that if we introduce a
$\X\times\X$ matrix $K$ which has the block form
$$
K=\bmatrix K_+[0,0]&K_+[0,1]&K_+[0,2]&\dots&K_+[0,N]\\
           K_-[1,0]&K_+[1,1]&K_+[1,2]&\dots&K_+[1,N]\\
           K_-[2,0]&K_-[2,1]&K_+[2,2]&\dots&K_+[2,N]\\
           \dots&\dots&\dots&\dots&\dots\\
           K_-[N,0]&K_-[N,1]&K_-[N,2]&\dots&K_+[N,N]
\endbmatrix
$$
(the $(\sigma,\tau)$-block is equal to $K_+[\sigma,\tau]$ if $\sigma\le \tau$,
and to $K_-[\sigma,\tau]$ otherwise), then we have the matrix relation
$(\id+L)K=L$. In order to extract the probabilistic meaning of this relation we
need to introduce certain finite point approximations of the shift-mixed
periodic Schur process.

For any $m=1,2,\dots$ denote by $\Z'(m)$ the subset of $\Z'$ consisting of $2m$
half-integers situated symmetrically around 0:
$$
\Z'(m)=\left\{-m+\tfrac 12,-m+\tfrac 32,\dots,-\tfrac 12,\tfrac
12,\dots,m-\tfrac32,m-\tfrac 12\right\}.
$$

Denote also by $L^{(m)}[\sigma]$, $K_\pm^{(m)}[\sigma,\tau]$, and $Q^{(m)}$ the
restrictions of the $\Z'\times\Z'$ matrices $L[\sigma]$, $K_\pm[\sigma,\tau]$,
and $Q$ to $\Z'(m)\times\Z'(m)$, and denote by $L^{(m)}$ and $K^{(m)}$ the
block matrices built from $L^{(m)}[\sigma]$, $K_\pm^{(m)}[\sigma,\tau]$,
$Q^{(m)}$ in the same way as $L$ and $K$ are built from $L[\sigma]$,
$K_\pm[\sigma,\tau]$, $Q$. Set
$$
\X^{(m)}=\Z'(m) \sqcup\dots\sqcup\Z'(m)\quad \text{($N+1$ copies).}
$$
Note that
$$
L^{(m)}=L\bigl|_{\X^{(m)}\times \X^{(m)}}\,,\qquad
K^{(m)}=K\bigl|_{\X^{(m)}\times \X^{(m)}}\,.
$$

\proclaim{Lemma 2.11} (i) $\det(1+L^{(m)})=z^mt^{-\frac{m^2}2}\bigl(1+\Cal
O(1)\bigr)$ where $\Cal O(1)$ stands for an element of $\Cal A_{hol}(1)$.

(ii) All matrix elements of $(1+L^{(m)})^{-1}$ belong to $\Cal A_{hol}$.
Consequently, for any $X\subset \X^{(m)}$ we have $\det
L^{(m)}_X/\det(1+L^{(m)})\in \Cal A_{hol}$.
\endproclaim

\demo{Proof} (i) Since $\det(1+L^{(m)})=\sum_{X\subset \X^{(m)}}\det
L^{(m)}_X$, we need to see which minors $\det L^{(m)}_X$ yield terms of the
lowest possible degree. Since
$$
(L[k])_{xy}\in \Cal A_{hol}(|x-y|)\quad \text{and}\quad (L[k])_{xx}=1+\Cal
O(1),\qquad k=1,\dots,N,
$$
it is immediate that the only lowest degree term comes from $\det L^{(m)}_X$
with
$$
X=\left\{-m+\tfrac 12,\dots,-\tfrac 12\right\}\sqcup\dots\sqcup
\left\{-m+\tfrac 12,\dots,-\tfrac 12\right\}
$$
($N+1$ copies of the same set of negative elements in $\Z'(m)$), and it is
equal to $z^mt^{-\frac{m^2}2}$ (recall that $z$ does not contribute to the
degree count).

(ii) Matrix elements of $(1+L^{(m)})^{-1}$ are ratios of the linear
combinations of minors of $L^{(m)}$ and $\det(1+L^{(m)})$. Since all minors of
$L^{(m)}$ lie in $t^{-\frac{m^2}2}\cdot \Cal A_{hol}$, by (i) we see that the
matrix elements belong to $\Cal A_{hol}$.

As for the ratios $\det L_X/\det(1+L^{(m)})$, observe that they are linear
combinations of ratios of the form $\det (1+L^{(m)})_X/\det(1+L^{(m)})$ which
coincide, up to a sign, with minors of the inverse matrix
$(1+L^{(m)})^{-1}$.\qed
\enddemo

We are now in a position to prove that $L^{(m)}$-ensembles on $\X^{(m)}$
approximate the shift-mixed periodic Schur process on $\X$ as $m$ becomes
large.

\proclaim{Proposition 2.12} The values of the correlation functions of the
$L^{(m)}$-ensembles restricted to the last $N$ copies of $\Z'(m)$ in
$\X^{(m)}=\Z'(m)\sqcup\dots\sqcup\Z'(m)$ converge, as $m\to\infty$, to those of
the dynamical correlation functions of the shift-mixed periodic Schur process
in the topology of $\Cal A$.
\endproclaim

\demo{Proof} In order to prove this statement we will construct an injective
map that associates to any point configuration $X\subset \X^{(m)}$ of nonzero
weight in the $L^{(m)}$-ensemble a sequence of partitions
$(\lambda^{(1)},\dots,\lambda^{(N)})$ and an integer $S$ so that:

$\bullet$\  The intersection of the point configuration
$$
\bigl\{S+\lambda_i^{(1)}-i+\tfrac 12\bigr\}_{i\ge
1}\sqcup\dots\sqcup\bigl\{S+\lambda_i^{(N)}-i+\tfrac 12\bigr\}_{i\ge 1}
$$
with $\Z'(m)\sqcup\dots\sqcup\Z'(m)$ coincides with the restriction of $X$ to
the last $N$ copies of $\Z'(m)$ in $\X^{(m)}$.

$\bullet$\  The weight of $X$ in the $L^{(m)}$-ensemble and the weight of
$(\lambda,S)$ in the shift-mixed periodic Schur process are obtained from each
other by multiplication by a constant of the form $C_m=1+\Cal O(1)$ and by
addition of an element of a high enough degree:
$$
\frac {\det L_X^{(m)}}{\det(1+L_X^{(m)})}=C_m\cdot z^S t^{\frac
{S^2}2}\sum_{\mu}\w(\lambda,\mu) + \Cal O(d(m)),
$$
where $d(m)$ does not depend on $X$ and $d(m)\to\infty$ as
$m\to\infty$.\footnote{We should have used $C_m/(\theta_3(z;t)Z(N,t,a,b))$
instead of $C_m$ thus taking into account the partition function. However, this
is equivalent because $\theta_3(z;t)Z(N,t,a,b)$ and
$1/(\theta_3(z;t)Z(N,t,a,b))$ are both of the form $1+\Cal O(1)$.}

$\bullet$\ The lowest degree of the weights of pairs $(\lambda,S)$ not covered
by this map goes to infinity as $m\to\infty$.

Since all the weights add up to 1 in both the $L^{(m)}$-ensemble and the
shift-mixed periodic Schur process, the existence of such a map implies that
the degree of $C_m-1$ goes to infinity as $m\to \infty$, and the needed
convergence of the correlation functions readily follows.

Observe that if a set $X\subset \X^{(m)}$ has a nonzero weight in the
$L^{(m)}$-ensemble (i.e., $\det L_X^{(m)}\ne 0$) then its intersections with
all $N+1$ copies of $\Z'(m)$ in $\X^{(m)}$ must have the same cardinality
because of the specific block structure of $L^{(m)}$, and its intersections
with the first and the last copies of $\Z'(m)$ must coincide because the matrix
$Q$ is diagonal. Thus, without loss of generality we may assume that
$$
X=\{x_1^{(0)},\dots,x_l^{(0)}\}\sqcup \{x_1^{(1)},\dots,x_l^{(1)}\}\sqcup
\{x_1^{(2)},\dots,x_l^{(2)}\}\sqcup\dots\sqcup \{x_1^{(N)},\dots,x_l^{(N)}\}
$$
for some $l\le 2m$, where $x_i^{(0)}=x_i^{(N)}$ for all $1\le i\le l$. Our
notation  for $X$ means that the first group of points lies in the first copy
of $\Z'(m)$, the second one lies in the second copy of $\Z'(m)$ and so on.

The element $(\lambda^{(1)},\dots,\lambda^{(N)};S)$ of the shift-mixed process
corresponding to such $X$ is defined as follows: $S=l-m$ and
$$
\bigl\{S+\lambda^{(\tau)}_i-i+\tfrac 12\bigr\}_{i\ge
1}=\{x_1^{(\tau)},\dots,x_l^{(\tau)}\}\cup\bigl\{-m-\tfrac 12,-m-\tfrac
32,-m-\tfrac 52,\dots\bigr\}
$$
for any $\tau=1,\dots,N$. It is readily seen that this formula correctly
defines the partitions $\lambda^{(\tau)}$. Note that $\ell(\lambda^{(\tau)})\le
l$ for all $\tau$.

The pair $(\lambda^{(1)},\dots,\lambda^{(N)};S)$ is not covered by this map if
and only if for some $\tau$ the point configuration
$\bigl\{S+\lambda^{(\tau)}_i-i+\tfrac 12\bigr\}_{i\ge 1}$ either does not
contain the set $\Z'_{<-m}$ or has a nonzero intersection with $\Z'_{>m}$. In
the latter case we must have $\lambda_1^{(\tau)}+S\ge m+1$ and in the former
case we must have $(\lambda^{(\tau)})'_1-S\ge m+1$. This means that either
$|S|$ or $|\lambda^{(\tau)}|$ is $\ge \frac m2$.

The definition of the shift-mixed process implies that the weight of
$(\lambda,S)$ has degree at least
$$
\tfrac{S^2}2 +\bigl||\lambda^{(1)}|-|\lambda^{(2)}|\bigr|+\dots+
\bigl||\lambda^{(N-1)}|-|\lambda^{(N)}|\bigr|+|\lambda^{(N)}|\le
\tfrac{S^2}2+\max_{1\le \tau\le N}|\lambda^{(\tau)}|.
$$
Therefore, the pairs $(\lambda,S)$ not in the image of our map have weights of
degrees uniformly going to infinity as $m\to\infty$.

It remains to compare the weights of a point configuration $X\subset \X^{(m)}$
in the $L^{(m)}$-ensemble and of its image in the shift-mixed process.

Using Lemma 2.11(i) we obtain that the weight ${\det
L_X^{(m)}}/{\det(1+L^{(m)})}$ of $X$ equals
$$
\const\cdot z^{-m}t^{\frac {m^2}2} \det
\bigl[(L[1])_{x_i^{(1)},x_j^{(0)}}\bigr]_{i,j=1}^l\cdots \det
\bigl[(L[N])_{x_i^{(N)},x_j^{(N-1)}}\bigr]_{i,j=1}^l \prod_{i=1}^l
zt^{x^{(N)}_i}.
$$
with $\const=1+\Cal O(1)$. Collecting the powers of $z$ and $t$ yields
$$
\multline z^{l-m}t^{\frac {m^2}2+(S+\frac 12)l+\sum_{i=1}^l
(\lambda_i^{(N)}-i)}\\= z^{l-m}t^{\frac {m^2}2+(l-m+\frac
12)l-\frac{l(l+1)}2+\sum_{i=1}^l \lambda_i^{(N)}}=z^St^{\frac
{S^2}2+|\lambda^{(N)}|}
\endmultline
$$
where we used the definition of $S=l-m$ and the fact that
$\ell(\lambda^{(N)})\le l$.\footnote{This is the computation mentioned in
Remark 2.7 above.} Evaluating the minors of $L[\tau]$ by Lemma 2.9 (which is
applicable because $\ell(\lambda^{(\tau)})\le l$) we obtain
$$
\multline \det \bigl[(L[1])_{x_i^{(1)},x_j^{(0)}}\bigr]_{i,j=1}^l\cdots \det
\bigl[(L[N])_{x_i^{(N)},x_j^{(N-1)}}\bigr]_{i,j=1}^l=\const\\
\times\left(\sum_{\mu^{(1)}} s_{\lambda^{(0)}/\mu^{(1)}}(a[1])\,
s_{\lambda^{(1)}/\mu^{(1)}}(b[1])+
\Cal O\bigl(2l-\ell(\lambda^{(0)})-\ell(\lambda^{(1)})+2\bigr)\right)\\
\times\left(\sum_{\mu^{(2)}} s_{\lambda^{(1)}/\mu^{(2)}}(a[2])\,
s_{\lambda^{(2)}/\mu^{(2)}}(b[2])+ \Cal
O\bigl(2l-\ell(\lambda^{(1)})-\ell(\lambda^{(2)})+2\bigr)\right)\cdots \\\times
\left(\sum_{\mu^{(N)}} s_{\lambda^{(N-1)}/\mu^{(N)}}(a[N])\,
s_{\lambda^{(N)}/\mu^{(N)}}(b[N])+ \Cal
O\bigl(2l-\ell(\lambda^{(N-1)})-\ell(\lambda^{(N)})+2\bigr)\right)
\endmultline
$$
with $\const=1+\Cal O(1)$. To conclude the proof we need to show that all the
remainders $\Cal O(\,\cdot\,)$ in the expression above can be removed by adding
the correction of degree uniformly going to infinity as $m\to\infty$.

Observe that the degree of the factor $t^{\frac{S^2}2+|\lambda^{(N)}|}$ is
bounded if and only if $|S|=|l-m|$ and $|\lambda^{(N)}|$ are bounded. Then in
the last factor in the product above the term
$O\bigl(2l-\ell(\lambda^{(N-1)})-\ell(\lambda^{(N)})+2\bigr)$ has degree going
to infinity with $m$ (because $l-\ell(\lambda^{(N)})\ge m-S-|\lambda^{(N)}|$),
so if we want the weight of $X$ or the weight of $(\lambda,S)$ to be of bounded
degree, the degree of the sum over $\mu^{(N)}$ has to be bounded. Since
$$
\deg\sum_{\mu}s_{\lambda/\mu}(a)s_{\nu/\mu}(b)\ge \bigl ||\lambda|-|\nu|\bigr|,
$$
this means that $|\lambda^{(N-1)}|$ has to be bounded. Repeating the argument
with the second to last factor we conclude that $|\lambda^{(N-2)}|$ is bounded,
and so on.

The final conclusion is that for the corresponding $X$ and $(\lambda,S)$ such
that the minimum of the degrees of their weights is bounded, we must have
$|l-m|$ and $|\lambda^{(\tau)}|$, $\tau=1,\dots,N$, bounded. But then all $\Cal
O(\,\cdot\,)$'s in the product above have degrees uniformly going to infinity
as $m\to\infty$, and thus the difference of the weights of $X$ and
$(\lambda,S)$ has degree uniformly going to infinity as $m\to\infty$. The proof
of Proposition 2.12 is complete. \qed
\enddemo

Let us now conclude the proof of Theorem 2.2. From the determinantal formula
for the correlation functions of general $L$-ensembles, see the beginning of
the proof, we know that the correlation functions of the $L^{(m)}$-ensembles
are given by minors of the matrices $L^{(m)}(1+L^{(m)})^{-1}$. The last step of
the proof is to show that matrix elements of these matrices converge to those
of the kernel $K$ in the topology of $\Cal A$.

The definitions of matrices $L[\sigma]$ and $K[\sigma,\tau]$ imply that
$$
\deg (L[\sigma])_{xy}\ge |x-y|,\qquad (K[\sigma,\tau])_{xy}\ge |x-y|,\qquad
x,y\in\Z',
$$
for any $\sigma$ and $\tau$. Hence, from Lemma 2.10(ii) we obtain
($x,y\in\Z'(m)$)
$$
\multline \bigl(L^{(m)}[\sigma+1]K^{(m)}_\pm[\sigma,\tau]\bigr)_{xy}=
(K^{(m)}_\pm[\sigma+1,\tau])_{xy}\\ +\Cal O\bigl(\min\{m+\tfrac
12-x,\,x-m-\tfrac 12\}+\min\{m+\tfrac 12-y,\,y-m-\tfrac 12\}\bigr).
\endmultline
$$
The relations (i) and (iii) clearly remain unchanged when restricted to
$\Z'(m)$:
$$
K_+^{(m)}[\sigma,\sigma]=\id+K^{(m)}_-[\sigma,\sigma],\qquad
Q^{(m)}K^{(m)}_-[N,\tau]=-K^{(m)}_+[0,\tau].
$$
Hence,
$$
\bigl((1+L^{(m)})K^{(m)}\bigr)_{xy}=L^{(m)}_{xy}+\Cal O\bigl(\min\{m+\tfrac
12-y,\,y-m-\tfrac 12\}\bigr).
$$
Multiplying by $(\id+L^{(m)})^{-1}$ on the left and using Lemma 2.11(ii) we
obtain
$$
K^{(m)}_{xy}=\bigl((1+L^{(m)})^{-1}L^{(m)}\bigr)_{xy}+\Cal
O\bigl(\min\{m+\tfrac 12-y,\,y-m-\tfrac 12\}\bigr).
$$
Since $\min\{m+\tfrac 12-y,\,y-m-\tfrac 12\}\to\infty$ as $m\to\infty$ for any
fixed $y$, the needed convergence follows. The proof of Theorem 2.2 is
complete.\qed
\enddemo

\head 3. Bulk scaling limit
\endhead

As was mentioned in \S1, the presence of the parameter $t$ is crucial for
defining the periodic Schur process; at $t=1$ the partition function is,
generally speaking, infinite. Thus, one might expect that if $t\to 1-$ then the
random Young diagrams become large. Indeed, this is correct: In the analytic
setting we will show that as $t\to 1-$, the density function of the scaled
random point configuration $|\ln t|\cdot \L(\lambda^{(\tau)})$ tends to a
nontrivial limit, and this limit is independent of  $\tau$. The main result of
this section is the computation of the local limit of the correlation functions
of the periodic Schur process and its shift-mixed modification near points of
fixed global limit density.

Denote $r=\ln t^{-1}$. Throughout this section we assume that
$r>0$ (equivalently, $0<t<1$), and also that the convergence
conditions of Remark 2.3 are satisfied. Namely, we assume that $
|a_n[k]|,|b_n[l]|<\const\cdot R^n  \, \text{ for some $0<R<1$ and
all } \, 1\le k,l\le N$, and we also assume that $|\arg(z)|<\pi$.

Recall the notation $A_k=\sum_{m=1}^N a_k[m]$, $B_k=\sum_{m=1}^N b_k[m]$,
introduced in \S1.

\proclaim{Theorem 3.1} (i) Assume that $A_k=\overline{B}_k$ for
all $k=1,2,\dots$\,. Then, as $t\to 1$, the dynamical correlation
functions of the shift-mixed periodic Schur process have a limit
in the following sense: Choose $x_1(t),\dots,x_n(t)\in\Z'$ such
that as $t\to 1$, $r x_k(t)\to \gamma$ for all $k=1,\dots,n$ and
some $\gamma\in\R$, and all pairwise distances
$x_i-x_j=x_i(t)-x_j(t)$ are independent of $t$. Then for any $1\le
\tau_1,\dots,\tau_n\le N$
$$
\lim_{t\to
1}\rho_n^{\shift}(\tau_1,x_1(t);\dots,\tau_n,x_n(t))=\det\bigl[
 \K_{\tau_i,\tau_j}^{(z,\gamma)}(x_i-x_j)\bigr]_{i,j=1}^n,
$$
where the correlation kernel has the following form
$$
\K_{\sigma,\tau}^{(z,\gamma)}(d)=\cases\frac 1{2\pi i}\oint\limits_{|\zeta=1|}
\frac{\exp\bigl({-\sum_{n\ge
1}\sum_{k=\sigma+1}^{\tau}(a_n[k]\zeta^{-n}+b_n[k]\zeta^n)}\bigr)}
{1+z^{-1}\exp\bigl({\gamma-\sum_{n\ge
1}(A_n\zeta^{-n}+B_n\zeta^n)}\bigr)}\,\frac{d\zeta}{\zeta^{d+1}},&\sigma\le
\tau,\\
-\frac 1{2\pi i}\oint\limits_{|\zeta|=1} \frac{\exp\bigl({\sum_{n\ge
1}\sum_{k=\tau+1}^{\sigma}(a_n[k]\zeta^{-n}+b_n[k]\zeta^n)}\bigr)}
{1+z\exp\bigl({-\gamma+\sum_{n\ge
1}(A_n\zeta^{-n}+B_n\zeta^n)}\bigr)}\,\frac{d\zeta}{\zeta^{d+1}},&\sigma>\tau.
\endcases
$$
(ii) Under the same assumptions the dynamical correlation functions of the
periodic Schur process converge to the limiting expression above evaluated at
$z=1$:
$$
\lim_{t\to 1}\rho_n(\tau_1,x_1(t);\dots,\tau_n,x_n(t))=\det\bigl[
 \K_{\tau_i,\tau_j}^{(1,\gamma)}(x_i-x_j)\bigr]_{i,j=1}^n.
$$
\endproclaim

\demo{Comments} 1. The limit correlation functions as functions on
$\Z\sqcup\dots\sqcup\Z$ ($N$ copies) are invariant with respect to the
simultaneous shifts of all variables.

2. It will be clear from the proof that a slightly more general statement is
true: For several groups of variables $\{x_k^{(m)}(t)\}_{k=1}^{n_m}$,
$m=1,\dots,M$, such that $rx_k^{(m)}\to \gamma_m$, the pairwise distances
inside every group are independent of $t$, and the distances between different
groups tend to infinity:
$$
x_i^{(m)}(t)-x_j^{(m)}(t)=x_i^{(m)}-x_j^{(m)},\qquad \min_{i,j}
\bigl|x_i^{(m_1)}(t)-x_j^{(m_2)}(t)\bigr|\to\infty \text{  for  } m_1\ne m_2,
$$
the limit of the dynamical correlation functions of the shift-mixed process is
the product of determinants:
$$
\lim_{t\to
1}\rho_n^{\shift}\bigl(\{\tau_i^{(m)},x_i^{(m)}(t)\}_{1\le i\le
n_m,\, 1\le m\le M}\bigr)=\prod_{m=1}^M\det\Bigl[
 \K_{\tau_i^{(m)},\tau_j^{(m)}}^{(z,\gamma_m)}
 (x_i^{(m)}-x_j^{(m)})\Bigr]_{i,j=1}^{n_m},
$$
and the same is true for the correlation functions of the periodic Schur
process with $z=1$ in the right-hand side of the formula above. Roughly
speaking, this means that particles become independent as the distance between
them grows.

3. The global limit density function mentioned in the beginning of the section
is equal to
$$
\rho^{(z)}(\gamma)=\frac 1{2\pi i}\oint_{|\zeta|=1}
\frac{1}{1+z^{-1}e^{\gamma -\sum_{n\ge
1}(A_n\zeta^{-n}+B_n\zeta^n)}}\,\frac{d\zeta}{\zeta}
$$
for the shift-mixed process, and to $\rho^{(1)}(\gamma)$ for the
unmixed process. This formula has the following corollary: If one
assumes the existence of the limit shape, as $t\to 1$, of the
random Young diagrams $\lambda^{(\tau)}$ distributed according to
the periodic Schur process, then this limit shape can be easily
guessed. Denote by $i$ and $j$ the row and column coordinates on
the diagrams and introduce new coordinates $u=r(j-i)$ and
$v=r(i+j)$. Then the equation for the boundary of the hypothetical
limit shape has the form
$$
v(u)=u+2\oint_{|\zeta|=1}\ln\left(1+e^{-u+\sum_{n\ge
1}(A_n\zeta^{-n}+B_n\zeta^n)}\right)\,\frac{d\zeta}{\zeta}\,.
$$
This formula follows from the relation $\rho(u)=\frac12(1-v'(u))$,
see \cite{BOO, Remark 1.7} for an explanation.

4. In the case of the uniform measure on partitions, when the
specializations $a[k],b[k]$ are trivial, Theorem 3.1 (or rather
its extended version from Comment 1 above) coincides with Theorem
7 of \cite{O1}. In this case the limit correlation kernel
degenerates:
$$
\Cal K_{\sigma,\tau}^{(z,\gamma)}(d)=\cases
\delta_{0d}\cdot(1+z^{-1}e^\gamma)^{-1},& \sigma\le \tau,\\
-\delta_{0d}\cdot (1+ze^{-\gamma})^{-1},&\sigma>\tau.
\endcases
$$

For the uniform measure the limit shape does exist, see \cite{V},
and the formula for $v(u)$ from the previous comment produces the
correct answer: $v(u)=u+2\ln(1+e^{-u})$ or
$e^{\frac{v-u}2}+e^{\frac{u+v}2}=1$.
\enddemo

\demo{Proof of Theorem 3.1} (i) We start with the integral representation of
the correlation kernel for the shift-mixed process from Remark 2.4. Replacing
the integration variable $\eta$ by $\xi=\zeta\eta$ we obtain
$$
K(\sigma,x;\tau,y)=-\frac{\prod_{n\ge 1}(1-t^n)^3}{\theta_3(z;t)
\,(2\pi
i)^2}\oint\limits_{|\zeta|=1}\oint\limits_{|\xi|=1\pm\epsilon}
\dfrac{F(\sigma,\zeta)}{F(\tau,\zeta\xi^{-1})}
\,\frac{\theta_3(z\xi;t)}{\theta_3(-\xi t^{-\frac 12};t)}\,
\frac{d\zeta d\xi}{\zeta^{x-y+1}\xi^{y+\frac12}}\,,
$$
where $\epsilon>0$ is much smaller then $r$, and we choose $|\xi|=1+\epsilon$
for $\sigma\le \tau$ and $|\xi|=1-\epsilon$ for $\sigma>\tau$. The next step is
to fix an arbitrary $\zeta$ on the unit circle and to evaluate the asymptotics
of the integral over $\xi$.

\proclaim{Proposition 3.2} (i) Assume that $|\arg z|<\pi$. Then, as $t\to 1$,
the function
$$
f(\xi)=-\frac{\prod_{n\ge
1}(1-t^n)^3\,\theta_3(z\xi;t)}{\theta_3(z;t)\,\theta_3(-\xi t^{-\frac 12};t)}
$$
on the circle $|\xi|=1\pm\epsilon$ uniformly converges to $0$ on the complement
to any neighborhood of the point $\xi=1\pm\epsilon$. On the other hand, there
exists $\delta>0$ such that for $|\arg (\xi)|<\delta$
$$
f(\xi)=\frac{2\pi i\, {\xi}^{-\frac 12} }r \,\frac{
e^{\frac{\ln(z)\ln(\xi)}{r}}}{e^{\frac{\pi i\ln(\xi)}r}- e^{-\frac{\pi
i\ln(\xi)}r}}\cdot(1+f_0(\xi))
$$
where $f_0(\xi)$ is an analytic function which, as $t\to 1$, uniformly
converges to 0 while $rf'_0(\xi)$ remains uniformly bounded.\footnote{Here we
use the principal branch of the logarithm function for both $\ln(z)$ and
$\ln(\xi)$.} All the estimates are uniform in $z$ varying in a compact set of
the complex plane bounded away from the negative real semiaxis.

(ii) Assume that $z$ lies on the unit circle $|z|=1$. Then on the
circle $|\xi|=1\pm\epsilon$, as $t\to 1$, we have the bound
$|f(\xi)|\le \const\cdot \epsilon^{-1}$ which is uniform in both
$z$ and $\xi$.

\endproclaim
\demo{Proof} (i) Denote $v=\frac 1{2\pi i} \ln\xi$ and $w=\frac
1{2\pi i}\ln z$. Applying the imaginary Jacobi transform, see e.g.
\cite{Er, 13.22(8)}, we obtain
$$
\gathered \theta_3(z;t)=\left(\tfrac{2\pi} r\right)^\frac12
e^{-\frac{2\pi^2w^2}{r}}\theta_3\Bigl(e^{\frac{4\pi^2 w}r}\,;
e^{-\frac{4\pi^2}r}\Bigr),\\
\theta_3(z\xi;t)=\left(\tfrac{2\pi} r\right)^\frac12
e^{-\frac{2\pi^2(v+w)^2}{r}}\theta_3\Bigl(e^{\frac{4\pi^2 (v+w)}r}\,;
e^{-\frac{4\pi^2}r}\Bigr),\\
\theta_3(-\xi t^{-\frac12};t)=\left(\tfrac{2\pi} r\right)^\frac12
e^{-\frac{2\pi^2(v-\frac{ir}{4\pi}+\frac12)^2}{r}}\theta_3\Bigl(-e^\frac{4\pi^2
(v+\frac12)}r\,; e^{-\frac{4\pi^2}r}\Bigr).
\endgathered
$$
Hence,
$$
f(\xi)= \frac{-\left(\tfrac r{2\pi}\right)^\frac 12 \prod\limits_{n\ge
1}(1-t^n)^3\, e^{-\frac{2\pi^2}r((v+w)^2-w^2-(v-\frac{ir}{4\pi}+\frac12)^2)}\,
\theta_3\Bigl(e^{\frac{4\pi^2 (v+w)}r}\,;
e^{-\frac{4\pi^2}r}\Bigr)}{\theta_3\Bigl(e^{\frac{4\pi^2 w}r}\,;
e^{-\frac{4\pi^2}r}\Bigr)\,\theta_3\Bigl(-e^\frac{4\pi^2 (v+\frac12)}r\,;
e^{-\frac{4\pi^2}r}\Bigr)}\,.
$$
We have
$$
\gathered e^{-\frac{2\pi^2}r((v+w)^2-w^2-(v-\frac{ir}{4\pi}+\frac12)^2)}=-i
e^{\frac{2\pi^2}{r}(\frac14+v(1-2w))-\pi i v}\cdot e^{-\frac r8}\,,\\
\lim_{t\to 1} \frac{r^{\frac 32}e^{\frac{\pi^2}{2r}}}{(2\pi)^{\frac
32}}\prod_{n\ge 1}(1-t^n)^3=1.
\endgathered
$$
The last relation can be obtained, for example, from the imaginary Jacobi
transform of $\theta_1(x;t)$ because $\theta_1(x;t)\sim (x^\frac 12-x^{-\frac
12})\,(t;t)_\infty^3$ as $x\to 1$. Thus,
$$
-\left(\tfrac r{2\pi}\right)^\frac 12 \prod\limits_{n\ge 1}(1-t^n)^3\,
e^{-\frac{2\pi^2}r((v+w)^2-w^2-(v-\frac{ir}{4\pi}+\frac12)^2)}=\frac{2\pi i \,
e^{\frac{2\pi^2v(1-2w)}{r}-\pi i v}(1+o(1))}r\,.
$$
Writing down the products for $\theta_3$'s explicitly we obtain
$$
\gathered \theta_3\Bigl(e^{\frac{4\pi^2 w}r}\,;
e^{-\frac{4\pi^2}r}\Bigr)=\Pi_0\cdot\prod_{n\ge 0}
\left(1+e^{\frac{4\pi^2(w-n-\frac12)}r}\right)
\left(1+e^{\frac{4\pi^2(-w-n-\frac12)}r}\right),\\
\theta_3\Bigl(e^{\frac{4\pi^2 (v+w)}r}\,; e^{-\frac{4\pi^2}r}\Bigr)=
\Pi_0\cdot\left(1+e^{\frac{4\pi^2
(v+w-\frac12)}r}\right)\left(1+e^{\frac{4\pi^2 (-v-w-\frac12)}r}\right)\\
\qquad\qquad\qquad\qquad\times\prod_{n\ge 0}
\left(1+e^{\frac{4\pi^2(v+w-n-\frac32)}r}\right)
\left(1+e^{\frac{4\pi^2(-v-w-n-\frac32)}r}\right),\\
\theta_3\Bigl(-e^\frac{4\pi^2 (v+\frac12)}r\,; e^{-\frac{4\pi^2}r}\Bigr)=
\Pi_0\cdot\left(1-e^{\frac{4\pi^2
v}r}\right)\qquad\qquad\qquad\qquad\qquad\qquad\\
\qquad\qquad\qquad\times\prod_{n\ge 0} \left(1-e^{\frac{4\pi^2(v-n-1)}r}\right)
\left(1-e^{\frac{4\pi^2(-v-n-1)}r}\right),
\endgathered
$$
where $\Pi_0=\prod_{n\ge 1}\bigl(1-e^{-\frac{4\pi^2n}r}\bigr)\to 1$ as $r\to
0$. Using the principal branch of the logarithm we may assume that
$$
\Re v=\tfrac 1{2\pi}\arg\xi\in[-\tfrac12,\tfrac12]\quad\text{  and  }\quad\Re
w=\tfrac 1{2\pi}\arg z\in(-\tfrac 12,\tfrac12)
$$
(recall that $|\arg z|<\pi$ by assumption).

It is readily seen that under these conditions, as $r\to 0$, all
three products over $n\ge 0$ in the formulas above uniformly tend
to 1 while their derivatives with respect to $v$ multiplied by $r$
are uniformly bounded. The remaining factors are
$$
f(\xi)\sim \frac{2\pi i \, e^{\frac{2\pi^2v(1-2w)}{r}-\pi i
v}\left(1+e^{\frac{4\pi^2 (v+w-\frac12)}r}\right)\left(1+e^{\frac{4\pi^2
(-v-w-\frac12)}r}\right)}{r\left(1-e^{\frac{4\pi^2 v}r}\right)}\,.
$$

Assume first that $\Re v$ is bounded away from 0 (that is, $\xi$
is bounded away from $\xi=1\pm\epsilon$). Note that since
$\epsilon\ll r$, we have $\Im v=\frac 1{2\pi}
\ln|1\pm\epsilon|=o(r)$ and
$$
\Re (v(1-2w))=\Re v(1-2\Re w)+o(r).
$$
If $\Re v<0$ and $\Re(v+w)\ge-\frac 12$ then the absolute value of
the asymptotic expression for $f(\xi)$ above is bounded by
$${\const}\cdot r^{-1}{e^{{2\pi^2\Re v(1-2\Re w)}r^{-1}}}$$
which uniformly converges to 0 as $r\to 0$. If $\Re v<0$ and $\Re(v+w)<-\frac
12$ then the bound is
$$
{\const}\cdot r^{-1}e^{{2\pi^2(\Re v(1-2\Re w)-2\Re v-2\Re w-1)r^{-1}}}
={\const}\cdot r^{-1}e^{{-2\pi^2(\Re v+1)(2\Re w+1)r^{-1}}}
$$
which also goes to 0 as $r\to 0$. If $\Re v>0$ and $\Re(v+w)\le \frac 12$ then
$|f(\xi)|$ is bounded by
$$
\const\cdot r^{-1}e^{{2\pi^2(\Re v(1-2\Re w)-2\Re v)r^{-1}}}= \const\cdot
r^{-1}e^{-2\pi^2\Re v(1+2\Re w)r^{-1}}
$$
which goes to 0 as $r\to 0$. Finally, if $\Re v>0$ and
$\Re(v+w)>\frac 12$ then the bound has the form
$$
\const\cdot r^{-1}e^{{2\pi^2(\Re v(1-2\Re w)+2(\Re v+\Re w-\frac12)-2\Re
v)r^{-1}}}= \const\cdot r^{-1}e^{-2\pi^2(\Re v-1)(1-2\Re w)r^{-1}}
$$
which is also small as $r\to 0$. This takes care of the first statement (about
uniform convergence to 0) of Proposition 3.2(i).

Since $|\Re w|<\frac 12$, we can choose $\delta>0$ such that $|\Re w\pm \frac
\delta{2\pi}|<\frac 12$. Then for $\xi$ with $|\arg(\xi)|=2\pi |\Re v|<\delta$
we have $\Re(v+w-\frac 12)<0$, and the two factors in the numerator of the
approximation for $f(\xi)$ above uniformly converge to 1 as $r\to 0$ with their
derivatives with respect to $v$ multiplied by $r$ are uniformly bounded. Thus,
$$
f(\xi)\sim \frac{2\pi i \, e^{-\frac{4\pi^2vw}{r}-\pi i
v}}{r\Bigl(e^{\frac{-2\pi^2 v}r}-e^{\frac{2\pi^2 v}r}\Bigr)}\,
$$
as desired. It is also immediately visible that all the estimates above are
uniform in $z$ varying in a compact set bounded away from the negative
semiaxis. The proof of (i) is complete.

The proof of (ii) follows the arguments used above to prove the first part of
(i). There are two differences in estimates: the range of $\Re w$ is now the
whole segment $[-\frac 12,\frac 12]$ rather than a closed subset of the
interval $(-\frac 12,\frac 12)$, and since $\Re v$ is not bounded away from 0,
we have an additional potentially small factor of the form
$\Bigl(1-e^{\pm\frac{4\pi^2 v}r}\Bigr)$ in the denominator. We have
$$
\arg \bigl(e^{\pm\frac{4\pi^2 v}r}\bigr)=\pm\tfrac{4\pi^2\Im v}{r}=\pm \tfrac
{2\pi}r\ln(1\pm\epsilon)\sim \pm\tfrac {2\pi \epsilon}r
$$
and hence $\Bigl| 1-e^{\pm\frac{4\pi^2 v}r}\Bigr|\ge
\const\frac{2\pi\epsilon}r$. Adding this estimate to those derived above yields
the needed bound. \qed
\enddemo

Our next step is to compute the asymptotics of the ratio
${F(\sigma,\zeta)}/{F(\tau,\zeta\xi^{-1})}$ from the integral representation of
the correlation kernel given in the beginning of the proof of Theorem 3.1.
Recall that the definition of the function $F(\tau,\zeta)$ was given just
before Theorem 2.2. After simple manipulations, using the assumption that
$B_k=\overline{A}_k$ for $k\ge 1$, we obtain
$$
\multline \frac {F(\sigma,\zeta)}{F(\tau,\zeta\xi^{-1})}=\exp\sum_{n\ge 1}\frac
1{1-t^n}\Bigl(\overline{A}_n{\zeta^{n}(1-\xi^{-n})} -A_n\zeta^{-n}(1-\xi^n)
\Bigr)\\ \times\,\exp\sum_{n\ge
1}\Biggl({\zeta^{-n}}\sum_{k=1}^{\sigma}a_n[k]-\zeta^n\sum_{k=\sigma+1}^N
b_n[k]-{(\zeta/\xi)^{-n}}\sum_{k=1}^{\tau}a_n[k]+(\zeta/\xi)^n
\sum_{k=\tau+1}^N b_n[k]\Biggr).
\endmultline
$$
For any $\sigma$ and $\tau$ this expression viewed as a function in
$(\zeta,\xi)$ ranging over the circles $|\zeta|=1$, $|\xi|=1\pm\epsilon$,
remains uniformly bounded away from $0$ and $\infty$ as $t\to 1$ if
$0<\epsilon\ll r$. Indeed, the boundedness of the second factor is obvious,
while for the first factor we have
$$
1-t^n\ge nrt^n,\qquad \bigl|\xi^{\pm n}-(\xi/|\xi|)^{\pm n}\bigr|= |(1\pm
\epsilon)^{\pm n}-1|\le n\epsilon(1+\epsilon)^n
$$
and hence
$$
\multline \sum_{n\ge 1}\left|\Re \,\frac {\overline{A}_n{\zeta^{n}(1-\xi^{-n})}
-A_n\zeta^{-n}(1-\xi^n)}{1-t^n} \right|=\sum_{n\ge 1}\left|\Re \,\frac
{\overline{A}_n{(\zeta/\xi)^{n}} -A_n(\zeta/\xi)^{-n}}{1-t^n} \right|\\ \le
2\epsilon r^{-1}\sum_{n\ge 1} t^{-n}(1+\epsilon)^n|A_n|
\endmultline
$$
is bounded as $t\to 1$ because $\epsilon r^{-1}\to 0$, and $|A_n|=O(R^n)$ for a
fixed $R<1$ as $n\to\infty$. Similar arguments imply that the derivative $\frac
d{d\xi}\tfrac{F(\sigma,\zeta)}{F(\tau,\zeta\xi^{-1})}$ multiplied by $r$ is
uniformly bounded.

Let us now look more carefully at the ratio in question when $\xi$ is close to
the real line. Assume that $|\arg \xi|\le \alpha$ with $r\ll\alpha\ll r^{\frac
12}$. Then, using the estimates
$$
\gathered |\xi^{\pm n}-1\mp n\ln\xi|\le (n|\ln\xi|)^2
\max\{1,|\xi|^n\}\le \const\cdot n^2\alpha^2
(1+\epsilon)^n\ll  n^2r(1+\epsilon)^n,\\
\left|(1-t^n)^{-1}-\frac
1{nr}\right|=\frac{|t^n-1-nr|}{(1-t^n)\,nr}\le
\frac{n^2r^2}{nrt^n\,nr}=t^{-n},
\endgathered
$$
and the fact that $|A_n|=O(R^n)$ as $n\to\infty$, we obtain
$$
\multline \frac{F(\sigma,\zeta)}{F(\tau,\zeta\xi^{-1})}= (1+F_0(\xi))\cdot
\exp\Bigl( r^{-1}\ln \xi\sum_{n\ge
1}\bigl(\overline{A}_n{\zeta^{n}} +A_n\zeta^{-n}\bigr)\Bigr)\\
\times\,\exp\sum_{n\ge
1}\Biggl({\zeta^{-n}}\sum_{k=1}^{\sigma}a_n[k]-{\zeta^{-n}}\sum_{k=1}^{\tau}a_n[k]
-\zeta^n\sum_{k=\sigma+1}^N b_n[k]+\zeta^n \sum_{k=\tau+1}^N b_n[k]\Biggr),
\endmultline
$$
where $F_0(\xi)$ uniformly converges to 0 as $r\to 0$. Differentiating this
relation with respect to $\xi$ and using the boundedness of
$\tfrac{F(\sigma,\zeta)}{F(\tau,\zeta\xi^{-1})}$ and $r\frac
d{d\xi}\tfrac{F(\sigma,\zeta)}{F(\tau,\zeta\xi^{-1})}$ as well as the fact that
$$ \Biggl|\exp \Bigl(r^{-1}\ln \xi\sum_{n\ge 1}\bigl(\overline{A}_n{\zeta^{n}}
+A_n\zeta^{-n}\bigr)\Bigr)\Biggr|=\exp \Bigl(r^{-1}\ln(1\pm\epsilon) \sum_{n\ge
1}\bigl(\overline{A}_n{\zeta^{n}} +A_n\zeta^{-n}\bigr)\Bigr)\to 0
$$
as $r\to 0$, we see that $rF'_0(\xi)$ is uniformly bounded in both $\xi$
ranging over the arch $|\arg\xi|\le\alpha$, $|\xi|=1\pm\epsilon$, and $\zeta$
ranging over the unit circle.

We can now proceed to evaluating the asymptotics of the integral
over $\xi$. Set
$$
G_\pm(\zeta)=-\frac{\prod_{n\ge 1}(1-t^n)^3}{\theta_3(z;t) \,(2\pi
i)}\oint\limits_{|\xi|=1\pm\epsilon}
\dfrac{F(\sigma,\zeta)}{F(\tau,\zeta\xi^{-1})}
\,\frac{\theta_3(z\xi;t)}{\theta_3(-\xi t^{-\frac 12};t)}\,
\frac{d\xi}{\xi^{y+\frac12}}\,.
$$
According to the hypothesis of Theorem 3.1, we assume that
$y=\gamma r^{-1}(1+o(1))$. Note that $|\xi^{y+\frac 12}|\to 1$ as
$r\to 0$, because $\epsilon r\to 0$.

Let us split the circle $|\xi|=1\pm\epsilon$ into three parts:
$|\arg\xi|\ge \delta$ with $\delta>0$ taken from Proposition
3.2(i); $\alpha\le |\arg\xi|<\delta$ with $r\ll \alpha\ll
r^{\frac12}$; and $|\arg\xi|<\alpha$.

The integral over the first part tends to zero (uniformly in
$\zeta$ on the unit circle) because $f(\xi)$ of Proposition 3.1
converges to 0 and the factors
$\tfrac{F(\sigma,\zeta)}{F(\tau,\zeta\xi^{-1})}$ and
$\xi^{-y-\frac 12}$ remain bounded.

The integral over the second part, by virtue of the boundedness of
the same factors and the asymptotics of $f(\xi)$ from Proposition
3.1(i), is bounded by (recall the notation $v=\frac 1{2\pi
i}\ln\xi$)
$$
\frac\const r\int_{v+\frac{i\ln(1\pm\epsilon)}{2\pi
}\in(-\delta,-\frac\alpha{2\pi}]\cup [\frac\alpha{2\pi},\,\delta)}
\left|\frac{ e^{\frac{2\pi i \ln(z)v}{r}}}{e^{\frac{-2\pi^2 v}r}-
e^{\frac{2\pi^2 v}r}}\right|dv.
$$
Introducing a new variable $\hat v=2\pi^2 r^{-1} v$, we can
rewrite this expression as
$$
\const \int_{\hat v+\frac{ir\ln(1\pm\epsilon)}{4\pi^2
}\in(-\frac{2\pi^2 \delta }r,-\frac{\pi\alpha}{r}]\cup [\frac{\pi
\alpha}r,\,\frac{2\pi^2 \delta }r)} \left|\frac{ e^{\frac{i
\ln(z)\hat v}{\pi }}}{e^{{-\hat v}}- e^{\hat v}}\right|d\hat v,
$$
which tends to 0 as $r\to 0$ as long as $|\arg z|<\pi$, $\epsilon
r\to 0$, and $\alpha/r\to\infty$.

The only nonzero contribution comes from the third part. Using the
variable $\hat v =2\pi^2 r^{-1} v=-{i\pi}r^{-1} \ln\xi$, the
estimates for the ratio
$\tfrac{F(\sigma,\zeta)}{F(\tau,\zeta\xi^{-1})}$ obtained above,
and Proposition 3.1(i), we can write the integral over the third
part $|\arg\xi|<\alpha$ in the form (note that $d\xi=\frac {i
r\xi}{\pi} d\hat v$)
$$
\multline\exp\sum_{n\ge
1}\Biggl({\zeta^{-n}}\sum_{k=1}^{\sigma}a_n[k]-{\zeta^{-n}}\sum_{k=1}^{\tau}a_n[k]
-\zeta^n\sum_{k=\sigma+1}^N b_n[k]+\zeta^n \sum_{k=\tau+1}^N
b_n[k]\Biggr)\\ \times \frac{1}{i\pi} \int\limits_{\hat
v+\frac{ir\ln(1\pm\epsilon)}{4\pi^2
}\in(-\frac{\pi\alpha}{r},\frac{\pi \alpha}r)}\frac{ e^{\frac
{i\hat v}\pi\left( \ln(z)-\gamma+\sum_{n\ge
1}\bigl(\overline{A}_n{\zeta^{n}}
+A_n\zeta^{-n}\bigr)\right)}}{e^{{\hat v}}- e^{-\hat v}}\,(1+G_0(
\hat v))\,d\hat v,
\endmultline
$$
where the function $G_0(\hat v)$ uniformly converges to 0 as $t\to 1$, and the
derivative $G_0'(\hat v)$ remains uniformly bounded. The statement of Theorem
3.1(i) is a corollary of this formula, the above estimates, and the following
lemma.

\proclaim{Lemma 3.3} For any $a\in\C$ such that $|\Im a|<1$ the
following limit relation holds:
$$
\lim_{\varepsilon\to 0+}\int_{-\infty\pm i\varepsilon}^{+\infty\pm
i\varepsilon} \frac{e^{iax}}{e^x-e^{-x}}\,(1+g(\varepsilon,x))dx=\frac{\mp\,
i\pi}{1+e^{\pm \pi a}}\,,
$$
where we assume that the function $g(\varepsilon,x)$ on $\R \pm i\varepsilon$
uniformly tends to 0 as $\varepsilon\to 0+$, while its derivative $\frac d{dx}
g(\varepsilon,x)$ remains uniformly bounded. The convergence is uniform in $a$
varying over any compact subset of $\C$ with $\max |\Im a|<1$.
\endproclaim
\demo{Proof} The results with $x\in\R+i\varepsilon$ and $x\in\R-i\varepsilon$
are obtained from each other by the change of sign of the integration variable.
Hence, we may assume that $x$ varies over $\R-i\varepsilon$.

 Let us handle
the term with $g(\varepsilon,x)$ first. Split the integral into
three parts:
$$
\int_{-\infty-i\varepsilon}^{\infty-i\varepsilon}
\frac{e^{iax}g(x,\varepsilon)dx}{e^x-e^{-x}}=\int_{-\infty-i\epsilon}^
{-x_0-i\varepsilon}+\int_{-x_0-i\epsilon}^
{x_0-i\varepsilon}+\int_{x_0-i\epsilon}^{+\infty-i\varepsilon}
$$
for some $x_0=x_0(\varepsilon)\in\R$ that will be chosen later. Denote
$M(\varepsilon)=\sup_{x\in\R-i\varepsilon}|g(x,\varepsilon)|$. The third
integral can be estimated as follows:
$$
\multline \left|\int_{x_0-i\epsilon}^{+\infty-i\varepsilon}
\frac{e^{iax}g(x,\varepsilon)dx}{e^x-e^{-x}}\right|\le M(\varepsilon)
\int_{x_0-i\epsilon}^{+\infty-i\varepsilon} \frac{e^{\Re(ax)}dx}{e^{\Re
x}-e^{-\Re x}}\\=e^{\varepsilon\Re a}M(\varepsilon) \int_{x_0}^{+\infty}
\frac{e^{\Im(a)y}dy}{e^y-e^{-y}}\le e^{\varepsilon\Re a}M(\varepsilon)
\int_{x_0}^{+\infty} \frac{dy}{e^{(1-|\Im
a|)y}-1}\\
= \frac{e^{\varepsilon\Re a}M(\varepsilon)}{1-|\Im a|}\ln\left(\frac{e^{(1-|\Im
a|)x_0}}{e^{(1-|\Im a|)x_0}-1}\right) \,.
\endmultline
$$
If $x_0(\varepsilon)\to 0$ as $\varepsilon\to 0$, then the last expression is
bounded by $\const\cdot M(\varepsilon)|\ln x_0(\varepsilon)|$, which tends to 0
if $x_0(\varepsilon)>e^{-\frac1{M(\varepsilon)}}$. A similar estimate holds for
the first integral. As for the second integral, using the notation
$M'(\varepsilon)=\sup\limits_{(-x_0-i\varepsilon,x_0-i\varepsilon)} |\frac
d{dx}(e^{iax} g(x,\varepsilon))|$, we have
$$
\multline \left|\int_{-x_0-i\varepsilon}^{x_0-i\varepsilon}
\frac{e^{iax}g(x,\varepsilon)dx}{e^x-e^{-x}}\right| \le
M'(\varepsilon)\int_{-x_0}^{x_0} \frac{|y|dy}{|e^{y-i\varepsilon}-
e^{-y+i\varepsilon}|}+\left|\int_{-x_0}^{x_0}
\frac{g(-i\varepsilon,\varepsilon)dy}{e^{y-i\varepsilon}-
e^{-y+i\varepsilon}}\right|\\
\le \const
M'(\varepsilon)x_0+\frac{M(\varepsilon)}2\left|\ln(e^{x_0-i\varepsilon}-1)-\ln
(e^{-x_0+i\varepsilon}-1)-\ln\frac{e^{x_0-i\varepsilon}+1}{e^{-x_0+i\varepsilon}+1}
\right|\\ \le \const (M'(\varepsilon)x_0+M(\varepsilon)),
\endmultline
$$
which tends to 0 if $x(\varepsilon)\to 0$ as $\varepsilon\to 0$. Thus, we
proved that the term with $g(x,\varepsilon)$ tends to 0 as $\varepsilon \to 0$.

The remaining term gives (we use $\hat y=\cosh(y)/\sin\varepsilon$
below)
$$
\multline \int_{-\infty-i\varepsilon}^{\infty-i\varepsilon}
\frac{e^{iax}dx}{e^x-e^{-x}}=
e^{a\varepsilon}\int_{-\infty}^\infty\frac{e^{iay}dy}{e^{y-i\varepsilon}-
e^{-y+i\varepsilon}}\\
=ie^{a\varepsilon}\sin\varepsilon\int_{-\infty}^\infty\frac{
\cos(ay)(e^y+e^{-y})dy}{(e^{y}- e^{-y})^2+4\sin^2\varepsilon}+
ie^{a\varepsilon}\cos\varepsilon\int_{-\infty}^\infty\frac{
\sin(ay)(e^y-e^{-y})dy}{(e^{y}- e^{-y})^2+4\sin^2\varepsilon}\\
\to \frac{i}2\lim_{\varepsilon\to 0} \int_{-\infty}^\infty
\frac{\cos(a\cdot\operatorname{arccosh}(\hat
y\sin\varepsilon))d\hat y}{\hat y^2+1}+
\int_{-\infty}^\infty\frac{\sin(ay)dy}{e^y-e^{-y}}\\
=\frac{i\pi}2+\frac{i\pi}2\,\frac{e^{\pi a}-1}{e^{\pi a}+1}
=\frac{i\pi}{1+e^{-\pi a}}\,. \qed
\endmultline
$$
\enddemo

Lemma 3.3 implies that as $t\to 1$,
$$
\multline G_\pm(\zeta)\to\pm\exp\sum_{n\ge 1}\Biggl({\zeta^{-n}}
\sum_{k=1}^{\sigma}a_n[k]-{\zeta^{-n}}\sum_{k=1}^{\tau}a_n[k]
-\zeta^n\sum_{k=\sigma+1}^N b_n[k]+\zeta^n \sum_{k=\tau+1}^N b_n[k]\Biggr)
\\ \times\frac{1}{1+e^{\mp\left( \ln(z)-\gamma+\sum_{n\ge
1}\bigl(\overline{A}_n{\zeta^{n}} +A_n\zeta^{-n}\bigr)\right)}}\,,
\endmultline
$$
where the convergence is uniform in $\zeta$ varying over the unit circle, and
in $z$ varying in a compact subset of $\C$ not touching the negative real
semiaxis. Since
$$
K(\sigma,x;\tau,y)=\cases \frac 1{2\pi i}\oint_{|\zeta|=1}
G_+(\zeta)\,\frac{d\zeta}{\zeta^{x-y+1}}\,,&\sigma\le \tau,\\
\frac 1{2\pi i}\oint_{|\zeta|=1}
G_-(\zeta)\,\frac{d\zeta}{\zeta^{x-y+1}}\,,&\sigma<\tau,
\endcases
$$
this completes the proof of part (i) of Theorem 3.1.

In order to prove Theorem 3.1(ii) we will use Proposition 2.1, which (along
with Theorem 3.2) implies
$$
\rho_n(\tau_1,x_1;\dots;\tau_n,x_n)=\frac 1{2\pi i}\oint_{|z|=1}
\det[K(\tau_i,x_i;\tau_j,x_j)]_{i,j=1}^n \,\theta_3(z;t)\,\frac {dz}{z}\,.
$$
Using the familiar notation $w=\frac 1{2\pi i}\ln z=\tfrac
1{2\pi}\arg z\in[-\tfrac 12,\tfrac12]$ and applying the imaginary
Jacobi transform, we obtain
$$
\theta_3(z;t)=\left(\tfrac{2\pi} r\right)^\frac12
e^{-\frac{2\pi^2w^2}{r}}\theta_3\Bigl(e^{\frac{4\pi^2 w}r}\,;
e^{-\frac{4\pi^2}r}\Bigr)
$$
(this formula was already used in the proof of Proposition 3.2 above). The
product formula for $\theta_3\bigl(e^{\frac{4\pi^2 w}r}\,;
e^{-\frac{4\pi^2}r}\bigr)$ implies that this function of $w$ remains uniformly
bounded on $[-\frac 12,\frac 12]$ as $r\to 0$, and it uniformly converges to 1
on any open interval inside $[-\frac 12,\frac 12]$.

Let us split the domain of integration over $z$ in the formula for
$\rho_n(\tau_1,x_1;\dots;\tau_n,x_n)$ above into two arches:
$|\arg z|\le c<\pi$ and $|\arg z|>c$, where $c>0$ in a arbitrary
constant $<\pi$.

The integral over the first arch, by Theorem 3.1(i), is equal to
$$
\left(\tfrac{2\pi} r\right)^\frac12\int_{-\frac c{2\pi}}^{\frac c{2\pi}}
e^{-\frac{2\pi^2w^2}{r}} \det\bigl[ \K_{\tau_i,\tau_j}^{(e^{2\pi i
w},\gamma)}(x_i-x_j)\bigr]_{i,j=1}^n\, dw + o(1)
$$
as $r\to 0$. Since $\left(\tfrac{2\pi} r\right)^\frac12
e^{-\frac{2\pi^2w^2}{r}}$ converges to the delta-function at $w=0$ as $r\to 0$,
the integral above converges to
$\det\bigl[\K_{\tau_i,\tau_j}^{(1,\gamma)}(x_i-x_j)\bigr]_{i,j=1}^n$.

As for the integral over the arch $|\arg z|>c$, we use Proposition
3.2(ii) and the boundedness of the ratio
${F(\sigma,\zeta)}/{F(\tau,\zeta\xi^{-1})}$ proved earlier to see
that the absolute value of the integrand is bounded by
$\const\cdot\epsilon^{-1}r^{-\frac 12} e^{-\frac {c^2}{2r}}$,
which converges to 0 as $r\to 0$ as long as $\ln \epsilon^{-1}\ll
r$. Since this does not contradict our previous assumption that
$\epsilon\ll r$, we may ignore the integral over the second arch
in the limit $r\to 0$. This completes the proof of Theorem
3.1.\qed
\enddemo

\example{Example 3.4} Consider the periodic Schur process with $N=1$ and
$$
a_n[1]=b_n[1]=\cases \vartheta,& n=1,\\ 0,&n>1,\endcases
$$
with an arbitrary $\vartheta>0$. This is equivalent to considering a
probability measure on pairs $(\lambda\supset\mu)$ of partitions given by
$$
\p\{\lambda,\mu\}=\left(e^{\frac{c_1}{1-c_1c_2}}\,\prod_{n\ge 1}\frac
1{1-(c_1c_2)^n}\right)^{-1}\frac{\dim^2(\lambda/\mu)\,c_1^{|\lambda|}c_2^{|\mu|}}
{(|\lambda|-|\mu|)!^2}\,,
$$
where $c_1=t\vartheta^2$, $c_2=\vartheta^{-2}$, and $\dim(\lambda/\mu)$ is the
number of standard Young tableaux of shape $\lambda/\mu$. Here we used the
relation
$$
s_{\lambda/\mu}(a[1])=s_{\lambda/\mu}(b[1])=\frac{\dim(\lambda/\mu)
\vartheta^{|\lambda|-|\mu|}}{(|\lambda|-|\mu|)!}
$$
and the formula for the partition function from Proposition 1.1.

Theorem 3.1(ii) yields the following limit result for the local
correlation functions of this measure: If
$x_1(t),\dots,x_n(t)\in\Z'$ are such that
$$
\lim_{t\to 1}\left(\ln t^{-1}\cdot x_i(t)\right)=\gamma,\qquad i=1,\dots,n,
$$
and the pairwise distances $x_i-x_j$ do not depend on $t$, then
$$
\lim_{t\to
1}\p\bigl\{\{x_1(t),\dots,x_n(t)\}\subset\{\lambda_i-i+\tfrac12\}_{i\ge
1}\bigr\}=\det[\Cal K^{(\gamma)}(x_i-x_j)]_{i,j=1}^n,
$$
where
$$
\Cal K^{(\gamma)}(d)=\frac 1{2\pi i}\oint\limits_{|\zeta=1|} \frac{1}
{1+e^{\gamma-2\vartheta\Re\zeta}}\,\frac{d\zeta}{\zeta^{d+1}}\,.
$$
The formula for the limit correlation functions for $\mu$ is exactly the same.

In the limit $\vartheta\to 0$ the measure on $\lambda$'s approaches the uniform
measure $\p\{\lambda\}=\const\cdot t^{|\lambda|}$ on partitions, and,
correspondingly, $\Cal K^{(\gamma)}(d)\to \delta_{0d}\cdot (1+e^\gamma)^{-1}$,
cf. Comment 4 to Theorem 3.1.

As $\vartheta\to\infty$, the measure on $\lambda$'s looks more like the
poissonized Plancherel measure on partitions with weights of the form
$\p\{\lambda\}\sim {\dim^2\lambda \cdot {\const}^{|\lambda|}}/{|\lambda|!^2}$,
because as $c_2=\vartheta^{-2}\to 0$, the partition $\mu$ tends to be small.
When $\gamma$ is of the same order as $\vartheta$, the local correlation kernel
converges to the discrete sine kernel,
$$
\lim_{\Sb \vartheta\to\infty\\ \gamma/\vartheta\to \Theta\endSb}\Cal
K^{(\gamma)}(d) = \frac 1{2\pi
i}\oint\limits_{|\zeta=1|}\frac{\id_{2|\Re\zeta|>\Theta}\,
d\zeta}{\zeta^{d+1}}=\cases 0,&\Theta\ge 2,\\
\dfrac{\sin\bigl(\arccos(\Theta/2)\cdot d\bigr)}{\pi d},&|\Theta|<2,\\
\delta_{0d},&\Theta\le -2,\endcases
$$
which is exactly the bulk scaling limit of the correlation kernel
of the poissonized Plancherel measure with large poissonization
parameter, cf. \cite{BOO, Theorem 3}.
%Note, however, that this empirical argument does not
%provide a rigorous computation of the local correlation functions for our
%initial measure in the double limit $t\to 1-$, $\vartheta\to\infty$, or, in
%other terms, $c_1\to\infty$, $c_2\to 0$, $c_1c_2\to 1-$. For a rigorous
%argument one needs to go back to the formula of Theorem 2.2 and analyze it
%under these limiting assumptions.
\endexample

\head 4. Extensions of the discrete sine kernel
\endhead

The goal of this section is to construct an infinite-dimensional family of
determinantal point processes on $\Z\times\Z$ such that the restrictions of
their correlation kernels to copies of $\Z$ obtained by fixing the first
coordinate coincide with the discrete sine kernel.

Two such extensions of the discrete sine kernel have been
constructed previously: in \cite{OR1} a kernel called the {\it
incomplete beta kernel\/} was obtained in the bulk limit of the
large uniformly distributed plane partitions, and in \cite{BO2}
another extension was obtained in the bulk limit of the Markov
chains on partitions preserving the Plancherel measure
(equivalently, in the bulk of the multi-layer polynuclear growth
process with droplet initial conditions). Both these extensions
are included in the family that will be constructed below.

We start with some necessary generalities.

Set $\X=\Z\times\Z$ and let $2^{\X}$ be the set of all subsets of $\X$. The
first coordinate of points in $\X$ will be viewed as ``time'' while the second
coordinate will be viewed as ``space''. Pick any exhaustion of $\X$ by a
sequence of increasing finite sets:
$$
\b_1\subset \b_2\subset\b_3\subset\dots\subset \X,\qquad |\b_n|<\infty,\qquad
\bigcup_{n\ge 1} \b_n=\X.
$$
For instance, one can take $\b_n$'s to be growing boxes $\b_n=\{(x,t)\,:\,
|t|\le n,\, |x|\le n\}$. Clearly, $2^\X=\lim\limits_{\longleftarrow} 2^{\b_n}$,
and we equip $2^\X$ with the projective limit topology. In other words, a
sequence $\{X_n\}\subset 2^\X$ converges to $X\subset \X$ if and only if for
any $m\ge 1$ there exists $M=M(m)$ such that for $n>M$ we have $X_n\cap
\b_m=X\cap \b_m$. It is easy to see that the topology does not depend on the
choice of $\b_n$'s. Since projective limits of compact topological spaces are
compact, the space $2^\X$ is compact with respect to this topology.

A random point process on $\X$ is, by definition, a Borel
probability measure on $2^\X$.\footnote{In the conventional
terminology, this is actually the definition of a {\it simple} or
{\it multiplicity free\/} random point process. Since all our
processes are multiplicity free, the definition above is general
enough for our purposes.} One way to construct random point
processes on $\X$ is to provide a countable set $X_1,X_2,\dots$ of
subsets of $\X$ and the set of positive weights $w_1,w_2,\dots$;
$\p\{X_i\}=w_i$, satisfying the normalization condition
$\sum_{i\ge 1}w_i=1$.

The correlation functions $\rho_n$ of a random point process on $\X$ are
probabilities for random subsets of $\X$ to contain a given finite set:
$$
\rho_n(t_1,x_1;\dots;t_n,x_n)=\p\{X\subset \X\mid
\{(\tau_1,x_1),\dots,(\tau_n,x_n)\}\subset X\}.
$$

For a finite set $\b\subset \X$, the values of the correlation functions on
$\b$ define the projection of the measure on $2^\X$ to $2^\b$ uniquely. Indeed,
the two sets of $2^{|\b|}-1$ numbers giving for any nontrivial subset $\b_0$ of
$\b$ the probability that the intersection of the random set $X\subset\X$ with
$\b$ {\it contains\/} $\b_0$ or {\it coincides with\/} $\b_0$ are related by a
nondegenerate linear transformation obtained from the inclusion-exclusion
principle. One set consists of the values of the correlation functions while
the other one consists of the weights of the subsets with respect to the
projected measure on $2^\b$. For example, for $\b=\{b_1,b_2\}$ we have
$$
\gathered \p\{X\in 2^\X\mid X\cap \b=\{b_1\}\}=\rho_1(b_1)-\rho_2(b_1,b_2),\\
\p\{X\in 2^\X\mid X\cap \b=\{b_2\}\}=\rho_1(b_2)-\rho_2(b_1,b_2),\\ \p\{X\in
2^\X\mid X\cap \b=\{b_1,b_2\}\}=\rho_2(b_1,b_2).
\endgathered
$$

\proclaim{Lemma 4.1} Let $\P_1,\P_2,\dots$ be a sequence of random
point processes on $\Z^2$, and assume that, as $m\to\infty$, all
correlation functions of $\P_m$ converge pointwise: For any
$n=1,2,\dots$ and any $(\tau_i,x_i)\in\Z^2$, $i=1,\dots,n$, we
have
$$
\lim_{m\to\infty} \rho_n(\tau_1,x_1;\dots;\tau_n,x_n\mid
\P_m)=:r_n(\tau_1,x_1;\dots;\tau_n,x_n).
$$
Then there exists a unique random point process $\P$ on $\Z^2$ such that
$$
r_n(\tau_1,x_1;\dots;\tau_n,x_n)=\rho_n(\tau_1,x_1;\dots;\tau_n,x_n\mid \P).
$$
\endproclaim
\demo{Proof} The uniqueness of $\P$ was demonstrated above. The same argument
shows that for any finite set $\b$, the set of possible values of the
correlation functions of a probability measure on $2^\b$ is described by a
finite list of linear inequalities. For instance, in the case $\b=\{b_1,b_2\}$
these inequalities are
$$
\gathered
\p\{\varnothing\}=1-\rho_1(b_1)-\rho_1(b_2)+\rho_2(b_1,b_2)\ge 0,\\
\p\{b_1\}=\rho_1(b_1)-\rho_2(b_1,b_2)\ge 0,\quad
\p\{b_2\}=\rho_1(b_2)-\rho_2(b_1,b_2)\ge 0,\\
\p\{b_1,b_2\}=\rho_2(b_1,b_2)\ge 0.
\endgathered
$$
Clearly, these inequalities are preserved under limit transitions. Hence, for
any of the finite sets $\b_m$ used above there exists a unique probability
measure $\P_m$ on $2^{\b_m}$ such that its correlation functions coincide with
the restrictions of the limit functions $r_n$ to $\b_m$. Moreover, these
measures are consistent: for $m_1>m_2$ the projection of the measure $\P_{m_1}$
on $2^{\b_{m_1}}$ to $2^{\b_{m_2}}$ coincides with $\P_{m_2}$. (Indeed, both
these measures have the same correlation functions $r_n|_{\b_{m_2}}$.) Since
$2^\X=\lim\limits_{\longleftarrow} 2^{\b_m}$, we can take
$\P=\lim\limits_{\longleftarrow} \P_m$. \qed
\enddemo

Observe that all that was said above is not specific to the set $\X=\Z^2$ and
applies equally well to any discrete countable set $\X$. In particular, we will
consider random point processes on $\X_N:=\{1,\dots,N\}\times\Z$ below. Since
$\X_N\subset \X$, a random point process on $\X_N$ may also be viewed as a
process on $\X$ whose point configurations are contained in $\X_N$ almost
surely.

Let us now return to the periodic Schur process. So far we have
not discussed the positivity of the weights which we used to
define the process. One way (but not the only way, see e.g. \S8
below) to guarantee the nonnegativity of the weights
$\w(\lambda,\mu)$ introduced in the beginning of \S1 is to take
$t>0$ and to demand that all the specializations of the skew Schur
functions are nonnegative. The following classical result is
useful.

Due to the Jacobi-Trudi formula
$$
s_{\lambda/\mu}=\det[h_{\lambda_i-i-\mu_j+j}]_{i,j=1}^r,\qquad r\ge
\max\{\ell(\lambda),\ell(\mu)\},
$$
see \cite{Macd, \S{I} (5.4)}, we need to guarantee the
nonnegativity of the determinants in the right-hand side (here
$h_n$'s are the complete homogeneous symmetric functions). Recall
that a sequence $\{c_n\}_{n=0}^\infty$ is called {\it totally
positive} if all minors of the matrix $[c_{i-j}]_{i,j\ge 0}$ are
nonnegative. Here all $c_{-k}$ for $k>0$ are assumed to be equal
to zero. We will only consider totally positive sequences with
$c_0=1$; clearly, multiplication of all members of a sequence by
the same positive number does not affect total positivity.

The following statement was independently proved by
Aissen-Edrei-Schoenberg-Whitney in 1951 \cite{AESW}, \cite{Ed},
and by Thoma in 1964 \cite{Th}. An excellent exposition of deep
relations of this result to representation theory of the infinite
symmetric group can be found in Kerov's book \cite{Ke}.

\proclaim{Theorem 4.2} A sequence $\{c_n\}_{n=0}^\infty$, $c_0=1$, is totally
positive if and only if its generating series has the form
$$
\sum_{n=0}^\infty c_nu^n=e^{\gamma u}\frac{\prod_{i\ge 1}
(1+\beta_iu)}{\prod_{i\ge 1}(1-\alpha_iu)}=:TP_{\alpha,\beta,\gamma}(u)
$$
for certain nonnegative parameters $\{\alpha_i\}$, $\{\beta_i\}$ and $\gamma$
such that $\sum_i(\alpha_i+\beta_i)<\infty$.

Equivalently, an algebra homomorphism $\phi:\Lambda\to \C$ takes nonnegative
values on all skew Schur functions if and only if the sequence
$\{\phi(h_n)\}_{n\ge 0}$ is totally positive, that is,
$$
\sum_{n\ge 0}\phi(h_n)u^n=TP_{\alpha,\beta,\gamma}(u)
$$
for a suitable choice of parameters $(\alpha,\beta,\gamma)$.
\endproclaim

The function $TP_{\alpha,\beta,\gamma}(u)$ is meromorphic in $u$, and it is
holomorphic and nonzero in a small enough neighborhood of the origin. In order
to satisfy the convergence conditions, see Remark 2.3, we will actually have to
use only the specializations for which $TP_{\alpha,\beta,\gamma}(u)$ is
holomorphic and nonzero in a disc of radius greater than 1.

Let us call a specialization $a$ of the algebra of symmetric functions {\it
admissible\/} if $H(a;u)=TP_{\alpha,\beta,\gamma}(u)$ for certain parameters
$(\alpha,\beta,\gamma)$ as above, and $\alpha_i,\beta_i<1$ for all $i\ge 1$.

We say that a specialization $a$ is a union of specializations $a^{(1)},
a^{(2)}, \dots, a^{(m)}$ if
$$
H(a;u)=H\bigl(a^{(1)};u\bigr)H\bigl(a^{(2)};u\bigr)\cdots
H\bigl(a^{(m)};u\bigr).
$$
Note that unions of admissible specializations are admissible. We will use the
notation $a=\biguplus_{k=1}^m a^{(k)}$.

\proclaim{Proposition 4.3} For any $N\ge 1$ and any admissible specializations
$a[1],b[1]$, $\dots$, $a[N],b[N]$ of $\Lambda$, denote
$$
a(\sigma,\tau]=\biguplus_{k=\sigma+1}^\tau a{[k]},\quad
b(\sigma,\tau]=\biguplus_{k=\sigma+1}^\tau b{[k]},\qquad a=a(0,N],\quad
b=b(0,N].
$$
Assume that $a=b$. Then for any $C>0$ there exists a unique random point
process on $\{1,\dots,N\}\times\Z$ with determinantal correlation functions
$$
\rho_n(\tau_1,x_1;\dots;\tau_n,x_n)=\det[\Cal
K_{\tau_i,\tau_j}(x_i-x_j)]_{i,j=1}^n
$$
{\rm (}here $n\ge 1$, $\tau_i\in\{1,\dots,N\}$, $x_i\in\Z$ are arbitrary{\rm
)}, and the correlation kernel
$$
\Cal K_{\sigma,\tau}(d)=\cases\dfrac 1{2\pi i}\oint\limits_{|z|=1} \dfrac
{\left(H(a(\sigma,\tau];\zeta^{-1})H(b(\sigma,\tau];\zeta)\right)^{-1}}{1+ C
\left(H(a;\zeta^{-1})H(b;\zeta)\right)^{-1}}\,
\dfrac{d\zeta}{\zeta^{d+1}},&\sigma\le
\tau,\\
-\dfrac 1{2\pi i}\oint\limits_{|z|=1} \dfrac
{H(a(\tau,\sigma];\zeta^{-1})H(b(\tau,\sigma];\zeta)}{1+ C^{-1}
H(a;\zeta^{-1})H(b;\zeta)}\, \dfrac{d\zeta}{\zeta^{d+1}},&\sigma>\tau.
\endcases
$$
\endproclaim
\demo{Proof} Consider the periodic Schur process with $t\in (0,1)$ and
specializations $a[1],b[1],\dots,a[N],b[N]$ as in the hypothesis. By Theorem
4.2, all weights of this measure on $\Y^N$ are nonnegative. Mapping $\Y^N$ to
$2^{\{1,\dots,N\}\times\Z'}$ via
$$
(\lambda^{(1)},\dots,\lambda^{(N)})\mapsto
\Bigl(\bigl\{\lambda^{(1)}_i-i+\tfrac12\bigr\}_{i\ge
1},\dots,\bigl\{\lambda^{(N)}_i-i+\tfrac12\bigr\}_{i\ge 1}\Bigr)
$$
yields a random point process on $\{1,\dots,N\}\times\Z'$ whose correlation
functions were computed in \S2. Shifting the space variable $x\mapsto [\ln
C/\ln t^{-1}]+x+\frac 12$ and using Theorem 3.1(ii) we see that the correlation
functions of our process converge, as $t\to 1-$, to those given in the
hypothesis above. Lemma 4.1 completes the proof. \qed
\enddemo

The main result of this section is the following statement.

\proclaim{Theorem 4.4} For any doubly infinite sequences
$\{a[k],b[k]\}_{k\in\Z}$ of admissible specializations of $\Lambda$ and any
$c\in (0,\pi)$, there exists a determinantal\footnote{i.e., with determinantal
correlation functions.} point process on $\Z\times\Z$ with the correlation
kernel
$$
\Cal K_{\sigma,\tau}(x-y)=\cases\dfrac 1{2\pi i}\int\limits_{e^{-ic}}^{e^{ic}}
{\left(H(a(\sigma,\tau];\zeta^{-1})H(b(\sigma,\tau];\zeta)\right)^{-1}}\,
\dfrac{d\zeta}{\zeta^{x-y+1}},&\sigma\le
\tau,\\
-\dfrac 1{2\pi i}\int\limits_{e^{ic}}^{e^{-ic}}
{H(a(\tau,\sigma];\zeta^{-1})H(b(\tau,\sigma];\zeta)}\,
\dfrac{d\zeta}{\zeta^{x-y+1}},&\sigma>\tau,
\endcases
$$
where the integrals are taken over positively oriented arches of the unit
circle.
\endproclaim
\demo{Comments} 1. The equal time values of the kernel above are exactly those
of the discrete sine kernel on $\Z$; for any $\tau\in\Z$
$$
\Cal K_{\tau,\tau}(x-y)=\dfrac 1{2\pi
i}\int\limits_{e^{-ic}}^{e^{ic}}
\dfrac{d\zeta}{\zeta^{x-y+1}}=\frac{e^{ic(x-y)}-e^{-ic(x-y)}}{2\pi
i (x-y)}=\frac{\sin(c(x-y))}{\pi (x-y)}\,.
$$
In this sense the kernels $\Cal K_{\sigma,\tau}(x-y)$ provide
extensions of the discrete sine kernel.

2. The point processes in questions are clearly invariant with respect to the
shifts of the space coordinate. If one wants the processes to be invariant with
respect to the time shifts as well, one has to take all specializations $a[k]$
to be the same, and all specializations $b[k]$ to be the same.

3. Using Lemma 4.1 we can slightly relax the conditions on the specializations
$a[k]$, $b[k]$ by allowing the parameters $\alpha_i$ to be equal to 1. (Recall
that they are required to be strictly less than 1 by the definition of
admissible specializations.) Taking trivial specializations $a[k]$ and choosing
$b[k]$ so that $H(b[k];u)=(1-u)^{-1}$, we obtain the kernel
$$
\Cal K_{\sigma,\tau}(x-y)=\cases\dfrac 1{2\pi i}\int\limits_{e^{-ic}}^{e^{ic}}
(1-\zeta)^{\tau-\sigma}\, \dfrac{d\zeta}{\zeta^{x-y+1}},&\sigma\le
\tau,\\
-\dfrac 1{2\pi i}\int\limits_{e^{ic}}^{e^{-ic}} (1-\zeta)^{\tau-\sigma}\,
\dfrac{d\zeta}{\zeta^{x-y+1}},&\sigma>\tau,
\endcases
$$
which is the incomplete beta kernel of \cite{OR1}. We will also
see this kernel arising in the bulk limit of the cylindric
partitions in the subsequent sections.

4. The choice $H(a[k];u)=H(b[k];u)=e^{\const_k\cdot u}$ yields the
extension of the discrete sine kernel obtained in \cite{BO2,
Theorem 4.2}.
\enddemo

\demo{Proof} It suffices to prove the statement for the processes on
$\{1,\dots,N\}\times \Z$ and finite sequences $\{a[k],b[k]\}_{k=1}^N$ of
admissible specializations. Indeed, then one can just embed such processes in
$2^{\Z\times\Z}$ and take the limit as $N\to\infty$ using Lemma 4.1.

In the finite $N$ case we use Proposition 4.3 to construct a process on
$\{1,\dots,2N+1\}\times\Z$ with extra $2(N+1)$ specializations
$a[N+1],b[N+1],\dots,a[2N+1],b[2N+1]$ defined by
$$
\gathered a[N+k]=b[k],\quad b[N+k]=a[k], \qquad k=1,\dots,N,\\
H(a[2N+1];u)=H(b[2N+1];u)=e^{M u},\quad M>0.
\endgathered
$$
Then
$$
\gathered
\biguplus_{k=1}^{2N+1}a[k]=a=b=\biguplus_{k=1}^{2N+1}b[k],\\
H(a[2N+1];\zeta^{-1})H(b[2N+1];\zeta)=e^{2M\Re \zeta}.
\endgathered
$$
Choosing the constant $C$ in Proposition 4.3 to be $C=e^{2M\cos c}$, for
$\zeta$ on the unit circle we obtain
$$
\gathered \lim_{M\to \infty}\frac 1{1+C
\bigl(H(a;\zeta^{-1})H(b;\zeta)\bigr)^{-1}}=\cases 1,&\Re \zeta>\cos c,\\
0,&\Re\zeta<\cos c,\endcases
\\
\lim_{M\to \infty}\frac 1{1+C^{-1} H(a;\zeta^{-1})H(b;\zeta)}=\cases 0,
&\Re \zeta>\cos c,\\
1,&\Re\zeta<\cos c.\endcases
\endgathered
$$
Hence, as $M\to\infty$ we have the convergence of the correlation functions of
our $M$-dependent processes restricted to $\{1,\dots,N\}\times \Z$ to the
needed values, and Lemma 4.1 completes the proof.\qed
\enddemo

\head 5. Cylindric partitions
\endhead

Cylindric partitions were first introduced by I.~Gessel and C.~Krattenthaler in
[GK]. We will initially follow their paper in our exposition.

Let $\lambda$ and $\mu$ be two partitions. Assume that $\lambda\supset\mu$ and
denote the length (=number of nonzero parts) of $\lambda$ by $l$. A {\it plane
partition of shape} $\lambda/\mu$ is a planar array $\pi $ of integers of the
form
$$
\matrix &&&&\pi_{1,\mu_1+1}&\cdots&\pi_{1,\lambda_l}&\cdots&\cdots&\pi_{1,\lambda_1}\\
&&\pi_{2,\mu_2+1}&\cdots&\pi_{2,\mu_1+1}&\cdots&\pi_{2,\lambda_l}&\cdots&\pi_{2,\lambda_2}&&\\
&\cdots&\cdots&\cdots&\cdots&\cdots&\cdots&\cdots&&\\
\pi_{l,\mu_l+1}&\cdots&\pi_{l,\mu_2+1}&\cdots&\pi_{l,\mu_1+1}&\cdots&\pi_{l,\lambda_l}&&&
\endmatrix
$$
such that the rows and columns are weekly decreasing: $\pi_{i,j}\ge
\pi_{i,j+1}$ and $\pi_{i,j}\ge \pi_{i+1,j}$.

{\it A cylindric partition of shape} $\lambda/\mu/d$ can be viewed as a plane
partition with an additional relation between the first and the last rows. This
relation depends on an integral parameter $d$ and it has the form $\pi_{l,j}\ge
\pi_{1,j+d}$ for all $j$. In other words, a cylindric partition has to remain a
plane partition when the last row shifted by $d$ to the right is placed on top
of the first row. Here is an example of a cylindric partition of shape
$(8,6,3)/(3,1)/4$. On the second picture the shifted last row (in bold) is
placed on top of the first one.
$$
\matrix &&&7&5&2&1\\
&10&10&6&5&1&1\\
11&9&1
\endmatrix
\qquad\qquad \matrix &&&&\text{\bf 11}&\text{\bf 9}&\text{\bf 1}\\
&&&7&5&2&1\\
&10&10&6&5&1&1\\
11&9&1
\endmatrix
$$
Note that $d$ has to be greater or equal to $\mu_1$.

The {\it norm} $|\pi|$ of the cylindric partition $\pi$ is the sum of its
elements. In the example above $|\pi|=68$.

For our purposes it is more convenient to parameterize cylindric partitions
differently. Namely, let us read them along the lines with fixed content $j-i$
(these lines are parallel to the diagonal which has content zero). On each line
we observe an ordinary partition, and partitions on neighboring lines are
different by adding or removing a horizontal strip.\footnote{Recall that a
Young diagram (equivalently, a partition) $\nu$ can be obtained from another
Young diagram $\kappa$ by adding a horizontal strip (notation $\kappa\prec\nu$
or $\nu\succ\kappa$) if and only if $\nu_i\ge\kappa_i\ge \nu_{i+1}$ for all
$i\ge 1$.} The necessary number of fixed content lines to be taken into account
is equal to $N:=l+d$, and the content can be considered as an element of the
cyclic group $\Z/N\Z$ of order $N$.

The example above leads to the following (periodic with period $N=7$) sequence
of partitions:
$$
\ldots\prec(11,2,1)\succ (9,1)\prec (10,1)\succ (10)\succ
(6)\prec(7,5)\succ(5,1)\prec(11,2,1)\succ\ldots
$$

The ``$\prec$'' or ``$\succ$'' relation of the neighboring partitions depends
on the boundary of the Young diagram $\mu$: horizontal edges correspond to
``$\succ$'' and vertical edges correspond to ``$\prec$''.

It is impossible to reconstruct $\lambda$ and  $\mu$ from such a sequence of
partitions: in the example above if we remove the last row from the second
picture then the resulting partition sequence will be the same while we would
be looking at a cylindric partition of type $(6,6,6)/(3,2)/4$. However, $d$ and
$l$ remain invariant --- $d$ is equal to the total number of ``$\prec$'' and
$l$ is equal to the total number of ``$\succ$'' in a period.

Let us encode the sequence of ``$\succ$'' and ``$\prec$'' by assigning to a
cylindric partition a periodic sequence of 1's and $-1$'s; 1's correspond to
$\succ$'s and $-1$'s correspond to $\prec$'s. Thus, the example above produces
$(\dots,-1,1,-1,1,1,-1,1,-1,1,\dots)$ (this sequence is periodic with period
$N=7$). We will call this sequence the {\it profile\/} of the corresponding
cylindric partition. Clearly, the profile depends only on $\mu,d$, and $l$, but
not on $\lambda$.

We say that the profile is {\it marked} if there is a marked $-1$ in each
period, and the distance between any marked $-1$'s is a multiple of the period.
In other words, if we consider profiles as maps from $\Z/N\Z$ to $\{-1,1\}$
then marking corresponds to choosing a distinguished element of $\Z/N\Z$ in the
preimage of $-1$.

In order to associate a {\it marked} profile to any cylindric partition we will
mark the $-1$'s corresponding to the relation of partitions in lines with
content $d$ and $-l+1$. Thus, in our example the marked $-1$'s corresponds to
$(5,1)\subset (11,9,1)$ and the marked profile is
$$(\dots,-1^*,1,-1,1,1,-1,1,-1^*,1,\dots).$$
Marked profiles are in one-to-one correspondence with triples $(\mu,d,l)$ with
$d\ge \mu_1$ and $l\ge \ell(\mu)$.

Let us associate to any marked profile two periodic sequences $A[k]$, $B[k]$
with period $N$ and elements $0$ or $1$ such that the difference $A[k]-B[k]$
gives the element of the profile with distance $k$ from the marked $-1$. In the
example above
$$\align
A&=(\dots,0^*,1,0,1,1,0,1,0^*,1,\dots), \\
B&=(\dots,1^*,0,1,0,0,1,0,1^*,0,\dots),
\endalign
$$
where we marked $A[0\operatorname{mod}N]$ and $B[0\operatorname{mod}N]$. As was
mentioned above, we have
$$
d=\sum_{k=1}^N A[k],\qquad l=\sum_{k=1}^N B[k].
$$

For any integer $m$ let $m(N)$ be the smallest positive integer such that
$$
m\equiv m(N)\mod N.
$$
For instance, $1(N)=1$ and $-1(N)=N-1$.

In what follows we will be interested in probability measures on
cylindric partitions $\pi$ with fixed $\mu,d$, and $l$
(equivalently, fixed marked profile), whose weights are
proportional to $s^{|\pi|}$, where $s\in (0,1)$ is a parameter. We
will call these probability measures {\it uniform\/} as the
weights of the cylindric partitions with the same norm are equal.
The next statement provides the partition function for the weights
$s^{|\pi|}$.

\proclaim{Proposition 5.1} For any partition $\mu$ and integers $d\ge\mu_1$,
$\l\ge \ell(\mu)$, the following identity holds:
$$
\sum_{\{\lambda:\ell(\lambda)\le l\}}\ \sum_{\{\pi: \pi \text{ has shape
}\lambda/\mu/d\}}s^{|\pi|}=\prod_{n\ge 1}\frac 1{1-s^{nN}}\prod_{\Sb
p\in \overline{1,N}:\,A[p\,]=1\\
q\in \overline{1,N}:\,B[q\,]=1\endSb} \frac1{1-s^{(p-q)(N)+(n-1)N}}\,,
$$
where $N=d+l$, and $A[k]$ and $B[k]$ are sequences of 0's and 1's associated to
$(\mu,d,l)$ as described above.
\endproclaim
\demo{Comments} 1. It is convenient to view the function $(p-q)(N)$ as an array
on the $N\times N$ torus with rows and columns parameterized by $p$ and $q$:
$$
\bmatrix N&N-1&N-2&\dots&3&2&1\\
1&N&N-1&\dots &4&3&2\\
2&1&N&\dots&5&4&3\\
\dots&\dots&\dots&\dots&\dots&\dots&\dots\\
\dots&\dots&\dots&\dots&\dots&\dots&\dots\\
N-2&N-3&N-4&\dots&1&N&N-1\\
N-1&N-2&N-3&\dots&2&1&N
\endbmatrix
$$
In the formula above we choose rows $p$ such that $A[p]=1$ and columns $q$ such
that $B[q]=1$. Note that the sets of rows and columns chosen in such a way form
a disjoint splitting of $\{1,\dots,N\}$ into two sets; the first set contains
$d$ elements and the second set contains $l$ elements. Hence, the total number
of factors in the second product equals $dl$.

2. The right-hand side of the formula for the partition function depends on the
profile of $(\mu,d,l)$, but not on the marked profile. This is in agreement
with rotational symmetry of the problem. The sum in the left-hand side of the
formula may be viewed as the sum over all cylindric partitions with fixed
marked profile.

3. If $d=0$ then $\mu$ must be empty, and all rows of the
corresponding cylindric partitions become identical. Accordingly,
the second product in the right-hand side disappears, and we
recover the formula for the partition function of the weights
$s^{(\#\text{ of boxes})\cdot N}$ on ordinary partitions.

4. If $d,l\gg 1$ and $A=(1,1,\dots,1,0,0\dots,0^*)$,
$B=(0,0,\dots,0,1,1\dots,1^*)$ then $\mu=d^l$ and the
corresponding cylindric partitions of small enough norm are in
one-to-one correspondence with similar plane partitions.
Accordingly, looking at the top right corner of the matrix in
Comment 1 above, we see that the partition function for small
powers of $s$ looks like that for the plane partitions:
$\prod_{n\ge 1}(1-s^n)^{-n}$. In the limit $d,l\to\infty$ one
recovers the celebrated MacMahon's formula for the partition
function of the weights $s^{\text{norm}}$ on plane partitions.

5. The case of largest rotational symmetry corresponds to
$A=(1,0,1,0,\dots,1,0^*)$ and $B=(0,1,0,1,\dots,0,1^*)$. Then $l=d$, $N=2d$,
and $\mu=(d,d-1,\dots,2,1)$ is the staircase partition. In this case the
formula for the partition function slightly simplifies to give
$$
\prod_{n\ge 1}\frac
1{1-s^{2nd}}\prod_{m=1}^{d}\frac1{(1-s^{2m-1+2(n-1)d})^d}\,.
$$
\enddemo
\demo{Proof} It is well known that the skew Schur function $s_{\lambda/\mu}$
specialized at a single variable $x$ is nonzero if and only if
$\lambda\succ\mu$, in which case $s_{\lambda/\mu}(x)=x^{|\lambda|-|\mu|}$.

This observation immediately implies that cylindric partitions in question are
in one-to-one correspondence with trajectories of the periodic Schur process
with period $N$ determined by the specializations (see \S1 for notations)
$$
\wt a_n[k]=\tfrac1n\,A[k],\quad \wt b_n[k]=\tfrac1n\,B[k],\qquad k=1,\dots,N,
$$
or, in different terms,
$$
H(a[k];u)=(1-s^{k}u)^{-A[k]},\quad H(b[k];u)=(1-s^{-k}u)^{-B[k]},\qquad
k=1,\dots,N.
$$
 Namely, we read the cylindric partitions along the lines
with fixed content, and the resulting ordinary partitions form a periodic Schur
process. Recall that following the marking rule we denote the partition coming
from the line with content $(m \mod N)$ by $\lambda^{(m+l)}$. In particular,
the marked place of the profile corresponds to the relation $\lambda^{(0)}\prec
\lambda^{(1)}$.

The formula of Proposition 5.1 is then a direct corollary of Proposition 1.1.
\qed
\enddemo

The one-to-one correspondence of the cylindric partitions with a given profile
and trajectories of the periodic Schur process explained in the proof above
also implies that we can use the results of \S2 to obtain the correlation
functions of the uniform measure on cylindric partitions.

We define the correlation functions of the uniform measure on cylindric
partitions with a given profile as the dynamical correlation functions of the
corresponding periodic Schur process.

\proclaim{Proposition 5.2} The correlation functions of the uniform measure on
cylindric partitions with a given profile $\{A[k],B[k]\}_{k=1}^N$ have the
following form: For any $n\ge1$, $\tau_1,\dots,\tau_n\in\{1,\dots,N\}$,
$x_1,\dots,x_n\in\Z'$,
$$
\rho_n(\tau_1,x_1;\dots;\tau_n,x_n)=\oint_{|z|=1}
\det[K(\tau_i,x_i;\tau_j,x_j)]_{i,j=1}^n\frac{\theta_3(z;t)\,dz}z\,,
$$
where $t=s^N$,
$$
K(\sigma,x;\tau,y)=-\frac{\prod_{n\ge 1}(1-t^n)^3}{\theta_3(z;t) \,(2\pi
i)^2}\oint_\zeta\oint_\eta \dfrac{F(\sigma,\zeta)}{F(\tau,\eta^{-1})}
\,\frac{\theta_3(z\zeta\eta;t)}{\theta_3(-\zeta\eta\, t^{-\frac 12};t)}\,
\frac{d\zeta d\eta}{\zeta^{x+\frac12}\eta^{y+\frac12}}\,,
$$
the integrals are taken over circles centered at the origin with radii
satisfying
$$
s^{\sigma+1}<|\zeta|<s^\sigma,\quad s^{-\tau}<|\eta|<s^{-\tau-1},\qquad
|\zeta\eta|>1\ \text{  if  }\ \sigma=\tau,
$$
and
$$
\multline F(\tau,\zeta)=\exp\sum_{n\ge 1}\frac
1{n(1-s^{nN})}\Biggl(\sum_{k=1}^{\tau}B[k](s^{-k}\zeta )^n+ \sum_{k=\tau+1}^N
B[k]
(s^{N-k}\zeta)^n\\
-\sum_{k=1}^{\tau}A[k](s^{k+N}/\zeta)^n -\sum_{k=\tau+1}^N A[k](s^k/\zeta)^n
\Biggr).
\endmultline
$$
\endproclaim
\demo{Proof} The statement follows from Theorem 2.2, Proposition 2.1, and
Remark 2.4. The only new effect that we have to be careful about is the choice
of the contours of integration, because the function $F(\tau,\zeta)$ is
correctly defined by the formula above only in the ring
$s^{\tau+1}<|\zeta|<s^\tau$, which explains our conditions on the contours.
\qed
\enddemo
\example{Remark 5.3} The integration over $z$ in the formula for the
correlation functions above can be carried out explicitly to produce a
multivariate integral formula similar to that in Corollary 2.8.
\endexample
\example{Remark 5.4} The function $F(\tau,\zeta)$ can be written in terms of
infinite products:
$$
F(\tau,\zeta)=\frac{\prod\limits_{k\in\overline{1,\tau}:A[k]=1}\bigl(
s^{k+N}/\zeta;t\bigr)_\infty
\prod\limits_{k\in\overline{\tau+1,N}:A[k]=1}\bigl(s^{k}/\zeta;t\bigr)_\infty}{
\prod\limits_{k\in\overline{1,\tau}:B[k]=1}\bigl( s^{-k}\zeta;t\bigr)_\infty
\prod\limits_{k\in\overline{\tau+1,N}:B[k]=1}\bigl(s^{N-k}\zeta;t\bigr)_\infty}\,,
$$
where, as usual, $(a;t)_\infty=\prod_{n\ge 1}(1-at^n)$. This formula also
provides the analytic continuation of $F(\tau,\zeta)$ from the ring
$s^{\tau+1}<|\zeta|<s^\tau$.
\endexample

\head 6. Bulk of large cylindric partitions with finite or slowly growing
period
\endhead

In this section we compute the local limit of the correlation functions of the
uniform measure on cylindric partitions with fixed profile
$\{A[k],B[k]\}_{k=1}^N$ in two cases: when $N$ is fixed and the profile does
not depend on the small parameter $r=\ln t^{-1}=N\ln s^{-1}$, and when $N$ is
growing in such a way that $r$ still tends to 0 as $s\to 1-$.

In both cases an important role is played by the parameter
$$
\varkappa=\frac ld=\frac{\sum_{k=1}^N B[k]}{\sum_{k=1}^N A[k]}
$$
which we call the {\it slope\/} of the profile. In the case of growing $N$ in
order to have a limit we will need to assume that the slope tends to a limiting
value strictly between 0 and $\infty$.

We need to do some preliminary work before stating the results.

For any $\varkappa>0$, let $\Gamma_\varkappa$ be the contour in the complex
plane defined by
$$
\multline
\Gamma_\varkappa=\left\{z=1-\tfrac{\sin(1+\varkappa)\phi}{\sin\varkappa\phi}\,e^{i\phi}
\mid \phi\in(-\tfrac\pi{1+\varkappa},\tfrac\pi{1+\varkappa})\right\}
\\=\left\{-\frac{\sin\frac\varphi{1+\varkappa}}{\sin\frac\varphi{1+\varkappa^{-1}}}
\,e^{i\varphi}\mid \varphi\in(-\pi,\pi)\right\}.
\endmultline
$$

This is a piecewise smooth closed curve which has a corner-like singularity at
$z=1$ as $\phi\to\pm \tfrac\pi{1+\varkappa}$ or $\varphi\to\pm\pi$. The origin
is located inside $\Gamma_{\varkappa}$. The curve is also symmetric with
respect to the real axis, it intersects the negative semiaxis at the point
$z=-\varkappa^{-1}$, and $\Gamma_{\varkappa_1}$ is located inside
$\Gamma_{\varkappa_2}$ if $\varkappa_1>\varkappa_2$. Note that $\Gamma_1$ is
the unit circle, which is the only case when the curve is smooth at $z=1$.
Whenever we use $\Gamma_\varkappa$ as the integration contour, we will assume
that it is oriented counterclockwise.

Our interest in the family $\{\Gamma_{\varkappa}\}$ is explained by the
following statement. We use the principal branch of the logarithm function
below.

\proclaim{Lemma 6.1} For any $\varkappa>0$, set
$f_\varkappa(z)=\varkappa\ln(1-z)+\ln(1-z^{-1})$. This is a holomorphic
function on $\C\setminus \R_{\ge 0}$, and
$$
\{z\in\C\setminus \R_{\ge 0}\mid f_\varkappa(z)\in\R\}=\R_{<0}\cup
\Gamma_{\varkappa}\,.
$$
Further, $\Im f_{\varkappa}(z)<0$ if and only if $z$ is located inside
$\Gamma_{\varkappa}$ and below the real axis or outside $\Gamma_{\varkappa}$
and above the real axis.

On the curve $\Gamma_{\varkappa}$ the function $f_\varkappa$ equals
$$
f_\varkappa\bigl(1-\tfrac{\sin(\varkappa+1)\phi}{\sin\varkappa\phi}\,e^{i\phi}\bigr)
=\varkappa\ln\left(\tfrac{\sin(\varkappa+1)\phi}{\sin\varkappa\phi}\right)+
\ln\left(\tfrac{\sin(\varkappa+1)\phi}{\sin\phi}\right).
$$
This function increases on $(-\frac\pi{\varkappa+1},0)$, decreases on
$(0,\frac\pi{\varkappa+1})$, and its maximal value on $\Gamma_{\varkappa}$
equals
$$
f_\varkappa(-\varkappa^{-1})=\ln(1+\varkappa)+\varkappa\ln(1+\varkappa^{-1}).
$$
\endproclaim
\demo{Proof} We have
$$
\frac
{df_\varkappa(re^{i\varphi})}{dr}=r^{-1}\left(\varkappa-\frac{1+\varkappa}
{1-re^{i\varphi}}
\right)=r^{-1}\left(\varkappa-\frac{(1+\varkappa)(1-re^{-i\varphi})}
{|1-re^{i\varphi}|^2} \right),
$$
whence on any ray $\{z=re^{i\varphi}\mid r>0\}$, the function $\Im
f_{\varkappa}(z)$ is strictly decreasing as a function of $r$ for
$\varphi\in(0,\pi)$, and it is strictly increasing for $\varphi\in(-\pi,0)$.
Thus, on any such ray there cannot be more than one point $z$ such that
$f_\varkappa(z)\in\R$.

On the other hand, for
$z=1-\frac{\sin(\varkappa+1)\phi}{\sin\varkappa\phi}\,e^{i\phi}\in
\Gamma_\varkappa$ we compute
$$
z=-\tfrac{\sin\phi}{\sin\varkappa\phi}\,e^{i(\varkappa+1)\phi},\quad
z^{-1}=-\tfrac{\sin\varkappa\phi}{\sin\phi}\,e^{-i(\varkappa+1)\phi},\quad
1-z^{-1}=\tfrac{\sin(\varkappa+1)\phi}{\sin\phi}\,e^{-i\varkappa\phi}\,,
$$
which implies the first statement of the lemma and the formula for
$f_{\varkappa}$ on the curve. Since
$$
\gathered \left(\tfrac{\sin(\varkappa+1)\phi}{\sin\varkappa\phi}\right)'=
\tfrac{\sin(2\varkappa+1)\phi-(2\varkappa+1)\sin\phi}
{2(\sin\varkappa\phi)^2}<0,\\
\left(\tfrac{\sin(\varkappa+1)\phi}{\sin\phi}\right)'=
\tfrac{\varkappa\sin(2\varkappa+1)\phi-(2\varkappa+1)\sin\varkappa\phi}
{2(\sin\phi)^2}<0,
\endgathered
$$
for $\phi\in(0,\frac\pi{\varkappa+1})$, and $f_\varkappa$ restricted to
$\Gamma_\varkappa$ is an even function of $\phi$, the last claim of the lemma
also follows.\qed
\enddemo

For any profile $\{A[k],B[k]\}$ denote
$$
A(\sigma,\tau]=\sum_{k=\sigma+1}^\tau A[k],\qquad
B(\sigma,\tau]=\sum_{k=\sigma+1}^\tau B[k].
$$
Then $d=A(0,N]$, $l=B(0,N]$, and $\varkappa=B(0,N]/A(0,N]$.

Now we can state the main result of this section.

\proclaim{Theorem 6.2} In the limit $s\to 1-$, the correlation
functions of the uniform measure on cylindric partitions with a
given profile $\{A[k],B[k]\}_{k=1}^N$ have a limit in the
following sense: Choose $x_1(s),\dots,x_n(s)\in\Z+\frac 12$ such
that as $s\to 1-$, $r x_k(s)\to \gamma$ for all $k=1,\dots,n$ and
some $\gamma\in\R$, and all pairwise distances
$x_i-x_j=x_i(s)-x_j(s)$ remain constant. Then for any
$\tau_1,\dots,\tau_n\in\{1,\dots,N\}$
$$
\lim_{s\to 1-}\rho_n(\tau_1,x_1(s);\dots,\tau_n,x_n(s))=\det\bigl[
\K^{(\gamma)}_{\tau_i,\tau_j}(x_i-x_j)\bigr]_{i,j=1}^n,
$$
where the correlation kernel has the form
$$
\K^{(\gamma)}_{\sigma,\tau}(x-y)=\cases\frac 1{2\pi
i}\oint_{\Gamma_{\varkappa}}
\dfrac{(1-\zeta)^{B(\sigma,\tau]}(1-\zeta^{-1})^{A(\sigma,\tau]}}
{1+e^\gamma(1-\zeta)^l(1-\zeta^{-1})^d}\,\dfrac{d\zeta}{\zeta^{x-y+1}}\,,&\sigma\le
\tau,\\
-\frac 1{2\pi i}\oint_{\Gamma_\varkappa}
\dfrac{(1-\zeta)^{-B(\tau,\sigma]}(1-\zeta^{-1})^{-A(\tau,\sigma]}}
{1+e^{-\gamma}(1-\zeta)^{-l}(1-\zeta^{-1})^{-d}}\,\dfrac{d\zeta}{\zeta^{x-y+1}}\,,
&\sigma>\tau.
\endcases
$$
\endproclaim
\demo{Comments} 1. By Lemma 6.1, the expression $(1-\zeta)^l(1-\zeta^{-1})^d$
is nonnegative on $\Gamma_\varkappa$, which shows that the kernel
$\K^{(\gamma)}_{\sigma,\tau}$ is correctly defined for $\sigma\le \tau$. On the
other hand, for $\sigma>\tau$ we can rewrite the formula for the kernel as
$$
\K^{(\gamma)}_{\sigma,\tau}(x-y)=-\frac 1{2\pi i}\oint_{\Gamma_\varkappa}
\dfrac{e^\gamma(1-\zeta)^{l-B(\tau,\sigma]}(1-\zeta^{-1})^{d-A(\tau,\sigma]}}
{1+e^{\gamma}(1-\zeta)^{l}(1-\zeta^{1})^{d}}\,\dfrac{d\zeta}{\zeta^{x-y+1}}\,,
$$
and this integral also always makes sense.

2. The global limit density function
$$
\rho(\gamma)=\K^{(\gamma)}_{\tau,\tau}(0)=\frac 1{2\pi
i}\oint_{\Gamma_{\varkappa}} \dfrac{1}
{1+e^\gamma(1-\zeta)^l(1-\zeta^{-1})^d}\,\dfrac{d\zeta}{\zeta}
$$
does not depend on $\tau$, which means that it is invariant with respect to
rotations of the cylindric partitions. Note also that $\rho(\gamma)$ depends on
the fixed profile only through $d=A(0,N]$ and $l=B(0,N]$, or, in other words,
through the period $N=d+l$ and the slope $\varkappa=l/d$.

3. The simplest nontrivial (i.e., not coinciding with the uniform measure on
ordinary partitions) example is $d=l=1$. Then $\Gamma_{\varkappa}$ is the unit
circle, and it is not hard to evaluate the global limit density function and
the corresponding hypothetical limit shape (see Comment 3 after Theorem 3.1 for
explanations) explicitly:
$$
\rho(\gamma)=\frac 1{\sqrt{1+4e^\gamma}}\,,\qquad
v(u)=u+4\operatorname{arctanh}(\sqrt{1+4e^{u}})+\tfrac {i\pi}2\,.
$$
For larger values of $d$ and $l$ explicit integration is also possible but as
$d$ and $l$ grows it becomes increasingly tedious.

4. A formal application of Theorem 3.1 to the uniform measure on
cylindric partitions produces the correct integrand in the formula
for the limit kernel, but the integration contour is different
unless $\varkappa=1$. For $\varkappa\ne 1$ using the unit circle
(as Theorem 3.1 suggests) instead of the curve
$\Gamma_{\varkappa}$ may lead to a wrong answer!
\enddemo

\proclaim{Theorem 6.3} In the limit
$$
s\to 1-,\quad d\to\infty,\quad l\to\infty,\qquad r=-N\ln s\to 0,\quad l/d\to
\varkappa>0,
$$ the correlation functions of the uniform measure on cylindric
partitions with a given profile $\{A[k],B[k]\}_{k=1}^N$ have a limit in the
following sense:
 Choose
$x_1(s),\dots,x_n(s)\in\Z+\frac 12$ such that as $s\to 1-$, $r x_k(s)/N\to
\gamma$ for all $k=1,\dots,n$ and some $\gamma\in\R$, and all pairwise
distances $x_i-x_j=x_i(s)-x_j(s)$ remain constant. Also choose the time moments
$\tau_1(s),\dots,\tau_n(s)$ such that as $s\to 1$ the pairwise distances
$\tau_i(s)-\tau_j(s)$ remain uniformly bounded. Then
$$
\rho_n(\tau_1(s),x_1(s);\dots,\tau_n(s),x_n(s))=\det\bigl[
 \K^{(\gamma)}_{\tau_i(s),\tau_j(s)}(x_i-x_j)\bigr]_{i,j=1}^n +o(1),
$$
where the estimate is uniform over any set of profiles with the slope uniformly
convergent to $\varkappa$ as $s\to 1$, and the correlation kernel has the
following form:
$$
\K^{(\gamma)}_{\tau,\tau}(x-y)=\delta_{x-y,0}\ \text{  if  }\
\gamma\le\gamma_0(\varkappa):=-\frac{\ln(1+\varkappa)}{1+\varkappa}
-\frac{\ln(1+\varkappa^{-1})}{1+\varkappa^{-1}}\,,
$$
and for $\gamma>\gamma_0(\varkappa)$, using the notation
$\zeta(\phi)=1-\tfrac{\sin(1+\varkappa)\phi}{\sin\varkappa\phi}\,e^{i\phi}$, we
have
$$
\K^{(\gamma)}_{\sigma,\tau}(x-y)=\cases\frac 1{2\pi
i}\int\limits_{\zeta(\phi)}^{\overline{\zeta(\phi)}}
(1-\zeta)^{B(\sigma,\tau]}(1-\zeta^{-1})^{A(\sigma,\tau]}
\,\dfrac{d\zeta}{\zeta^{x-y+1}}\,,&\sigma\le
\tau,\\
-\frac 1{2\pi i}\int\limits_{\overline{\zeta(\phi)}}^{\zeta(\phi)}
(1-\zeta)^{-B(\tau,\sigma]}{(1-\zeta^{-1})^{-A(\tau,\sigma]} }
\,\dfrac{d\zeta}{\zeta^{x-y+1}}\,, &\sigma>\tau,
\endcases
$$
where both integration contours leave the origin on their left sides, and the
constant $\phi=\phi(\gamma)\in(0,\frac{\pi}{1+\varkappa})$ from the limits of
the integrals above is uniquely determined by the relation
$$
-(1+\varkappa)\gamma=\varkappa\ln\left(\tfrac{\sin(\varkappa+1)\phi}{\sin\varkappa\phi}\right)+
\ln\left(\tfrac{\sin(\varkappa+1)\phi}{\sin\phi}\right)\,.
$$
\endproclaim

\demo{Comments} 1. The existence of the unique $\phi(\gamma)$ for
$\gamma>\gamma_0$, as well as nonexistence of $\phi$ satisfying
the relation above for $\gamma<\gamma_0$, follows directly from
Lemma 6.1. Using this lemma it is also easy to see that the
statement of Theorem 6.3 can be formally deduced from Theorem 6.2.
Indeed,
$$
\multline \ln \bigl(e^\gamma
(1-\zeta)^l(1-\zeta^{-1})^d\bigr)=N\bigl(\tfrac\gamma N
+\tfrac{\varkappa}{1+\varkappa}\,\ln(1-\zeta)+\tfrac
1{1+\varkappa}\,\ln(1-\zeta^{-1})\bigr)\\
=\tfrac{N}{1+\varkappa}\bigl( \tfrac{(1+\varkappa)}N\,\gamma
+f_\varkappa(\zeta)\bigr).
\endmultline
$$
Thus, for $\zeta\in\Gamma_\varkappa$, $N\to\infty$, $\gamma\to\infty$, and
$\gamma/N\to\widehat\gamma$ we have
$$
\lim_{N\to\infty} e^\gamma (1-\zeta)^l(1-\zeta^{-1})^d=\cases 0,
&f_\varkappa(\zeta)<-(1+\varkappa)\widehat\gamma,\\
+\infty,& f_\varkappa(\zeta)>-(1+\varkappa)\widehat\gamma.\endcases
$$
Substituting this limit relation in the formulas of Theorem 6.2 we recover the
formulas of Theorem 6.3.

2. The global limit density function for $\gamma<\gamma_0(\varkappa)$ is
identically equal to 1, while for $\gamma>\gamma_0(\varkappa)$ it is equal to
$$
\rho(\gamma)=\K^{(\gamma)}_{\tau,\tau}=\frac 1{2\pi
i}\int\limits_{\zeta(\phi)}^{\overline{\zeta(\phi)}} \,\dfrac{d\zeta}{\zeta}=
-\frac {\arg\zeta(\phi)}\pi=1-\frac{(1+\varkappa)\phi}\pi\,.
$$
Note that this expression is independent of $\tau$, which reflects
the rotational invariance of the limit density. As
$\gamma\to\gamma_0$ we have $\phi(\gamma)\to 0$ and
$\rho(\gamma)\to 1$. On the other hand, as $\gamma\to+\infty$ we
have $\phi(\gamma)\to\frac\pi{1+\varkappa}$ and $\rho(\gamma)\to
0$.

The fact that $\rho(\gamma)\equiv 1$ for $\gamma<\gamma_0$ means
that the random Young diagrams in question have the lower edge ---
the event of having columns of length substantially greater than
$N\gamma_0(\varkappa)/r$ has vanishing probability.

3.  For $\varkappa=1$ the formulas simplify, and we obtain
$$
\gathered
\varkappa\ln\left(\tfrac{\sin(\varkappa+1)\phi}{\sin\varkappa\phi}\right)+
\ln\left(\tfrac{\sin(\varkappa+1)\phi}{\sin\phi}\right)=2\ln(2\cos
\phi),\quad\gamma_0=-\ln 2,\\
\rho(\gamma)=\cases \frac 2\pi \arcsin(\frac{e^{-\gamma}}2),& \gamma>-\ln 2,\\
1,&\gamma<-\ln 2.
\endcases
\endgathered
$$
Since $\frac 2\pi \arcsin(\frac{e^{-\gamma}}2)=\frac
1\pi\arccos(1-\frac{e^{-2\gamma}}2)$, in this case the global
limit density coincides (up to rescaling of $\gamma$ by 2) with
that for the largest section of random plane partitions with
uniform measure, see \S3.1.10 of \cite{OR1} with $\tau=0$.

4. Equal time values of the limit correlation kernel equal
$$
\K_{\tau,\tau}(x-y)=\frac 1{2\pi i}
\int\limits_{\zeta(\phi)}^{\overline{\zeta(\phi)}}
\dfrac{d\zeta}{\zeta^{x-y+1}}=\left(\tfrac{\sin\varkappa\phi}
{\sin\phi}\right)^{x-y}\frac{\sin\bigl((\pi-(1+\varkappa)\phi)(x-y)\bigr)}
{\pi(x-y)}\,.
$$
The prefactor $\left(\tfrac{\sin\varkappa\phi} {\sin\phi}\right)^{x-y}$ can be
ignored as it cancels out in the determinantal expression for the correlation
functions, and we obtain the discrete sine kernel. Hence, the full kernel
$\K_{\sigma,\tau}(x-y)$ can be viewed as an extension of the discrete sine
kernel. Let us show that this extension belongs to the family of extensions
constructed in Theorem 4.4 above. For $\sigma\le \tau$ we have
$$
\multline \K_{\sigma,\tau}(x-y)=(-1)^{A(\sigma,\tau]}\frac 1{2\pi
i}\int\limits_{\zeta(\phi)}^{\overline{\zeta(\phi)}}
\,\dfrac{(1-\zeta)^{\tau-\sigma}d\zeta}{\zeta^{x-y+1+A(\sigma,\tau]}}\\=
(-1)^{A(\sigma,\tau]}\left(\tfrac{\sin\varkappa\phi}
{\sin\phi}\right)^{x-y+A(\sigma,\tau]} \frac 1{2\pi
i}\int\limits_{e^{i((1+\varkappa)\phi-\pi)}}^{e^{i(\pi -(1+\varkappa)\phi)}}
\,\dfrac{(1-\frac{\sin\phi}{\sin\varkappa\phi}\,\zeta)^{\tau-\sigma}d\zeta}{\zeta^{x-y+1+A(\sigma,\tau]}}
\endmultline
$$
while for $\sigma>\tau$
$$
\K_{\sigma,\tau}(x-y)=-(-1)^{A(\tau,\sigma]}\left(\tfrac{\sin\varkappa\phi}
{\sin\phi}\right)^{x-y-A(\sigma,\tau]}\frac 1{2\pi i}\int\limits_{e^{i(\pi
-(1+\varkappa)\phi)}}^{e^{i((1+\varkappa)\phi-\pi)}}
\dfrac{(1-\frac{\sin\phi}{\sin\varkappa\phi}\,\zeta)^{\sigma-\tau}d\zeta}
{\zeta^{x-y+1-A(\tau,\sigma]}}\,.
$$
Once again, the prefactors can be ignored, and after shifting the space
variable $x$ at time $\tau$ by $x\mapsto x+A(T,\tau]$ for a fixed $T$ and all
$\tau$, we see that the correlation functions are determinants of the kernel
$$
\cases \frac 1{2\pi i}\int\limits_{e^{i((1+\varkappa)\phi-\pi)}}^{e^{i(\pi
-(1+\varkappa)\phi)}}
\,\dfrac{(1-\frac{\sin\phi}{\sin\varkappa\phi}\,\zeta)^{\tau-\sigma}d\zeta}
{\zeta^{x-y+1}}\,,&\qquad\sigma\le \tau,\\
-\frac 1{2\pi i}\int\limits_{e^{i(\pi
-(1+\varkappa)\phi)}}^{e^{i((1+\varkappa)\phi-\pi)}}
\dfrac{(1-\frac{\sin\phi}{\sin\varkappa\phi}\,\zeta)^{\sigma-\tau}d\zeta}
{\zeta^{x-y+1}}\,,&\qquad\sigma>\tau,
\endcases
$$
where the integrals are taken over positively oriented arches of the unit
circle.

For $\varkappa=1$ this is exactly the incomplete beta kernel
arising in the bulk limit of uniform measures on plane partitions,
see \cite{OR1}.

If $\varkappa>1$ then $\frac{\sin\phi}{\sin\varkappa\phi}<1$ and
the kernel coincide with one of the stationary kernels afforded by
Theorem 4.4; one has to take $H(a[k];u)\equiv 1$,
$H(b[k];u)=(1-\frac{\sin\phi}{\sin\varkappa\phi}\,u)^{-1}$ for all
$k\in\Z$.

Finally, if $\varkappa<1$ then
$\frac{\sin\phi}{\sin\varkappa\phi}>1$, and one more shift of the
space variable $x$ at time $\tau$ by $x\mapsto x-T+\tau$ (or,
equivalently, the shift $x\mapsto x-B(T,\tau]$ of the initial
space variable $x$ at time $\tau$) for a fixed $T$ and all $\tau$,
brings our kernel, up to irrelevant prefactors, to the kernel of
Theorem 4.4 with
$H(a[k];u)=(1-\frac{\sin\varkappa\phi}{\sin\phi}\,u)^{-1}$ and
$H(b[k];u)\equiv 1$ for all $k\in\Z$.

Observe that all these kernels are invariant with respect to time shifts. Thus,
the bulk behavior of our cylindric partitions with slowly growing period $N$ is
rotationally invariant and independent of the initial profile except for its
slope, after we perform the above shifting of the space variables. However,
these shifts do depend on the profile.
\enddemo

\demo{Proof of Theorems 6.2 and 6.3} As in the proof of Theorem 3.1, let us
introduce a new integration variable $\xi=\zeta\eta$ and rewrite the kernel
from Proposition 5.2 in the form
$$
K(\sigma,x;\tau,y)=-\frac{\prod_{n\ge 1}(1-t^n)^3}{\theta_3(z;t) \,(2\pi
i)^2}\oint\limits_{\xi}\oint\limits_{\zeta}
\dfrac{F(\sigma,\zeta)}{F(\tau,\zeta\xi^{-1})}
\,\frac{\theta_3(z\xi;t)}{\theta_3(-\xi t^{-\frac 12};t)}\, \frac{d\zeta
d\xi}{\zeta^{x-y+1}\xi^{y+\frac12}}\,,
$$
where, by Remark 5.4,
$$
\multline \dfrac{F(\sigma,\zeta)}{F(\tau,\zeta\xi^{-1})}=
\frac{\prod_{k\in\overline{1,\tau}:B[k]=1}(1-s^{-k}\zeta/\xi)}
{\prod_{k\in\overline{1,\sigma}:B[k]=1}(1-s^{-k}\zeta)}
\frac{\prod_{k\in\overline{\sigma+1,N}:A[k]=1}(1-s^k/\zeta)}
{\prod_{k\in\overline{\tau+1,N}:A[k]=1}(1-s^k\xi/\zeta)}
\\ \times \prod_{k\in\overline{1,N}:B[k]=1}\frac{(s^{-k}t\zeta/\xi;t)_\infty}
{(s^{-k}t\zeta;t)_\infty}\prod_{k\in\overline{1,N}:A[k]=1} \frac
{(s^{k}t/\zeta;t)_\infty}{(s^{k}t\xi/\zeta;t)_\infty}\,.
\endmultline
$$
Here, according to Proposition 5.2, the integrals are taken over circles
centered at the origin with radii satisfying $s^{\sigma+1}<|\zeta|<s^\sigma$
and $s^{\tau+1}<|\zeta|/|\xi|<s^{\tau}$ plus the additional condition $|\xi|>1$
if $\sigma=\tau$.

The proof proceeds along the same lines as that of Theorem 3.1.
The only major obstacle is the unboundedness of the ratio
${F(\sigma,\zeta)}/{F(\tau,\zeta\xi^{-1})}$ on the integration
contours. Recall that in Theorem 3.1 we imposed the condition
$A_k=\overline{B_k}$ which ensured the needed boundedness when
$|\zeta|$ and $|\xi|$ were close to 1. Here we only have such
``self-adjointness'' in the case $\varkappa=1$; for general
$\varkappa>0$ we need to deform the $\zeta$-integration contour.

Clearly, the first factor in the formula for
${F(\sigma,\zeta)}/{F(\tau,\zeta\xi^{-1})}$ above remains bounded
as $s\to 1$ as long as the factors in the denominator do not
approach zero. Furthermore, as $\xi\to 1$ and $s\to 1$, it tends
to $(1-\zeta)^{B(\sigma,\tau]}(1-\zeta^{-1})^{A(\sigma,\tau]}$ for
$\sigma\le \tau$, and to
$(1-\zeta)^{-B(\tau,\sigma]}(1-\zeta^{-1})^{-A(\tau,\sigma]}$ for
$\sigma>\tau$. As we will see, only integration over
infinitesimally close to 1 $\xi$'s yields a nonzero contribution,
and this asymptotics gives the corresponding factors in the final
formula for the correlation functions.

The second factor needs more attention. Recall that the dilogarithm function is
defined in the unit disc by the power series
$$
\di(z)=\sum_{n=1}^\infty \frac{z^n}{n^2}\,.
$$
Its analytic continuation is provided by the integral representation
$$
\di(z)=-\int_0^z\frac{\ln(1-x)\,dx}{x}\,,
$$
which shows that $z=1$ is the branching point of this function. We will
consider $\di(z)$ as the holomorphic function on $\C\setminus(1,+\infty)$
defined by the integral above with the principal branch of $\ln(1-x)$. Note
that the jump of $\di(z)$ across $(1,+\infty)$ is purely imaginary, which means
that $\Re\di(z)$ is a continuous function on $\C$.

Our interest in $\di(\,\cdot\,)$ is explained by the fact that
$$
\ln(x;t)_\infty=\sum_{n=0}^\infty \ln(1-xt^n)=\frac 1{1-t}\sum_{n=0}^\infty
\frac{(xt^n-xt^{n+1})\ln(1-xt^n)}{xt^n}
$$
is, up to the factor $(1-t)^{-1}$, a Riemannian sum for the integral
representation of $\di(x)$. Thus, using the mean value theorem to estimate the
remainder, for $x\in\C\setminus [1,+\infty)$ we obtain
$$
\ln(x;t)_\infty=-(1-t)^{-1}\di(x)+O\bigl(\operatorname{dist}(1,\{rx\mid 0\le
r\le 1\})^{-1}\bigr),\qquad t\to 1,
$$
and the estimate is uniform on compact sets. (The remainder may become large as
$x$ approaches the cut $[1,+\infty)$.)

This asymptotic relation shows that the second factor in the formula for the
ratio ${F(\sigma,\zeta)}/{F(\tau,\zeta\xi^{-1})}$, as $t\to 1$, is approximated
by
$$
\exp\left(-\frac
{l(\di(\zeta/\xi)-\di(\zeta))+d(\di(1/\zeta)-\di(\xi/\zeta))}{1-t} \right).
$$
Note also that, since $\di'(x)=-\ln(1-x)/x$, as $\xi\to 1$ we have
$$
l(\di(\zeta/\xi)-\di(\zeta))+d(\di(1/\zeta)-\di(\xi/\zeta))\sim
\ln\xi(l\ln(1-\zeta)+d\ln(1-\zeta^{-1})),
$$
which coincides, up to multiplication by $d\ln\xi$, with function
$f_\varkappa(\zeta)$ introduced in Lemma 6.1 above.

\proclaim{Lemma 6.4} For any $\varkappa>0$ and $\xi\ne 1$ on the unit circle,
set
$$
f_{\varkappa,\xi}(\zeta)=\varkappa(\di(\zeta/\xi)-\di(\zeta))+\di(1/\zeta)-\di(\xi/\zeta).
$$
Then $f_{\varkappa,\xi}$ is a holomorphic function on $\C\setminus\{\R_+\cup
\xi\R_+\}$, its real part is continuous on $\C$, and
$$
\{\zeta\in\C\mid \Re f_{\varkappa,\xi}(\zeta)=0\}=\sqrt{\xi}\cdot \R\cup
\Gamma_{\varkappa,\xi},
$$
where $\Gamma_{\varkappa,\xi}$ is a piecewise smooth closed curve which
encloses the origin and intersects each ray $e^{i\varphi}\R_+\nsubseteq
\sqrt{\xi}\cdot \R$ at a single point. Further, $\Re
f_{\varkappa,\xi}(\zeta)<0$ if and only if $\zeta$ is inside
$\Gamma_{\varkappa,\xi}$ and to the left of the line $\sqrt{\xi}\cdot \R$, or
$\zeta$ is outside $\Gamma_{\varkappa,\xi}$ and to the right of the line
$\sqrt{\xi}\cdot \R$.

The curve $\Gamma_{\varkappa_1,\xi}$ lies inside
$\Gamma_{\varkappa_2,\xi}$ if $\varkappa_1>\varkappa_2$; the curve
$\Gamma_{1,\xi}$ is the unit circle. Also, as $\xi\to 1$, we have
$\Gamma_{\varkappa,\xi}\to\Gamma_{\varkappa}$ in the sense that
the intersection points of $\Gamma_{\varkappa,\xi}$ with rays
$e^{i\varphi}\R_+$ converge to the corresponding intersection
points of $\Gamma_{\varkappa}$.
\endproclaim
\demo{Proof} The fact that $f_{\varkappa,\xi}$ takes real values on
$\sqrt{\xi}\cdot\R$ follows from the relation
$\overline{\di(x)}=\di(\overline{x})$. On the other hand, on any ray we compute
$$
\Re\,\frac{df_{\varkappa,\xi}(re^{i\varphi})}{dr}=\frac{1+\varkappa}r\,\ln\left|
\frac {1-re^{i\varphi}}{\xi-re^{i\varphi}}\right|.
$$
This expression is negative if $\zeta$ is to the right of the line
$\sqrt{\xi}\cdot \R$, and it is positive if $\zeta$ is to the left of the line
$\sqrt{\xi}\cdot \R$.

Using the asymptotic relation
$$
\lim_{r\to+\infty}\Re\bigl(\di(-re^{i\psi_1})-\di(-re^{i\psi_2})\bigr)=
\tfrac12({\psi_1^2-\psi_2^2}),\qquad -\pi\le \psi_1,\psi_2\le \pi,
$$
we see that the limits of $f_{\varkappa,\xi}(re^{i\varphi})$ as $r\to 0+$ and
$r\to+\infty$ exist, and if $e^{i\varphi}\R_+\nsubseteq \sqrt{\xi}\cdot \R$
then they are nonzero and have different signs. Since the real part of the
derivative along any such ray is sign definite, this implies the existence and
uniqueness of the needed intersection points as well as the inequalities on
$\Re f_{\varkappa,\xi}(\zeta)$.

Similar arguments show that $\frac{df_{\varkappa,\xi}(\zeta)}{d\zeta}$ vanishes
at exactly two points located in $\sqrt{\xi}\cdot\R_+$ and
$\sqrt{\xi}\cdot\R_-$. These are the points where $\Gamma_{\varkappa,\xi}$
intersects the critical line $\sqrt{\xi}\cdot\R$. Thus, the curve is closed,
and a standard implicit function theorem argument implies the smoothness.

Since the limit values $\lim_{r\to 0+}f_{\varkappa,\xi}(re^{i\varphi})$ do not
depend on $\varkappa$, while the absolute value of the derivative
$|\Re\,\frac{df_{\varkappa,\xi}(re^{i\varphi})}{dr}|$ is an increasing (and
linear) function of $\varkappa$, we see that the larger $\varkappa$ the sooner
$f_{\varkappa,\xi}(re^{i\varphi})$ reaches 0 as $r$ increases. This implies the
inclusion property of the curves $\Gamma_{\varkappa,\xi}$. Finally, the
convergence of $\Gamma_{\varkappa,\xi}$ to $\Gamma_\varkappa$ follows from the
asymptotic relation $f_{\varkappa,\xi}(\zeta)\sim \ln\xi f_{\varkappa}(\zeta)$
as $\xi\to 1$, Lemma 6.1, and the fact that $\ln\xi$ is purely imaginary as
$|\xi|=1$.\qed
\enddemo

From this moment the proof very much resembles that of Theorem 3.1. For $\xi$
bounded away from the point $1$ we use Lemma 6.4 to deform the
$\zeta$-integration contour to the one where the ratio
${F(\sigma,\zeta)}/{F(\tau,\zeta\xi^{-1})}$ remains bounded as $t\to 1$. The
result of Proposition 3.2(i) then implies that the integral over such $\xi$
tends to zero. On the other hand, for $\xi$ close to 1 we deform the
$\zeta$-contour to $\Gamma_\varkappa$ and use Proposition 3.2(i) together with
the asymptotics of ${F(\sigma,\zeta)}/{F(\tau,\zeta\xi^{-1})}$ to see that (cf.
the formula before Lemma 3.3)
$$
\multline K(\sigma,x;\tau,y)\sim \frac{1}{2(i\pi)^2}
\int_{\Gamma_\varkappa}\frac{\prod_{k\in\overline{1,\tau}:B[k]=1}(1-\zeta)}
{\prod_{k\in\overline{1,\sigma}:B[k]=1}(1-\zeta)}
\frac{\prod_{k\in\overline{\sigma+1,N}:A[k]=1}(1-1/\zeta)}
{\prod_{k\in\overline{\tau+1,N}:A[k]=1}(1-1/\zeta)}\\ \times \int_{-\infty\pm
i\varepsilon}^{+\infty\pm i\varepsilon} \frac{e^{\frac{iu}{\pi}\bigl(\ln
z-\gamma-l\ln(1-\zeta)-d\ln(1-\zeta^{-1})\bigr)}}{e^u-e^{-u}}\,\frac{dud\zeta}
{\zeta^{x-y+1}}\,,
\endmultline
$$
as $t\to 1$, $\varepsilon>0$ is small enough, the sign in $\pm i\varepsilon$ is
``$+$'' for $\sigma\le \tau$ and ``$-$'' for $\sigma>\tau$, and $\gamma$ is
replaced by $\gamma N$ in the case of Theorem 6.3. The integration variable $u$
is related to $\xi$ via $u=-{i\pi}r^{-1} \ln\xi$.

The integral over $u$ is evaluated by Lemma 3.3 to give
$$
\frac {\pm i\pi}{1+\bigl(z^{-1}e^\gamma(1-\zeta)^l(1-\zeta^{-1})^d\bigr)^{\pm
1}}\,.
$$
Finally, similarly to the proof of Theorem 3.1, the integration over $z$ (see
the formula of Proposition 5.2) is shown to be asymptotically equivalent to
substituting $z=1$.

In the case of Theorem 6.3 the argument is exactly the same except
for the fact that $e^\gamma(1-\zeta)^l(1-\zeta^{-1})^d$ converges
to either 0 or $+\infty$, see Comment 1 to Theorem 6.3. However,
Lemma 3.3 can still be applied because its estimates remain
uniform as long as (in the notation of Lemma 3.3) $\varepsilon\Re
a$ is uniformly bounded and $|\Im a|$ is bounded away from 1.
Since $\varepsilon$ can be made arbitrarily small, $\Re a$ is
allowed to converge to $\pm\infty$.

We skipped the discussion of several technical issues here. One
should make sure that the deformation of $\zeta$-contours happens
in such a way that the poles of the ratio
${F(\sigma,\zeta)}/{F(\tau,\zeta\xi^{-1})}$ are not passed in the
process; the remainder in the approximation of $\ln(x,t)_\infty$
by $\di(x)$ may grow as $x$ approaches cut $x\in(1,+\infty)$, and
one needs to control the growth of the derivative of the remainder
term in order to apply Lemma 3.3. All these issues can be
resolved. However, the arguments are rather tedious although
fairly straightforward, and we omit them.\qed
\enddemo

\head 7. Bulk of large cylindric partitions with period of intermediate growth
\endhead

In the previous section we have seen that if the period $N$ is small comparing
to $|\ln s|^{-1}$ then the global density function for the uniform measure on
cylindric partitions with a fixed profile, when the partitions are scaled by
$|\ln s|$, converges to a rotationally invariant limit which depends only on
the average slope of the profile.

In this section we consider the case when $t=s^N$ has a nontrivial limit
strictly between $0$ and $1$. Equivalently, $N\ln s^{-1} \to \ln t^{-1}\in
(0,+\infty)$. This situation is much more complex because, as we will see, the
global density is no longer rotationally invariant.

It is worth noting that if $N$ grows fast enough (so that $t\to 0$) then the
random cylindric partitions split into disjoint independent random plane
partitions located at the ``corners'' of the profile provided that these
corners are deep enough.

We will consider in detail two cases: the most rotationally symmetric one with
$\{A[k]\}$ and $\{B[k]\}$ being periodic of finite period (in the notation of
\S5 this means that $\mu$ is the staircase-like partition), and the most
rotationally asymmetric one with $\{A[k]\}$ and $B\{[k]\}$ consisting of two
blocks of zeroes and ones (this means that $\mu$ is either empty or it is the
rectangle $d\times l$).

In both cases we will assume that the average slope $\varkappa$ is equal to 1.
This will allow us to obtain asymptotic results without substantially deforming
the contours in our integral representation for the correlation kernel. We hope
to consider more general cases in a subsequent publication.

\proclaim{Proposition 7.1} Consider the uniform measure on
cylindric partitions with a given profile $\{A[k],B[k]\}_{k=1}^N$
and assume that the sequences $\{A[k]\}$ and $\{B[k]\}$ are
periodic with an even period $M$ and $A(0,M]=B(0,M]=M/2$. Then as
$$
s\to 1-, \quad N\to\infty, \quad s^N\to t\in (0,1),
$$
the limit of the correlation functions of this measure is described by Theorem
6.3 above with $\varkappa=1$ (see also Comment 3 after Theorem 6.3).
\endproclaim
\demo{Proof} For the sake of convenience let us consider the case when the
large period $N$ is a multiple of the small period $M$; $N/M\in\Z$. This
assumption is by no means necessary and can be easily removed.

Also, using the rotational invariance, we may assume without loss of generality
that all the time moments $\tau_1(s),\dots,\tau_n(s)$ from the statement of
Theorem 6.3 remain positive and uniformly bounded as $s\to 1$.

As in the proof of Theorems 6.2 and 6.3, we start with the integral
representation for the correlation kernel $K(\sigma,x;\tau,y)$ given in
Proposition 5.2. We have
$$
\multline \dfrac{F(\sigma,\zeta)}{F(\tau,\eta^{-1})}=
\frac{\prod_{k\in\overline{1,\tau}:B[k]=1}(1-s^{-k}/\eta)}
{\prod_{k\in\overline{1,\sigma}:B[k]=1}(1-s^{-k}\zeta)}
\frac{\prod_{k\in\overline{1,\tau}:A[k]=1}(1-s^k\eta)}
{\prod_{k\in\overline{1,\sigma}:A[k]=1}(1-s^k/\zeta)}
\\ \times \prod_{k\in\overline{1,N}:B[k]=1}\frac{(s^{N-k}/\eta;t)_\infty}
{(s^{N-k}\zeta;t)_\infty}\prod_{k\in\overline{1,N}:A[k]=1} \frac
{(s^{k}/\zeta;t)_\infty}{(s^{k}\eta;t)_\infty}\,.
\endmultline
$$
Using the periodicity assumption, we can rewrite the second factor as
$$
\prod_{k\in\overline{1,M}:B[k]=1}\frac{(s^{M-k}/\eta;s^M)_\infty}
{(s^{M-k}\zeta;s^M)_\infty}\prod_{k\in\overline{1,M}:A[k]=1} \frac
{(s^{k}/\zeta;s^M)_\infty}{(s^{k}\eta;s^M)_\infty}
$$
(here the base $t$ of all products $(x;t)_\infty$ changed from $t$ to $s^M$).

Starting from this moment the proof almost literally repeats the
proof of Theorem 2 in \cite{OR1}, see \S\S3.1.3-3.1.9.
Approximating $(x;s^M)$ by the dilogarithm function (see the
previous section) and keeping in mind that $A(0,M]=B(0,M]=M/2$, we
obtain that as $s\to 1$ the second factor in the formula for the
ratio $F(\sigma,\zeta)/F(\tau,\eta)$ above is approximated by
$$
\exp\left(\frac {M(\di(\zeta)-\di(1/\zeta)+\di(\eta)-\di(1/\eta))}{2(1-s^M)}
\right).
$$
Note that $1-s^M$ can be replaced by $M|\ln s|$ as $s\to 1$.

The only other two factors in the integral representation of
$K(\sigma,x;\tau,y)$ which may not have a finite limit as $s\to 1$ are
$\zeta^{-x-\frac 12}\eta^{- y-\frac 12}$. Since $x,y\sim\gamma/|\ln s|$, we can
include these factors into the exponential above to obtain
$$
\exp\left(\frac{\di(\zeta)-\di(1/\zeta)-2\gamma\ln\zeta+\di(\eta)-\di(1/\eta)-
2\gamma\ln\eta}{2|\ln s|} \right).
$$

Following \S3.1 of \cite{OR1} we see that the function
$$
\Re\bigl(\di(\zeta)-\di(1/\zeta)-2\gamma\ln\zeta\bigr),
$$
which clearly vanishes on the unit circle $|\zeta|=1$, is negative
slightly inside the unit circle if $-\frac 12\arcsin(\frac
{e^{-\gamma}}2)\le\arg\zeta\le\frac 12\arcsin(\frac
{e^{-\gamma}}2)$ and slightly outside the unit circle if
$|\arg\zeta|>\frac 12\arcsin(\frac {e^{-\gamma}}2)$; the points
$\arg\zeta=\pm \frac 12\arcsin(\frac {e^{-\gamma}}2)$ with
$|\zeta|=1$ are the critical points of the function
$\di(1/\zeta)-\di(\zeta)-2\gamma\ln\zeta$.

We would like to deform both $\zeta$- and $\eta$-integration contours to put
them in the domain where the exponential above converges to zero, i.e. where
$$
\Re\bigl(\di(\zeta)-\di(1/\zeta)-2\gamma\ln\zeta\bigr)<0,\qquad
\Re\bigl(\di(\eta)-\di(1/\eta)-2\gamma\ln\eta\bigr)<0
$$
(this corresponds to using the contour $\gamma_<$ from
\cite{OR1}). However, there is an obstacle: the function
$\theta_3(-\zeta\eta\, t^{-\frac 12};t)$, which enters the
denominator of the integrand in the formula for
$K(\sigma,x;\tau,y)$, vanishes when $\zeta\eta=1$. Thus, while
deforming the contours we have to add the residues corresponding
to the poles $\zeta\eta=1$.

Consider the case $\sigma\le \tau$ first. Then Proposition 5.2 states that the
$\zeta$- and $\eta$-contours are such that $\eta$-contour contains the
$1/\zeta$-contour. Since our desired contours go inside the unit circle when
the arguments of $\zeta$ and $\eta$  are inside the interval $(-\frac
12\arcsin(\frac {e^{-\gamma}}2),\frac 12\arcsin(\frac {e^{-\gamma}}2))$, this
is where the residues have to be added. For $\gamma<-\ln 2$ the residues have
to be taken along the whole integration contour.

Observe that
$$
\multline
 -\operatorname{Res}_{\eta=1/\zeta} \frac{\prod_{n\ge
1}(1-t^n)^3\,\theta_3(z\zeta\eta;t)}{\theta_3(z;t)\,
\theta_3(-\zeta\eta\,t^{-\frac
12};t)}\dfrac{F(\sigma,\zeta)}{F(\tau,\eta)}\,\frac 1{\zeta^{x+\frac
12}\eta^{y+\frac 12}}\\= \frac
{\prod_{k\in\overline{\sigma+1,\tau}:B[k]=1}(1-s^{-k}\zeta)
\prod_{k\in\overline{\sigma+1,\tau}:A[k]=1}(1-s^k/\zeta)} {\zeta^{x-y+1}}\,,
\endmultline
$$
which asymptotically equals
$(1-\zeta)^{B(\sigma,\tau]}(1-\zeta^{-1})^{A(\sigma,\tau]}/\zeta^{x-y+1}$ as
$s\to 1$. (Recall that we assumed in the beginning of the proof that the time
moments $\sigma,\tau$ remain bounded as $s\to 1$.) Thus, the residue
contribution is asymptotically equal to
$$
\frac 1{2\pi i}\int_{e^{-\frac i2\arcsin(\frac {e^{-\gamma}}2)}}^{e^{\frac
i2\arcsin(\frac
{e^{-\gamma}}2)}}(1-\zeta)^{B(\sigma,\tau]}(1-\zeta^{-1})^{A(\sigma,\tau]}
\,\frac{d\zeta}{\zeta^{x-y+1}}\,.
$$
On the other hand, the integral when $\zeta$ and $\eta$ range over
the desired contour converges to zero as $s\to 1$, see \S 3.1.1 of
\cite{OR1} for additional explanations.

The arguments for $\sigma>\tau$ are very similar.

Note that all the estimates are uniform in the parameter $z$,
$|z|=1$, and the final asymptotic expressions are independent of
$z$. This means that the integration over $z$ in the formula of
Proposition 5.2 in the limit $s\to 1$ can be simply removed. \qed
\enddemo

Now let us consider the non-symmetric case. Take $d=l=N/2$ and
$$
\{A[k]\}_{k=1}^N=(\underbrace{1,1,\dots,1}_{N/2},\underbrace{0,\dots,0,0}_{N/2}),\qquad
\{B[k]\}_{k=1}^N=(\underbrace{0,0,\dots,0}_{N/2},\underbrace{1,\dots,1,1}_{N/2})
$$
(we are assuming that $N$ is even).  This corresponds to the square partition
$\mu=(\frac N2)^{\frac N2}$, and the random sequences of ordinary partitions of
the form
$$
\lambda^{(1)}\succ\lambda^{(2)}\succ \dots\succ\lambda^{(\frac N2+1)}\prec
\lambda^{(\frac N2+2)}\prec\dots\prec\lambda^{(N)}\prec
\lambda^{(N+1)}=\lambda^{(1)}.
$$
We will describe the limit of the correlation functions near the largest
partition $\lambda^{(1)}$ and near the smallest partition $\lambda^{(\frac
N2+1)}$. (Recall that the time variable of the periodic Schur process or,
equivalently, the index of the Young diagrams above, may be viewed as an
element of $\Z/N\Z$.)

\proclaim{Theorem 7.2} In the limit
$$
s\to 1-,\quad N\to\infty,\quad s^N\to t\in(0,1)
$$
the correlation functions of the uniform measure on cylindric partitions with
the profile as above have limits in the following sense:
 Choose
$x_1(s),\dots,x_n(s)\in\Z+\frac 12$ such that as $s\to 1-$, $|\ln s|\cdot
x_k(s)\to \gamma$ for all $k=1,\dots,n$ and some $\gamma\in\R$, and all
pairwise distances $x_i-x_j=x_i(s)-x_j(s)$ remain constant.

(i) Choose the time moments $\tau_1(s),\dots,\tau_n(s)\in \{-\frac
N2,\dots,\frac N2\}$ such that as $s\to 1$ the absolute values $|\tau_i(s)|$
remain uniformly bounded. Then
$$
\rho_n(\tau_1(s),x_1(s);\dots,\tau_n(s),x_n(s))=\det\bigl[
 K^{(\gamma)}_{\tau_i(s),\tau_j(s)}(x_i-x_j)\bigr]_{i,j=1}^n +o(1),
$$
where the correlation kernel has the following form:
$$
\K^{(\gamma)}_{\tau,\tau}(x-y)=\delta(x-y)\ \text{
 if  }\
 \gamma\le\gamma_1^{\min}(t):=2\ln\frac{(-\sqrt{t};t)_\infty}{(-1;t)_\infty}\,,
$$
and for $\gamma>\gamma_1^{\min}(t)$
$$
\K^{(\gamma)}_{\sigma,\tau}(x-y)=\cases\frac 1{2\pi
i}\int\limits_{e^{-ic}}^{{e^{ic}}}
(1-\zeta)^{B(\sigma,\tau]}(1-\zeta^{-1})^{A(\sigma,\tau]}
\,\dfrac{d\zeta}{\zeta^{x-y+1}}\,,&\sigma\le
\tau,\\
-\frac 1{2\pi i}\int\limits_{e^{-ic}}^{e^{ic}}
(1-\zeta)^{-B(\tau,\sigma]}{(1-\zeta^{-1})^{-A(\tau,\sigma]} }
\,\dfrac{d\zeta}{\zeta^{x-y+1}}\,, &\sigma>\tau,
\endcases
$$
where the constant $c=c(\gamma)\in(0,\pi)$ from the limits of integration is
uniquely determined by the relation
$$
\gamma=2\ln\left|\frac{(e^{ic}\sqrt{t};t)_\infty}{(e^{ic};t)_\infty}\right|.
$$

(ii) Choose the time moments $\tau_1(s),\dots,\tau_n(s)\in \{1,\dots,N\}$ such
that as $s\to 1$ the absolute values $|\tau_i(s)-\frac N2|$ remain uniformly
bounded. Then the asymptotics of the correlation functions has the
determinantal form as in (i) above with the correlation kernel given by
$$
K_{\tau,\tau}^{(\gamma)}(x-y)= \cases \delta(x-y),&\
\text{if}\quad \gamma\le
\gamma_2^{\min}(t):=2\ln\dfrac{(-t;t)_\infty}{(-\sqrt{t};t)_\infty}\,,
\\0,&\ \text{if}\quad
\gamma\ge\gamma_2^{\max}(t):=2\ln\dfrac{(t;t)_\infty}{(\sqrt{t};t)_\infty}\,,
\endcases
$$
and for $\gamma\in(\gamma_2^{\min}(t),\gamma_2^{\max}(t))$ it is given by the
formula of (i) above, where the constant $c=c(\gamma)\in(0,\pi)$ from the
limits of integration is uniquely determined by the relation
$$
\gamma=2\ln\left|\frac{(e^{ic}t;t)_\infty}{(e^{ic}\sqrt{t};t)_\infty}\right|.
$$
\endproclaim
\demo{Comments} 1. For any $t\in (0,1)$ the functions
$$
\gamma_1(t, c)=2\ln
\left|\frac{(e^{ic}\sqrt{t};t)_\infty}{(e^{ic};t)_\infty}\right|,\qquad
\gamma_2(t,c)=2\ln\left|\frac{(e^{ic}t;t)_\infty}{(e^{ic}\sqrt{t};t)_\infty}\right|
$$
as functions of $c$ are strictly decreasing on $(0,\pi)$, and
$$
\gathered \lim_{c\to 0}\gamma_1(t,c)=+\infty, \qquad \lim_{c\to
\pi}\gamma_1(t,c)=\gamma_1^{\min}(t),\\
\lim_{c\to 0}\gamma_2(t,c)=\gamma_2^{\max}(t), \qquad \lim_{c\to
\pi}\gamma_2(t,c)=\gamma_2^{\min}(t).
\endgathered
$$
Indeed, this follows from the fact that for any $\alpha,\beta\in(0,1)$ and
$c\in (0,\pi)$, we have
$$
\frac d{dc}\left|\frac{1-\alpha\beta e^{ic}}{1-\alpha
e^{ic}}\right|^2=-\frac{2\alpha(1-\beta)(1-\alpha^2\beta)\sin(c)}{|1-\alpha
e^{ic}|^4}<0.
$$

The decay of $\gamma_1(t,c)$ and $\gamma_2(t,c)$ guarantees the existence of
the unique $c\in(0,\pi)$ satisfying the needed relation.

2. The global density function in both (i) and (ii) is given by
$$
\rho(\gamma)=\cases 1, &\gamma\le\gamma_{1,2}^{\min}(t),\\
\frac1\pi \,c(\gamma),& \gamma_{1,2}^{\min}(t)\le \gamma\le
\gamma_{1,2}^{\max}(t),
\\0,& \gamma\ge\gamma_{1,2}^{\max}(t),
\endcases
$$
where $\gamma_1^{\max}(t)=+\infty$. In particular, in (i) we see that the
random Young diagrams have only the lower edge, while in (ii) there are both
the lower and the upper edges.

3. The functions $\gamma_{1,2}(t,c)$ can be written in terms of the Jacobi
elliptic sine function:
$$
\gathered \left|\frac{(e^{ic}\sqrt{t};t)_\infty}{(e^{ic};t)_\infty}\right|^2=
\frac{t^{\frac 18}}{2\sqrt{k}}\,\frac 1{\operatorname{sn}\bigl(\frac
{\operatorname{K}c}\pi ,k\bigr)\sin\bigl(\frac c2\bigr)}\,, \qquad
\left|\frac{(e^{ic}{t};t)_\infty}{(e^{ic}\sqrt{t};t)_\infty}\right|^2=
\frac{\sqrt{k}} {2\,t^{\frac 18}}\,\frac {\operatorname{sn}\bigl(\frac
{\operatorname{K}c}\pi ,k\bigr)}{\sin\bigl(\frac
c2\bigr)}\,,\\
k=\frac{\theta_2^2(0;t)}{\theta_3^2(0;t)},\qquad \operatorname{K}=\tfrac \pi 2
\,\theta_3^2(0;t).
\endgathered
$$

Note that as $t\to 1$, the elliptic sine tends to 1, and the formulas of both
(i) and (ii) degenerate to those of Comment 3 after Theorem 6.3. This agrees
with the fact that if $N$ is growing slowly then the global limit density
function is rotationally invariant.

On the other hand, as $t\to 0$ in the setting of (i) we have
$$
e^\gamma=\frac{|(e^{ic(\gamma)}\sqrt{t};t)_\infty|^2}
{|(e^{ic(\gamma)};t)_\infty|^2}\to \frac 1{|1-e^{ic(\gamma)}|^2}=\frac
{1}{2(1-\cos c(\gamma))}\ \text{  or  }\
c(\gamma)\to\arccos\bigl(1-\tfrac{e^{-\gamma}}2\bigr)
$$
which is exactly the behavior of the largest section of the random
plane partitions, cf. \S3.1.10 of \cite{OR1} with $\tau=0$. This
agrees with the fact that if $N$ grows too fast then the random
cylindric partitions split into disjoint independent plane
partitions located in deep corners of the profile (in this case we
have only one such corner).
\enddemo

\demo{Proof} The proof is very similar to that of Proposition 7.1
and \cite{OR1, Theorem 2}. The only difference from the arguments
used for Proposition 7.1 is in the specifics of the functions
$F(\tau,\zeta)$. Using Remark 5.4, for bounded $\tau>0$ we have
$$
F(\tau,\zeta)=\prod\limits_{k=1}^\tau\bigl(1- {s^{k}}/\zeta\bigr)^{-1}
\frac{\prod\limits_{k=1}^{N/2}\bigl(s^{k}/\zeta;t\bigr)_\infty}
{\prod\limits_{k=0}^{N/2-1}\bigl(s^{k}\zeta;t\bigr)_\infty}\,.
$$
The first factor has the finite limit $(1-\zeta^{-1})^{-\tau}$, while for the
second factor we have
$$
\multline
\prod\limits_{k=0}^{N/2-1}\bigl(s^{k}x;t\bigr)_\infty=\frac{(x;s)_\infty
(xt;s)_\infty(xt^2;s)_\infty\cdots}
{(x\sqrt{t};s)_\infty(xt\sqrt{t};s)_\infty(xt^2\sqrt{t};s)_\infty\cdots}\\ \sim
\exp\frac 1{|\ln s|}\sum_{m=0}^\infty(\di(xt^m\sqrt{t})-\di(xt^m))\,.
\endmultline
$$
Thus, in $F(\tau,\zeta)$ we have the factor
$$
\exp\frac 1{|\ln s|}\sum_{m=0}^\infty(\di(\zeta t^m)-\di(\zeta
t^m\sqrt{t})-\di(t^m/\zeta)+\di(t^m\sqrt{t}/\zeta))
$$
times the part which has a finite limit as $s\to 1$.

The real part of the above sum is identically equal to zero on the
unit circle $|\zeta|=1$. Furthermore, the derivative of this sum
with respect to $\zeta$ multiplied by $\zeta$ equals (for
$|\zeta|=1$)
$$
-2\ln \prod_{m\ge 0}\left|\frac{1-\zeta t^m}{1-\zeta t^m\sqrt{t}}\right|=2\ln
\left|\frac{(\zeta\sqrt{t};t)_\infty}{(\zeta;t)_\infty}\right|.
$$
From Comment 1 above we know that for any $\gamma\in(\gamma_1^{\min},+\infty)$
this expression is strictly greater than $\gamma$ on the arch $\arg\zeta\in
(-c(\gamma),c(\gamma))$, where $c=c(\gamma)$ satisfies the relation in the
hypothesis (i) above, and it is strictly less than $\gamma$ on the
complementary arch. This means that if we deform the $\zeta$-integration
contour so that it goes slightly inside the unit circle for $\arg\zeta\in
(-c(\gamma),c(\gamma))$, and slightly outside the unit circle for
$|\arg\zeta|>c(\gamma)$, then the factor $F(\tau,\zeta)\zeta^{-x-\frac 12}$
with $x\sim \gamma/|\ln s|$ will exponentially decay on this contour.

Similarly, for $\tau=N-\widehat\tau$ with nonnegative and bounded
$\widehat\tau$ we have
$$
F(\tau,\zeta)=\prod_{k=0}^{\widehat
\tau}\bigl(1-s^k\zeta/t\bigr)\,\frac{\prod\limits_{k=1}^{N/2}\bigl(
ts^{k}/\zeta;t\bigr)_\infty }{ \prod\limits_{k=0}^{N/2-1}\bigl(
s^{k}\zeta/t;t\bigr)_\infty }\,
$$
which has the same asymptotically nontrivial part as before up to the change
$\zeta\to t\zeta$. This means that in this case the desired $\zeta$-integration
contour is obtained from the one in the previous case by multiplication by $t$.

The proof of (i) proceeds in the same way as the proof of Proposition 7.1: The
deformation of contours to the desired ones meets an obstacle of the theta
function $\theta_3(-\zeta\eta\, t^{-\frac 12};t)$ in the denominator of the
integrand vanishing for $\zeta\eta=1$ or $\zeta\eta=t$; evaluating the
corresponding residue yields the limit correlation kernel, while the remaining
integrals over the constructed contours tend to zero as $s\to 1$. The estimates
are uniform in $z$, $|z|=1$, and the limiting kernel is independent of $z$;
thus, the integration over $z$ in the formula for the correlation functions in
Proposition 5.2 can be removed.

In order to prove (ii) we need to consider $F(\tau,\zeta)$ where
$\tau=N/2+\widetilde\tau$ with bounded $\widetilde\tau$. The asymptotically
nontrivial part is independent of the sign of $\widetilde\tau$, so let us take
$\widetilde \tau>0$ to be concrete. We have
$$
F(\tau,\zeta)=\prod\limits_{k=1}^{\widetilde
\tau}\bigl(1-\sqrt{t}s^k\zeta\bigr)\,\frac{\prod\limits_{k=1}^{N/2}\bigl(
ts^{k}/\zeta;t\bigr)_\infty}{
\prod\limits_{k=0}^{N/2-1}\bigl(s^{k}\zeta;t\bigr)_\infty}\,.
$$
The first factor has a finite limit as $s\to 1$, while the second factor (the
ratio) produces
$$
\exp\frac 1{|\ln s|}\sum_{m=0}^\infty(\di(\zeta t^m)-\di(\zeta
t^m\sqrt{t})-\di(t^{m+1}/\zeta)+\di(t^{m+1}\sqrt{t}/\zeta)).
$$
The real part of this sum vanishes on the circle $|\zeta|=\sqrt{t}$. Further,
on this circle the derivative of the sum multiplied by $\zeta$ equals
$$
2\ln \left|\frac{(\zeta_0 t;t)_\infty}{(\zeta_0\sqrt{t};t)_\infty}
\right|,\qquad \zeta_0=\frac{\zeta}{|\zeta|}\,.
$$
The remaining part of the proof is just as in the case (i) considered above. It
is worth noting that for $\gamma>\gamma_2^{\max}$ the deformation of the
contours to the domain of exponential decay of the integrand does not meet any
obstacles, and thus the limiting kernel is identically equal to zero.\qed
\enddemo

\head 8. On a measure of Nekrasov and Okounkov
\endhead

A large part of the material presented in this section is the result of joint
discussions of the author and Grigori Olshanski.

Our goal in this final section is to demonstrate how the
techniques developed in the first three sections apply to a
remarkable measure on partitions introduced by Nekrasov and
Okounkov in \cite{NO}.

Fix $\mu\in\C$ and $t\in (0,1)$. The object of interest is a (generally
speaking, complex) measure on the set of all partitions given by the formula,
see \cite{NO, \S6.2},
$$
M_{\mu,t}(\lambda)=\prod_{n\ge 1}(1-t^n)^{1-\mu^2}\cdot
t^{|\lambda|}\prod_{\square\in\lambda} \frac
{h(\square)^2-\mu^2}{h(\square)^2}\,,\qquad \lambda\in\Y.
$$
Here the product is taken over all boxes of the Young diagram $\lambda$, and
$h(\square)$ denotes the length of the hook rooted at the given box. The sum of
the weights $M_{\mu,t}(\lambda)$ over all partitions $\lambda$ is identically
equal to 1, see \cite{NO, (6.12)}. Note that this measure becomes a {\it
probability\/} measure (meaning that all weights are nonnegative) if $\mu\in
i\R$.

One interesting feature of measures $M_{\mu,t}$ is that they
interpolate between the uniform measure on partitions arising at
$\mu=0$, and the poissonized Plancherel measure which is obtained
by the limit transition
$$
\lim_{\Sb \mu\to\infty,\,t\to 0\\
\mu^2t\to -\theta\endSb}
M_{\mu,t}(\lambda)=e^{-\theta}\theta^{|\lambda|}\prod_{\square\in\lambda}\frac
1{h(\square)^2}=e^{-\theta}\left(\frac
{\dim\lambda\,\theta^{\frac{|\lambda|}2}}{|\lambda|!}\right)^2.
$$
Here $\dim\lambda$ is the dimension of the irreducible
representation of the symmetric group $S_{|\lambda|}$. We refer
the reader to \cite{BOO} for details and further references on the
Plancherel measures.

Let us denote by $\rho_\mu$ the specialization of the algebra of symmetric
functions $\Lambda$ such that
$$
h_n(\rho_\mu)=\frac{(\mu)_n}{n!},\quad n\ge 0,\quad \text{  or  }\quad
H(\rho_\mu,u)=(1-u)^{-\mu}.
$$
Here $(a)_k=a(a+1)=\dots (a+k-1)$ denotes the Pochhammer symbol.

The applicability of the periodic Schur process results follows
from

\proclaim{Lemma 8.1} For any $\kappa,\lambda\in\Y$ choose
$r\ge\max\{\ell(\kappa),\ell(\lambda)\}$ and set
$$
k_i=\kappa_i-i+r,\quad l_i=\lambda_i-i+r,\qquad i=1,\dots,r.
$$
Then
$$
\multline \sum_{\nu\in\Y}s_{\kappa/\nu}(\rho_\mu)s_{\lambda/\nu}(\rho_{-\mu})=
\prod_{\square\in\kappa}\frac{h(\square)+\mu}{h(\square)}
\prod_{\square\in\lambda}\frac{h(\square)-\mu}{h(\square)}
\\ \times \frac{(-1)^{\frac{r(r-1)}2}\mu^r\prod_{1\le i<j\le
r}(k_i-k_j+\mu)(l_i-l_j-\mu)}{\prod_{i,j=1}^r(\mu+k_i-l_j)}\,.
\endmultline
$$
In particular, for $\kappa=\lambda$ the last factor turns into 1 and we obtain
$$
\sum_{\nu\in\Y}s_{\lambda/\nu}(\rho_\mu)s_{\lambda/\nu}(\rho_{-\mu})=
\prod_{\square\in\lambda}\frac{h(\square)^2-\mu^2}{h(\square)^2}\,.
$$
\endproclaim
\demo{Comments} 1. Only finitely many terms in the sums over $\nu$ above are
nonzero.

2. The formula for $\kappa=\lambda$ can be easily extracted from the Fock space
representation of the measure $M_{\mu,t}$ given in \cite{NO}.
\enddemo
\demo{Proof} Since the terms of the sum vanish unless
$\nu\subset\kappa$ and $\nu\subset\lambda$, we can restrict the
sum to $\nu$ with $\ell(\nu)\le r$. Denote $n_i=\nu_i-i$,
$i=1,\dots,r$.

The Jacobi-Trudi formula for the skew Schur functions, see
\cite{Macd, \S{I} (5.4)}, gives
$$
s_{\kappa/\nu}=\det[h_{k_i-n_j}]_{i,j=1}^r,\qquad
s_{\lambda/\nu}=\det[h_{l_i-n_j}]_{i,j=1}^r.
$$
Applying the Cauchy-Binet summation formula, we obtain
$$
\multline
\sum_{\nu\in\Y}s_{\kappa/\nu}(\rho_\mu)s_{\lambda/\nu}(\rho_{-\mu})=
\sum_{n_1>\dots>n_r\ge 0}\det[h_{k_i-n_j}(\rho_\mu)]_{i,j=1}^r
\det[h_{l_i-n_j}(\rho_{-\mu})]_{i,j=1}^r\\
=\det\left[\sum_{m=0}^{\min\{k_i,l_j\}}h_{k_i-m}(\rho_\mu)h_{l_j-m}
(\rho_{-\mu})\right]_{i,j=1}^r=
\det\left[\sum_{m=0}^{\min\{k_i,l_j\}}\frac{(\mu)_{k_i-m}}{(k_i-m)!}
\frac{(\mu)_{l_j-m}}{(l_j-m)!}\right]_{i,j=1}^r
\endmultline
$$
Using the identity
$$
\sum_{m=0}^{\min\{k,l\}}\frac{(\mu)_{k-m}}{(k-m)!}
\frac{(\mu)_{l-m}}{(l-m)!}=
\frac{\mu(\mu+1)_k(-\mu+1)_l}{k!\,l!\,(\mu+k-l)}\,,
$$
which can be proved by induction on $\min\{k,l\}$ with $k-l$
fixed, and the formula for the Cauchy determinant
$$
\det \left[\frac
1{\mu+k_i-l_j}\right]_{i,j=1}^r=(-1)^{\frac{r(r-1)}2}\,\frac{\prod_{1\le i<j\le
r} (k_i-k_j)(l_i-l_j)}{\prod_{i,j=1}^r(\mu+k_i-l_j)}\,,
$$
we obtain
$$
\multline \sum_{\nu\in\Y}s_{\kappa/\nu}(\rho_\mu)s_{\lambda/\nu}(\rho_{-\mu})\\
= (-1)^{\frac{r(r-1)}2}\mu^r\prod_{i=1}^r\frac{(\mu+1)_{k_i}(-\mu+1)_{l_i}}
{k_i!\,l_i!}\,\frac{\prod_{1\le i<j\le r}
(k_i-k_j)(l_i-l_j)}{\prod_{i,j=1}^r(\mu+k_i-l_j)}\,.
\endmultline
$$
Finally, we use the fact that the set of hook lengths
$\{h(\square)\}_{\square\in\kappa}$ can be obtained as the union
of sets $\{1,\dots,k_i\}$ for $i=1,\dots,r$ minus the set of
numbers $\{k_i-k_j\}_{1\le i<j\le r}$, and similarly for
$\lambda$. This brings us to the desired formula. \qed
\enddemo

Let us now associate to the measure $M_{\mu,t}$ a measure on point
configurations (subsets) in $\Z'=\Z+\frac 12$ via the map
$$
\Y\longrightarrow 2^{\Z'},\qquad \lambda\mapsto \bigl\{\lambda_i-i+\tfrac
12\bigr\}_{i\ge 1}.
$$
The correlation functions of this measure on point configurations are given by
$$
\rho_n(x_1,\dots,x_n)=\sum_{\lambda\in\Y:\{\lambda_i-i+\frac 12\}_{i\ge
1}\supset\{x_1,\dots,x_n\}} M_{\mu,t}(\lambda).
$$
We will also need the {\it shift-mixed\/} version of $M_{\mu,t}$. For any
$z\in\C\setminus\{-t^{\pm\frac 12},-t^{\pm\frac 32},\dots\}$ define the measure
$M_{\mu,t,z}$ on $\Y\times\Z$ as the product-measure
$$
M_{\mu,t,z}(\lambda,S)=\frac{z^S
t^{\frac{S^2}2}}{\theta_3(z;t)}\,M_{\mu,t}(\lambda)\,,
$$
cf. \S2. We associate to this measure a measure on point
configurations in $\Z'$ via the map
$$
\Y\times \Z\longrightarrow 2^{\Z'},\qquad \lambda\mapsto
\bigl\{S+\lambda_i-i+\tfrac 12\bigr\}_{i\ge 1}.
$$
Its correlation functions are given by
$$
\rho_n^{\shift}(x_1,\dots,x_n)=\sum_{(\lambda,S)\in\Y\times\Z:\{S+\lambda_i-i+\frac
12\}_{i\ge 1}\supset\{x_1,\dots,x_n\}} M_{\mu,t,z}(\lambda,S).
$$
The two sets of correlation functions are related as described in Proposition
2.1.

\proclaim{Proposition 8.2} The correlation functions of the
shift-mixed measure $M_{\mu,t,z}$ have determinantal form: For any
$n=1,2,\dots$ and $x_1,\dots,x_n\in\Z'$,
$$
\rho_n^{\shift}(x_1,\dots,x_n)=\det[K_{\mu,t,z}(x_i,x_j)]_{i,j=1}^n,
$$
where the correlation kernel has the form
$$
\multline K_{\mu,t,z}(x,y)=-\frac{\prod_{n\ge 1}(1-t^n)^3}{\theta_3(z;t)
\,(2\pi i)^2}\\ \times \oint_\zeta\oint_\eta \prod_{m\ge 0}\dfrac{(1-\zeta
t^m)^\mu (1-t^{m+1}/\zeta)^\mu} {(1-t^m/\eta)^\mu (1-\eta t^{m+1})^\mu}
\,\frac{\theta_3(z\zeta\eta;t)}{\theta_3(-\zeta\eta\, t^{-\frac 12};t)}\,
\frac{d\zeta d\eta}{\zeta^{x+\frac12}\eta^{y+\frac12}}\,.
\endmultline
$$
Here both integration contours are simple positively oriented loops going
around the origin such that
$$
t<|\zeta|<1<|\eta|<t^{-1},\qquad 1<|\zeta\eta|<t^{-1},
$$
and we use
the principal branch of the logarithm to define $(\,\cdot\,)^\mu$.
\endproclaim

\demo{Proof} Consider the periodic Schur process with $N=1$ and specializations
$a=a[1]$ and $b=b[1]$ of $\Lambda$ defined by
$$
H(a;u)=(1-\alpha u)^{-\mu},\qquad H(b;u)=(1-\alpha u)^{\mu},\qquad 0<\alpha<1.
$$
Then by Proposition 1.1
$$
\sum_{\lambda\in\Y}t^{|\lambda|}\sum_{\nu\in\Y}s_{\lambda/\nu}(a)
s_{\lambda/\nu}(b)=\prod_{n\ge 1} \frac{(1-\alpha^2 t^n)^{\mu^2}}{1-t^n}\,,
$$
which converges to $\prod_{n\ge 1}(1-t^n)^{\mu^2-1}$ as $\alpha\to 1$. Since
$f(a)\to f(\rho_\mu)$ and $f(b)\to f(\rho_{\mu})$ for any $f\in\Lambda$ as
$\alpha\to 1$, we see that the weights of the periodic Schur process above as
well as those of its shift-mixed version converge to the corresponding weights
of $M_{\mu,t}$ and $M_{\mu,t,z}$. Then the fact that the sum of the weights in
all cases is identically equal to one, implies the convergence of the
correlation functions of the (shift-mixed) Schur process to those of
$M_{\mu,t}$ and $M_{\mu,t,z}$.

On the other hand, the analytic version of Theorem 2.2 (see also Remarks 2.3
and 2.4) gives the formula for the correlation functions of the shift-mixed
Schur process introduced above. Taking the limit $\alpha\to 1$ in the integral
representation of the kernel given in Remark 2.4 completes the proof. \qed
\enddemo

Let us now restrict our attention to the case when the measures $M_{\mu,t}$ and
$M_{\mu,t,z}$ become probability measures (in other words, all weights are
nonnegative). This happens when $\mu=i\mu_0$ with $\mu_0\in\R$, and $z\in\R_+$.

\proclaim{Theorem 8.3} Assume that $\mu=i\mu_0$ with $\mu_0\in\R$,
and $z>0$. Then as $t\to 1$ the correlation functions of the
shift-mixed measure $M_{\mu,t,z}$ have the following limit: Choose
$x_1(t),\dots,x_n(t)\in\Z'$ such that as $t\to 1$, $|\ln t|\cdot
x_k(t)\to \gamma$ for all $k=1,\dots,n$ and some $\gamma\in\R$,
and all pairwise distances $x_i-x_j=x_i(t)-x_j(t)$ are independent
of $t$. Then
$$
\lim_{t\to 1}\rho_n^{\shift}(x_1(t),\dots,x_n(t))=\det\bigl[
 \K^{(z,\gamma,\mu)}(x_i-x_j)\bigr]_{i,j=1}^n,
$$
where the correlation kernel has the following form
$$
\K^{(z,\gamma,\mu)}(d)=\frac 1{2\pi i}\oint\limits_{|\zeta=1|} \frac{1}
{1+z^{-1}e^{\gamma}(1-\zeta)^{-\mu}(1-\zeta^{-1})^{\mu}}\,
\frac{d\zeta}{\zeta^{d+1}},\qquad d\in\Z.
$$
Under the same assumptions the correlation functions of the measure $M_{\mu,t}$
converge to the limiting expression above evaluated at $z=1$.
\endproclaim
\demo{Comments} 1. For $|\zeta|=1$ we have
$$
(1-\zeta)^{-i\mu_0}(1-\zeta^{-1})^{i\mu_0}=\cases e^{\mu_0(\arg\zeta-\pi)},&
0<\arg\zeta\le\pi,\\
e^{\mu_0(\arg\zeta+\pi)},&-\pi\le\arg\zeta<0.
\endcases
$$
The change of sign of $\mu_0$ is equivalent to the change of variable
$\zeta\to\zeta^{-1}$ in the integral above, which in its turn is equivalent to
transposing the correlation kernel. Clearly, this operation does not change the
correlation functions.

2. The global limit density function for the shift-mixed case is equal to
$$
\multline \rho(\gamma)=\K^{(z,\gamma,\mu)}(0)=\frac 1{2\pi}\left(\int_0^\pi
\frac{d\phi} {1+z^{-1}e^{\gamma+\mu_0(\phi-\pi)}}+\int_{-\pi}^0 \frac{d\phi}
{1+z^{-1}e^{\gamma+\mu_0(\phi+\pi)}}\right)\\=\frac1{2\pi\mu_0}
\ln\frac{e^\gamma+ze^{\pi\mu_0}}{e^\gamma+ze^{-\pi\mu_0}}\,,
\endmultline
$$
and one has to substitute $z=1$ to get the formula corresponding
to the non-mixed measure $M_{\mu,t}$.

3. In the limit $\mu_0\to +\infty$, $\gamma\to\infty$ so that
$\gamma/\mu_0\to\widehat\gamma$, the correlation kernel $\K^{(z,\gamma,\mu)}$
becomes equivalent to the discrete sine kernel, cf. Example 3.4:
$$
\lim_{\Sb\mu_0\to+\infty, \gamma\to\infty,\\
\gamma/\mu_0\to\widehat\gamma\endSb} \K^{(z,\gamma,\mu)}(x-y) =\cases
0,&\widehat\gamma\ge\pi,\\e^{\frac i2{(\pi-\widehat\gamma)(x-y)}}\cdot
\dfrac{\sin\bigl(\frac{\pi-\widehat\gamma}{2}\,(x-y)\bigr)}{\pi(x-y)} ,&-\pi<\widehat\gamma<\pi,\\
\delta(x-y),&\gamma\le-\pi.
\endcases
$$
\enddemo

The proof of Theorem 8.3 is completely analogous to that of
Theorem 3.1 and we omit it.

\end